\documentclass[a4paper,10pt,oneside]{article}

\usepackage[utf8]{inputenc}
\usepackage{amsmath}
\usepackage{bm}
\usepackage{overpic}
\usepackage{multirow}
\usepackage{amsfonts}
\usepackage{amssymb}
\usepackage{amsthm}
\usepackage[english]{babel}
\usepackage{fontenc}
\usepackage{enumerate}
\usepackage{color}
\usepackage{dirtytalk}
\usepackage{subcaption}

\usepackage{geometry}
\newgeometry{vmargin={15mm}, hmargin={12mm,17mm}}

\usepackage{listings}
\newcommand{\vpp}{\texttt{Vem++}}
\newcommand{\cpp}{\texttt{c++}}

\newcommand{\dO}{~\text{d}\Omega}
\newcommand{\dF}{~\text{d}F}
\newcommand{\dP}{~\text{d}P}
\newcommand{\de}{~\text{d}e}
\newcommand{\dE}{~\text{d}E}
\newcommand{\dx}{~\text{d}\mathbf{x}}
\newcommand{\R}{\mathbb{R}}
\renewcommand{\P}{\mathbb{P}}
\newcommand{\nn}{\mathbf{n}}
\newcommand{\xx}{\mathbf{x}}
\newcommand{\uu}{\mathbf{u}}
\newcommand{\HH}{\mathbf{H}}
\newcommand{\B}{\mathbf{B}}
\newcommand{\jj}{\mathbf{j}}
\newcommand{\PiNMat}{\bm{\Pi}^\nabla}
\newcommand{\PiNStarMat}{\bm{\Pi}^\nabla_*}
\renewcommand{\div}{\textrm{div}}
\newcommand{\bcurl}{\textrm{\textbf{curl}}}
\newcommand{\PiNabla}[1]{\Pi^\nabla_{#1}}
\newcommand{\veOne}{\texttt{v1}}
\newcommand{\veTwo}{\texttt{v2}}

\newcommand\Tstrut{\rule{0pt}{2.6ex}\\}         
\newcommand\Bstrut{\rule[-0.9ex]{0pt}{0pt}\\}   

\newtheorem{code}{Coding hint}[section]
\newtheorem{rem}{Remark}[section]

\title{\texttt{Vem++}, a c++ library to handle and play\\
with the Virtual Element Method}
\author{F.~Dassi}

\begin{document}

\maketitle 

\begin{abstract}
The Virtual Element Method (VEM) is an extension of the Finite Element Method (FEM) used for handling polytopal meshes. 
This paper provides a brief introduction to the VEM for a two-dimensional Laplacian problem. 
Additionally, it highlights the differences between a VEM implementation and a FEM code, 
emphasising the main challenges associated with the VEM.

Furthermore, this paper presents a possible approach to address these challenges: 
\vpp{}, a \cpp{} library specifically developed for working with VEM discretisation. 
The library is designed to handle various partial differential equations in two or three dimensions, 
arising from both academic and engineering problems. 
Its flexible design allows for the seamless integration of new features, 
such as novel polytopes' quadrature rules, solvers, and virtual element spaces.\\

\end{abstract}

\section{Introduction}

The Virtual Element Method (VEM) was introduced for the first time in 2013~\cite{volley}.
It is a generalization of the Finite Element Method (FEM) that provides some advantages.
First of all it allows the presence of general polytopal elements inside the mesh~\cite{ANTONIETTI2022110900,ANTONIETTI2022111531}.
The VEM naturally handles hanging nodes, making the refinement process and mesh merging more straightforward~\cite{DASSI2022249,freeCutting}.
Additionally, since the VEM is an extension of the FEM, these two methods are fully compatible. 
For instance, if we have a two-dimensional mesh consisting of triangles, squares, and polygons, 
we can use FEM spaces for triangles and squares, while exploiting VEM spaces 
for polygons \emph{without any} special techniques~\cite{freeCutting}.

However, the advantages of the VEM extend beyond just meshing.
Indeed, the VEM spaces can be built in a way that the discrete functions inherit some useful properties.
In~\cite{vaccaDivFree} the authors introduced a family of the VEM spaces specifically for the Stokes problem 
that are exactly \say{divergence-free}.
This is in contrast of the classical FEM approaches 
where such a property is achieved in a relaxed/weak sense.
Another example about the flexibility of the VEM is given in~\cite{DASSI2020112910}. 
In this work the VEM spaces are utilized 
to solve elasticity equations based on the Hellinger–Reissner variational principle. 
The FEM enforces the symmetry of the stress tensor variable in a weak sense~\cite{boffi2013mixed}, but 
such constraint increases the dimension of the linear system at hand
and so the complexity in its resolution.
The VEM overcomes this issue:
the symmetry of the tensor field is imposed within the discrete functional spaces.
The result is a cheaper linear system to solve
since the VEM does not have additional equations 
to impose the symmetry of the tensor field.
The resulting linear system is more cost-effective 
since the VEM does not require additional equations to enforce this property.

Due to all these advantages,
the VEM has achieved significant success among researchers in various fields.
It has been successfully applied in
solid mechanics~\cite{artioli2017high,curviMec,PartI,PartII,DASSI2020112910},
Maxwell equations~\cite{BEIRAODAVEIGA202282,freeCutting},
wave equations~\cite{dassi2022virtual}
aeronautics~\cite{brick,beirao2021vorticity},
reservoir simulations~\cite{andersen2017virtual,fumagalli2019dual} and 
discrete fractional network~\cite{benedetto2014virtual,benedetto2016hybrid}.

However, despite the VEM being an extension of the FEM, 
it presents several coding challenges.
Firstly, since it can deal with polygons or polyhedrons 
the elements' data structure has 
to be as general as possible.
Furthermore, since classical quadrature rules are designed for specific element types, 
such as triangles or tetrahedrons,
the VEM needs an extension of such formulas to polygons or polyhedrons 
to compute integrals.
Finally, basis functions are no more polynomials, 
but functions known only via their degrees of freedom.
As a consequence, a way to manage such \emph{virtual functions} needs to be found.  

All of these challenges were addressed in the development of \cpp{} library \vpp{}.
The main goal of this paper is to provide readers and future code developers 
a possible solution to these issues 
which have been implemented within \vpp{}.
Additionally, \vpp{} was designed to be flexible and customizable,
allowing researchers to \say{play} with VEM.
Specific guidelines are provided to assist users in integrating their own polytopal quadrature rules, 
devising strategies for solving linear systems arising from the discretization of PDEs, 
and implementing VEM discretisations for specific problems.
While these guidelines are briefly described in the paper, 
there are also ten tutorials available 
to help users understand the underlying philosophy of the code and 
empower them to develop their desired functionalities within \vpp{}. 

Before going into the detail of the current paper,
it is worth mentioning that there are several implementations of the VEM documented in the literature.
In~\cite{Sutton2016}, Sutton proposes a concise implementation of the VEM using \emph{only} 50 lines of Matlab code.
Another Matlab implementation of the high-order VEM is presented in~\cite{Herrera2022}.
In this case, the goal of the authors is to furnish an easy-to-follow implementation
to make a first step in the VEM framework. 
Additionally, an important code to mention is the one proposed by A. Russo during the Dobbiaco Summer School  
\say{Theory and Practice of the Virtual Element Methods} in 2018,
which is freely available in~\cite{ezvem}.
Concerning \cpp{} development of the VEM,
two object-oriented libraries, \texttt{VEMLab} and~\texttt{Veamy}, deserve mention~\cite{Vemlab,VeamySoft}.
Furthermore, since the VEM is an extension of the FEM,
some research groups have integrated the VEM into existing FEM codes,
such as Dune, Deal II and OpenFOAM~\cite{Arndt2021,Dedner2010,Weller1998}.

The paper is organised as follows. 
In Section~\ref{sec:dics} we provide a brief description of the VEM 
using an easily understandable Poisson problem as an example. 
Section~\ref{sec:VEMFEM} is the core of the paper
where we discuss the similarities and differences between the implementation of the VEM and the FEM.
Then, Section~\ref{sec:play} focuses on explaining how the structure of~\vpp{} is designed 
to ease the creation of new local matrices, cubature formulas and linear system solvers.

Finally, in Section~\ref{sec:exe}, 
we describe the partial differential equations implemented in \vpp{} and 
present practical examples where the library has already been applied. 
It is important to note that this section does not provide an exhaustive coverage of all topics, 
as each example deserves a dedicated paper for a comprehensive understanding. 
However, we provide specific references for each example, 
enabling readers to access more detailed descriptions.

\paragraph{Notations.}
Throughout this paper we will follow the usual notations for Sobolev and Hilbert spaces as well as their norms~\cite{Adams:1975}.

When we are referring to polyhedrons 
we use $P$; we denote their faces with $F$. 
The outward-pointing normal of the face $F$ with respect to the polyhedron $P$ is denoted by $\nn_F$.
The volume of a polyhedron is $|P|$, 
its diameter is $h_P$ and the coordinates of its barycenter are $x_P$, $y_P$ and $z_P$.  
A generic two dimensional polygon is referred by $E$
while $e$ is always a generic edge.
As for the polyhedron, we refer to the area of a polygon $E$ as $|E|$, $h_E$ is its diameter while 
the coordinates of its barycenter will be $x_E$ and $y_E$.
We will use $\partial$ operator to denote the boundary of a polygon or polyhedron,
more specifically $\partial P$ is the set of faces that compose the boundary of $P$
while $\partial E$ is the set of edges of a generic polygon $E$.

Given a generic open domain $\mathcal{O}\subset\R^d$, with $d=1,2$ and $3$,
we refer to the polynomials of degree $k$ defined on $\mathcal{O}$ as $\P_k(\mathcal{O})$
and we denote $\pi_k$ its dimension.
In these polynomial spaces $\{m_\alpha\}$ will be the base of the scaled monomials~\cite{autostoppisti}.
Moreover, a discrete function will always be referred with the subscript \say{$h$},
for instance $v_h$, 
and the basis functions coming from the virtual element approximation are always denoted by~$\phi_i$.

The Virtual Element Method (VEM) was introduced for the first time in 2013~\cite{volley}.
It is a generalization of the Finite Element Method (FEM) that provides some advantages.
First of all it allows the presence of general polytopal elements inside the mesh~\cite{ANTONIETTI2022110900,ANTONIETTI2022111531}.
The VEM naturally handles hanging nodes, making the refinement process and mesh merging more straightforward~\cite{DASSI2022249,freeCutting}.
Additionally, since the VEM is an extension of the FEM, these two methods are fully compatible. 
For instance, if we have a two-dimensional mesh consisting of triangles, squares, and polygons, 
we can use FEM spaces for triangles and squares, while exploiting VEM spaces 
for polygons \emph{without any} special techniques~\cite{freeCutting}.

However, the advantages of the VEM extend beyond just meshing.
Indeed, the VEM spaces can be built in a way that the discrete functions inherit some useful properties.
In~\cite{vaccaDivFree} the authors introduced a family of the VEM spaces specifically for the Stokes problem 
that are exactly \say{divergence-free}.
This is in contrast of the classical FEM approaches 
where such a property is achieved in a relaxed/weak sense.
Another example about the flexibility of the VEM is given in~\cite{DASSI2020112910}. 
In this work the VEM spaces are utilized 
to solve elasticity equations based on the Hellinger–Reissner variational principle. 
The FEM enforces the symmetry of the stress tensor variable in a weak sense~\cite{boffi2013mixed}, but 
such constraint increases the dimension of the linear system at hand
and so the complexity in its resolution.
The VEM overcomes this issue:
the symmetry of the tensor field is imposed within the discrete functional spaces.
The result is a cheaper linear system to solve
since the VEM does not have additional equations 
to impose the symmetry of the tensor field.
The resulting linear system is more cost-effective 
since the VEM does not require additional equations to enforce this property.

Due to all these advantages,
the VEM has achieved significant success among researchers in various fields.
It has been successfully applied in
solid mechanics~\cite{artioli2017high,curviMec,PartI,PartII,DASSI2020112910},
Maxwell equations~\cite{BEIRAODAVEIGA202282,freeCutting},
wave equations~\cite{dassi2022virtual}
aeronautics~\cite{brick,beirao2021vorticity},
reservoir simulations~\cite{andersen2017virtual,fumagalli2019dual} and 
discrete fractional network~\cite{benedetto2014virtual,benedetto2016hybrid}.

However, despite the VEM being an extension of the FEM, 
it presents several coding challenges.
Firstly, since it can deal with polygons or polyhedrons 
the elements' data structure has 
to be as general as possible.
Furthermore, since classical quadrature rules are designed for specific element types, 
such as triangles or tetrahedrons,
the VEM needs an extension of such formulas to polygons or polyhedrons 
to compute integrals.
Finally, basis functions are no more polynomials, 
but functions known only via their degrees of freedom.
As a consequence, a way to manage such \emph{virtual functions} needs to be found.  

All of these challenges were addressed in the development of \cpp{} library \vpp{}.
The main goal of this paper is to provide readers and future code developers 
a possible solution to these issues 
which have been implemented within \vpp{}.
Additionally, \vpp{} was designed to be flexible and customizable,
allowing researchers to \say{play} with VEM.
Specific guidelines are provided to assist users in integrating their own polytopal quadrature rules, 
devising strategies for solving linear systems arising from the discretization of PDEs, 
and implementing VEM discretisations for specific problems.
While these guidelines are briefly described in the paper, 
there are also ten tutorials available 
to help users understand the underlying philosophy of the code and 
empower them to develop their desired functionalities within \vpp{}. 

Before going into the detail of the current paper,
it is worth mentioning that there are several implementations of the VEM documented in the literature.
In~\cite{Sutton2016}, Sutton proposes a concise implementation of the VEM using \emph{only} 50 lines of Matlab code.
Another Matlab implementation of the high-order VEM is presented in~\cite{Herrera2022}.
In this case, the goal of the authors is to furnish an easy-to-follow implementation
to make a first step in the VEM framework. 
Additionally, an important code to mention is the one proposed by A. Russo during the Dobbiaco Summer School  
\say{Theory and Practice of the Virtual Element Methods} in 2018,
which is freely available in~\cite{ezvem}.
Concerning \cpp{} development of the VEM,
two object-oriented libraries, \texttt{VEMLab} and~\texttt{Veamy}, deserve mention~\cite{Vemlab,VeamySoft}.
Furthermore, since the VEM is an extension of the FEM,
some research groups have integrated the VEM into existing FEM codes,
such as Dune, Deal II and OpenFOAM~\cite{Arndt2021,Dedner2010,Weller1998}.

The paper is organised as follows. 
In Section~\ref{sec:dics} we provide a brief description of the VEM 
using an easily understandable Poisson problem as an example. 
Section~\ref{sec:VEMFEM} is the core of the paper
where we discuss the similarities and differences between the implementation of the VEM and the FEM.
Then, Section~\ref{sec:play} focuses on explaining how the structure of~\vpp{} is designed 
to ease the creation of new local matrices, cubature formulas and linear system solvers.

Finally, in Section~\ref{sec:exe}, 
we describe the partial differential equations implemented in \vpp{} and 
present practical examples where the library has already been applied. 
It is important to note that this section does not provide an exhaustive coverage of all topics, 
as each example deserves a dedicated paper for a comprehensive understanding. 
However, we provide specific references for each example, 
enabling readers to access more detailed descriptions.

\paragraph{Notations.}
Throughout this paper we will follow the usual notations for Sobolev and Hilbert spaces as well as their norms~\cite{Adams:1975}.

When we are referring to polyhedrons 
we use $P$; we denote their faces with $F$. 
The outward-pointing normal of the face $F$ with respect to the polyhedron $P$ is denoted by $\nn_F$.
The volume of a polyhedron is $|P|$, 
its diameter is $h_P$ and the coordinates of its barycenter are $x_P$, $y_P$ and $z_P$.  
A generic two dimensional polygon is referred by $E$
while $e$ is always a generic edge.
As for the polyhedron, we refer to the area of a polygon $E$ as $|E|$, $h_E$ is its diameter while 
the coordinates of its barycenter will be $x_E$ and $y_E$.
We will use $\partial$ operator to denote the boundary of a polygon or polyhedron,
more specifically $\partial P$ is the set of faces that compose the boundary of $P$
while $\partial E$ is the set of edges of a generic polygon $E$.

Given a generic open domain $\mathcal{O}\subset\R^d$, with $d=1,2$ and $3$,
we refer to the polynomials of degree $k$ defined on $\mathcal{O}$ as $\P_k(\mathcal{O})$
and we denote $\pi_k$ its dimension.
In these polynomial spaces $\{m_\alpha\}$ will be the base of the scaled monomials~\cite{autostoppisti}.
Moreover, a discrete function will always be referred with the subscript \say{$h$},
for instance $v_h$, 
and the basis functions coming from the virtual element approximation are always denoted by~$\phi_i$.


\section{Discretization of a Poisson Problem in 2d}\label{sec:dics}

The \vpp{} library can handle various types of two and three dimensional PDEs;
in Secion~\ref{sec:exe} we will collect part of them.
Here, we focus on the Poisson equation in two dimensios, i.e.,
\begin{equation}
\left\{
\begin{array}{rl}
-\Delta u=f     &\text{in}\:\Omega\\
u=0             &\text{on}\:\partial\Omega 
\end{array}
\right.\,,
\label{eqn:poiss}
\end{equation}
where $\Omega$ is a domain in $\R^2$ and $f\in L^2(\Omega)$.
Although this problem is easy,
it clearly shows the essential features of VEM so 
that the reader can focus on the implementation details and issues.

The numerical approximation of the solution to problem~\eqref{eqn:poiss} in the VEM is obtained 
via a standard Galerkin framework~\cite{rabier2007theory}.
Firstly we introduce a continuous functional space $V$ and 
construct its variational formulation:
\begin{center}
\begin{minipage}{0.7\textwidth}
\emph{find $u\in H^1(\Omega)$ such that}
\end{minipage}
\end{center}
\begin{equation}
\int_\Omega \nabla u \cdot \nabla v \dO = \int_\Omega f v \dO\qquad \forall v\in  H^1_0(\Omega)\,.
\label{eqn:varForm}
\end{equation}
Let $\Omega_h$ be a decomposition of the domain $\Omega$ into general polygons and
we define a discrete space $V_h\subset V$ over such a mesh.
Then, the continuous bilinear form at the left hand-side of Equation~\eqref{eqn:varForm} 
is split on the sum over each mesh element, i.e., 
\begin{equation}
a(\cdot,\cdot)\approx a_h(\cdot,\cdot)=\sum_{E\in\Omega_h} a_h^E(\cdot,\cdot)\,,
\label{eqn:localGradOp}
\end{equation}
where $E$ is a polygon of $\Omega_h$.
Regarding the right-hand side,
the scalar product is also split over the mesh elements and 
a suitable approximation of the load term $f$, denoted as $f_h$, is constructed.

In the subsequent paragraphs, we will provide a concise overview of 
how the virtual element method constructs the discrete space $V_h$ and the bilinear form $a_h^E(\cdot,~\cdot)$.
For a more comprehensive explanation, we refer to~\cite{autostoppisti,volley}.

\subsection{Discrete space $V_h$ in two dimensions}\label{sec:dof}

In the FEM, the discrete space $V_h$ is split among the mesh elements, and 
such a space restricted on a mesh element is a two dimensional polynomial space of a specific degree.
In the VEM, $V_h$ is still split over mesh elements, but 
its restriction on an element is no longer a polynomial space, 
it is defined as
\begin{equation}
V_h(E) :=\left\{v_h\in H^1(E)\cap C^0(\partial E)\::\: \Delta v_h\in\mathbb{P}_{k-2}(E), \quad v_h|_e\in\mathbb{P}_k(e)\:\forall e\in\partial E\right\}\,.
\label{eqn:vemLocalSpace}
\end{equation}
However, if we look at the definition of $V_h(E)$, 
we notice that this space still contains polynomials of degree $k$,
but it is \say{enriched} by other unknown functions
that satisfy two conditions:
their Laplacian is a polynomial of degree $k-2$ and 
they are polynomials of degree $k$ on $\partial E$.
The inclusion of these unknown (virtual) functions is a fundamental aspect of the VEM.
Indeed, their presence enables the method 
to effectively handle arbitrarily shaped polygons, including those with hanging nodes.

Since basis functions in a FEM framework are polynomials of a given degree, 
point-wise evaluations are sufficient to uniquely determine them.
However, in the VEM, basis functions are solution of a PDE, c.f. Equation~\eqref{eqn:vemLocalSpace}, 
so a more general concept of degrees of freedom is considered.
In the VEM, degrees of freedom are proper linear functionals that uniquely determine the function~$v_h$.

In this specific case the function $v_h$ is uniquely determined 
by the following degrees of freedom:
\begin{itemize}
    \item[\texttt{D1:}] the value of $v_h$ at the vertexes of $E$;
    \item[\texttt{D2:}] the value of $v_h$ at $k-1$ points inside each edge $e\in \partial E$;
    \item[\texttt{D3:}] the moments of $v_h$ inside the polygon $E$
    \begin{equation}
    \frac{1}{|E|}\int_E v_h\,m_\alpha\dE\,.
    \label{eqn:mom}  
    \end{equation}
\end{itemize}

\begin{figure}[!htb]
    \centering
    \begin{tabular}{ccc}
    $k=1$     &$k=2$  &$k=3$ \\ 
    \includegraphics[width=0.22\textwidth]{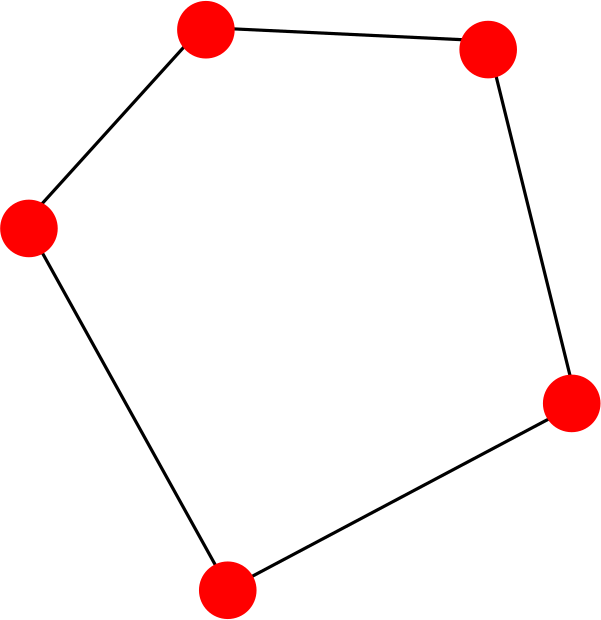} &
    \includegraphics[width=0.22\textwidth]{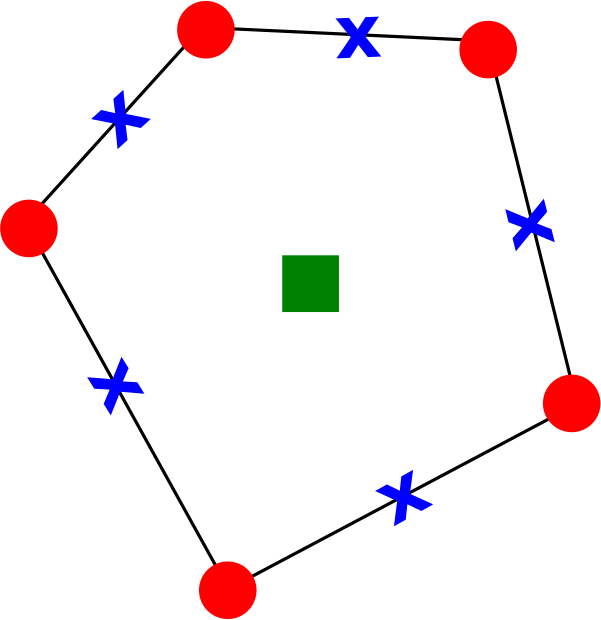} &
    \includegraphics[width=0.22\textwidth]{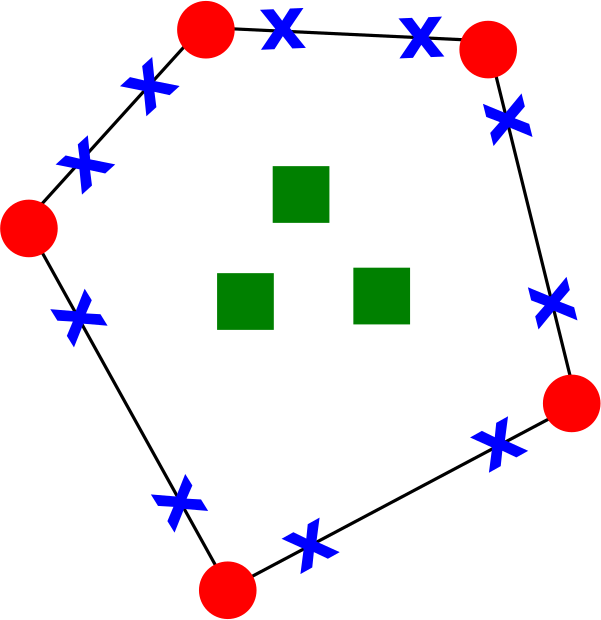} \\
    (a) &(b) &(c)
    \end{tabular}
    \caption{An example of a polygon where we highlight the dof for a VEM approximation degree~$k=1,2$ and 3.}
    \label{fig:polyWitdofs}
\end{figure}

In Figure~\ref{fig:polyWitdofs} we show the degrees of freedom for a function $v_h\in V_h(E)$ 
with varying values of $k$.
For $k=1$ the virtual functions are both harmonic and polynomial of degree 1 over $\partial E$.
As a consequence the function $v_h$ is uniquely determined by the values of $v_h$ at the vertexes of $E$, 
see the red bullets in Figure~\ref{fig:polyWitdofs}~(a).

For $k=2$, $v_h$ is a polynomial of degree 2 over edges and its Laplacian is a constant polynomial.
Then, to uniquely determine $v_h$ we need 3 points over edges, 
the values of $v_h$ at the edge endpoints, the red bullets in Figure~\ref{fig:polyWitdofs}~(b), 
and the value of $v_h$ at the edge midpoints, the blue crosses Figure~\ref{fig:polyWitdofs}~(b).
Additionally, we need an internal moment to fix the constant Laplacian, the green square in Figure~\ref{fig:polyWitdofs}~(b). 
In this specific case this degrees of freedom is the cell average of the function $v_h$.

Then, in Figure~\ref{fig:polyWitdofs}~(c) we show the degrees of freedom for the case $k=3$.
Since $v_h$ is a polynomial of degree 3 over each edge of $\partial E$.
Therefore, we need four points evaluations to determine such polynomials,
two red bullets and two blue crosses. 
Furthermore, 
to determine the Laplacian of $v_h$, which is a polynomial of degree 1,
three more conditions are required.
Such additional degrees of freedom are the cell average and the first order moments.

\begin{rem}
The choice of \texttt{D2} is not mandatory.
Indeed, one could have used moments up to degree $k-2$.
In general, any set of parameters that, together with the vertex values, 
identify polynomial of degree $k$ on each edge can replace the degrees of freedom \texttt{D2}.
\label{rem:femVem}
\end{rem}

\begin{rem}
As for the FEM, in the virtual element setting a function of the \say{canonical} basis 
is associated with a specific degree of freedom, i.e.,
there is one particular degree of freedom that is 1, 
while all the other ones are set to zero.
However, since \texttt{D1}, \texttt{D2} and  \texttt{D3} are not only point-wise evaluations,
the virtual basis functions have specific characteristics according to their definition.
Consider, for instance, a pentagon and the VEM approximation degree $k=3$.
The basis function $\phi_A$ 
that is associated with the vertex $A$,
has the following characteristics:
\begin{itemize}
    \item it is zero in all the degrees of freedom \texttt{D1} but the one associated with $A$;
    \item it is zero in all the degrees of freedom \texttt{D2};
    \item the following relations hold
    $$
    \frac{1}{|E|}\int_E \phi_A \dE = 0\,,\qquad
    \frac{1}{|E|}\int_E \phi_A \left(\frac{x-x_E}{h_E}\right)\dE = 0\,\quad\text{and}\quad
    \frac{1}{|E|}\int_E \phi_A \left(\frac{y-y_E}{h_E}\right)\dE = 0\,.
    $$
\end{itemize}
Instead, the basis function $\phi_{\mathcal{M}}$ 
associated with the \say{$x$-moment} has the following properties:
\begin{itemize}
    \item its degrees of freedom \texttt{D1} are zero;
    \item it is zero in all the degrees of freedom \texttt{D2};
    \item the following relations hold
    $$
    \frac{1}{|E|}\int_E \phi_A \dE = 0\,,\qquad
    \frac{1}{|E|}\int_E \phi_A \left(\frac{x-x_E}{h_E}\right)\dE = 1\,\quad\text{and}\quad
    \frac{1}{|E|}\int_E \phi_A \left(\frac{y-y_E}{h_E}\right)\dE = 0\,.
    $$
\end{itemize}
We compute both $\phi_A$ and $\phi_{\mathcal{M}}$
and display their shape in Figure~\ref{fig:basis} (a) and (b), respectively.
\begin{figure}[!htb]
    \centering
    \begin{tabular}{cc}
    $\phi_A$     &$\phi_{\mathcal{M}}$  \\ 
    \includegraphics[width=0.37\textwidth]{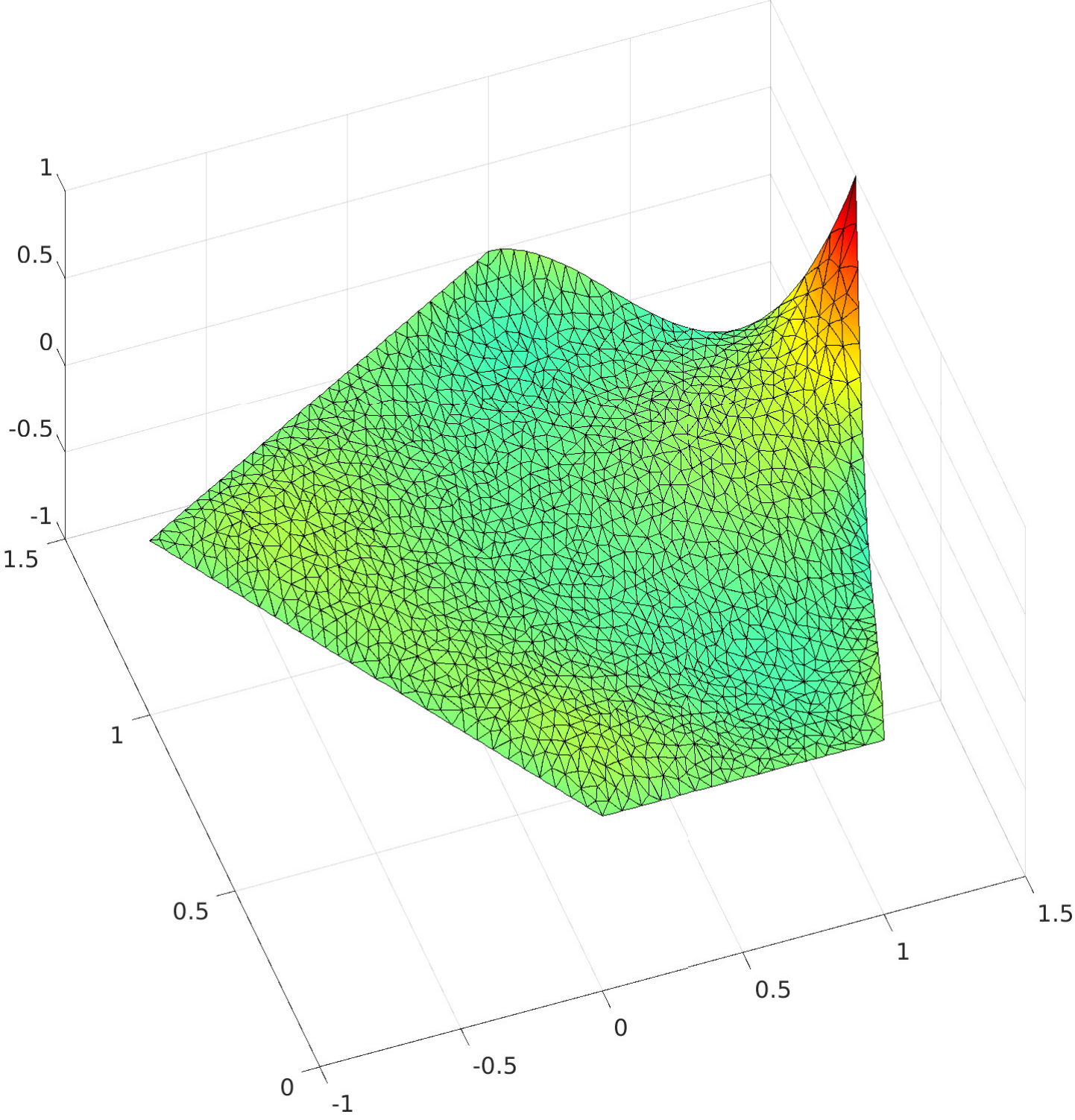} &
    \includegraphics[width=0.37\textwidth]{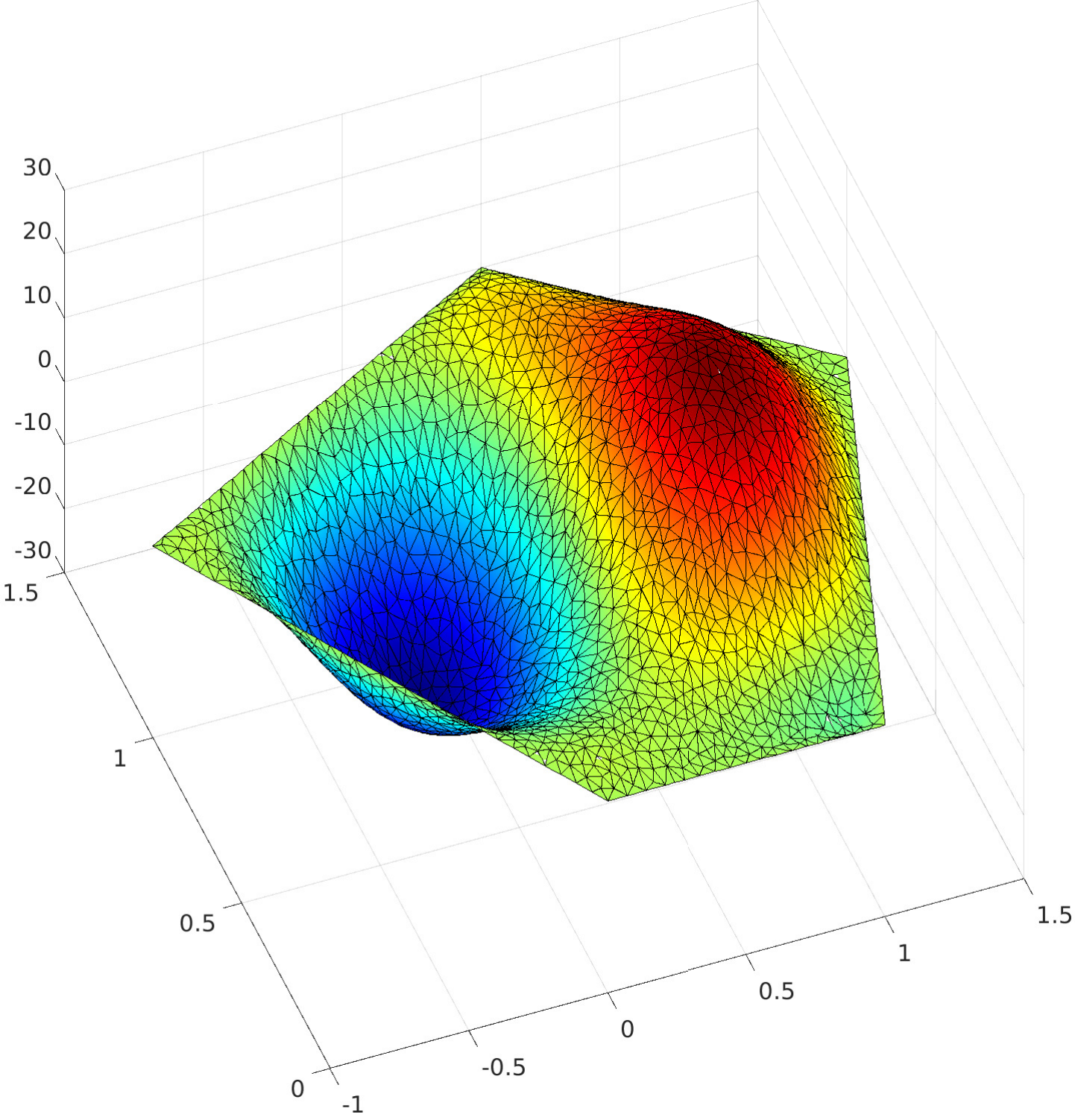} \\
    (a) &(b) 
    \end{tabular}
    \caption{The shapes of the basis functions $\phi_A$ and $\phi_{\mathcal{M}}$.}
    \label{fig:basis}
\end{figure}
\label{rem:dofs}
In the VEM virtual basis functions are typically considered virtual and are not explicitly computed. 
However, in this particular case, we compute them 
to gain an understanding of their shape and characteristics.
\end{rem}

In the upcoming sections, we will observe that even without an explicit expression of the functions in $V_h(E)$,
we can still discretise the problem defined in Equation~\eqref{eqn:varForm} based on their degrees of freedom.

\begin{figure}[!htb]
    \centering
    \includegraphics[width=0.25\textwidth]{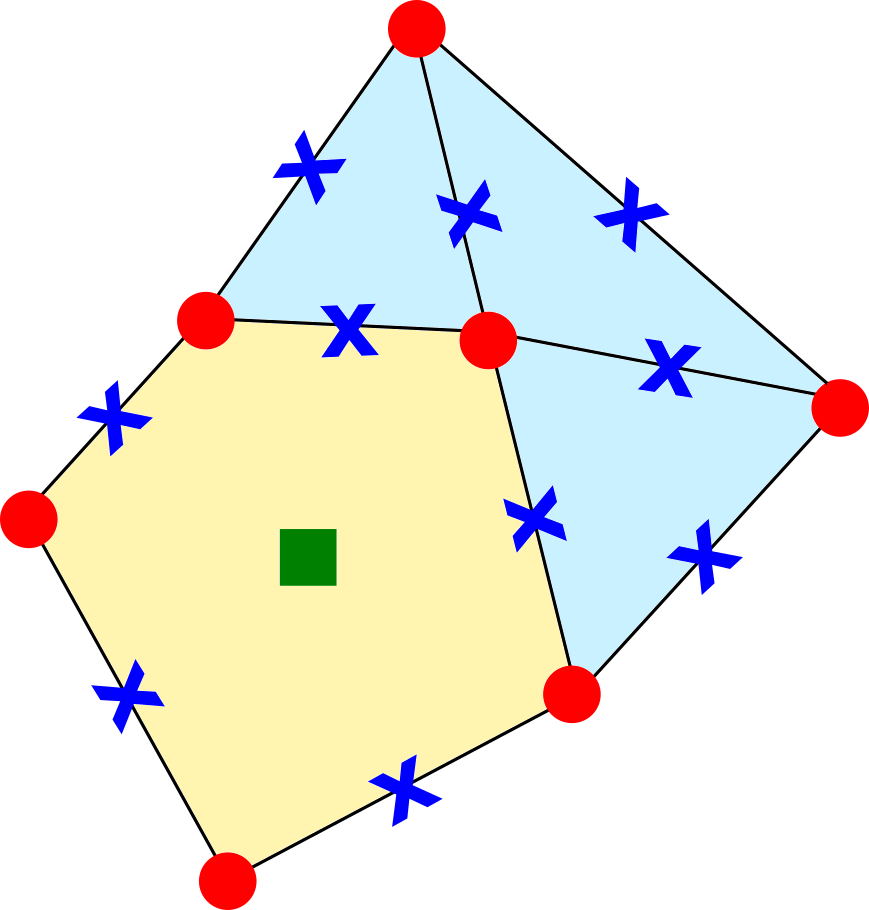}
    \caption{Merging between a VEM (light yellow) and FEM  (light blue) element when we are considering an approximation degree $k=2$.
    We further underline that the compatibility is possible \emph{for each} degree $k$.}
    \label{fig:FEMVEMele}
\end{figure}

Before going into the detail on how we define the bilinear forms,
it is important to highlight a crucial aspect of the VEM. 
Based on the definition of $V_h(E)$ and its degrees of freedom,
we observe that a function $v_h\in V_h(E)$ is a polynomial of degree $k$ on $\partial E$.
Then, according to Remark~\ref{rem:femVem}, 
one can take the FEM degrees of freedom on edges as degrees of freedom of the type \texttt{D2}.
This choice enables the combination of FEM and VEM spaces. 
In a PDE discretization, one can utilise VEM exclusively for polygonal elements with more than four edges, 
while employing FEM for all triangular or quadrilateral elements, 
see Figure~\ref{fig:FEMVEMele}.

\subsection{Projector Operator}\label{sec:proj}

In this section, we will introduce the projection operator $\PiNabla{k}$.
We will see that such projection operator plays a crucial role in constructing the bilinear operator $a_h^E(\cdot,\cdot)$.
Given a virtual function $\phi_i\in V_h(E)$, we define the projection $\PiNabla{k}: V_h(E)\to\P_k(E)$ as
\begin{equation}
\left\{
\begin{array}{rcl}
     \int_E \nabla \PiNabla{k}v_h\cdot\nabla p_k\dE &=& \int_E \nabla v_h\cdot\nabla p_k\dE\qquad\forall p_k\in\P_k(E) \\[1em]
     \int_{\partial E} \PiNabla{k}v_h \de&=&  \int_{\partial E}v_h \de
\end{array}
\right.\,.
\label{eqn:piNablaDef}
\end{equation}
The projection $\PiNabla{k} v_h$ is a polynomial of degree $k$ and 
it can be computed with the degrees of freedom of $v_h$, i.e.,
we can compute $\PiNabla{k}v_h$ \emph{without explicitly} knowing the function $v_h$ itself.
This fact is well-established in the VEM literature, see for instance~\cite{volley,autostoppisti}.
We will present a comprehensive explanation of the computation process for $\PiNabla{k}v_h$ 
in Section~\ref{sub:DiffVem}, including its implementation in \vpp{}.

We can also define an $L^2$ projection operator $\Pi^0_{k-2}$
which will play a crucial role in approximating the right-hand side of Equation~\eqref{eqn:poiss}.
Given a function $v_h\in V_h(E)$ we define the operator $\Pi^0_{k-2}:V_h(E)\to\P_{k-2}(E)$ as 
\begin{align}
\Pi^0_{0}v_h &= \frac{1}{N_V}\sum_{i=1}^{N_V} v_h(V_i) &\text{for }k=1\,,  \label{eqn:piZeroDefK1} \\
\int_E \Pi^0_{k-2}v_h p_{k-2}\dE &= \int_E v_h p_{k-2}\dE\qquad \forall p_{k-2}\in\P_{k-2}(E)&\text{for }k\geq 2,
\label{eqn:piZeroDefKgen}
\end{align}
where $N_V$ is the number of vertexes of $E$ and $V_i$ is the $i-$th vertex of the polygon $E$. 
The computability of the projection operator $\Pi^0_{k-2}v_h$ is well-established in the literature~\cite{volley,autostoppisti}.
However, it is worth to notice that the definition of such a projection operator depends on the degree $k$. 
This fact is mainly due to the computability of $\Pi^0_{k-2}v_h$. 
Indeed, if $k\geq 2$ the right-hand sides of Equation~\eqref{eqn:piZeroDefKgen} are exactly the degrees of freedom of $v_h$ so we know their values, see Figure~\ref{fig:polyWitdofs} (b) and (c).
While, if $k=1$, we do have the moments degrees of freedom,
we know only the value of the virtual function $v_h$ at the vertexes of the polygon, 
see Figure~\ref{fig:polyWitdofs} (a).
As a consequence, since one of the basic principle of the VEM is \emph{not} explicitly compute the virtual function, 
we need a different (and computable) way to define the $L^2$ projection operator.
One of the most easy way to accomplish this goal is defining such operator as 
the mean value of the function over the polygon vertexes.
Furthermore,  by appropriately modifying the virtual element spaces,
it becomes feasible to define an $L^2$ projector onto $\P_{k}$ instead of $\:\P_{k-2}$.
This technique is commonly referred to as the \emph{enhancing procedure}. 
From a theoretical standpoint,
such procedure is a complex topic that cannot be succinctly explained.
However, from a coding perspective, once you have a grasp of how to construct the $\PiNabla{k}$ projection,
it becomes apparent how to build the $\Pi_k^0$ as well~\cite{autostoppisti}.

\begin{rem}
The VEM always involves the use of a problem-specific projection operator.
For example, in the discretisation of a Stokes problem, 
a suitable $L^2$ projection operator $\bm{\Pi}^0_k$ is defined to approximate the mass operator~\cite{brick}.
On the other hand, for an elasticity problem, a different projection operator $\bm{\Pi}^\epsilon_k$
is employed to obtain the virtual element approximation of the symmetric gradient~\cite{brick}.
However, in all these cases, the projection operator can be computed solely based on the degrees of freedom, without requiring an explicit expression of the virtual function~$v_h$.
\end{rem}

\subsection{Problem Operators}\label{sec:operator}

After introducing the projection operator we are now ready to present the VEM approximation of the local bilinear operator $a_h^E(\cdot,\cdot)$, Equation~\eqref{eqn:localGradOp}, and the right hand side of Equation~\eqref{eqn:varForm}.
Let us consider the projection operator $\PiNabla{k}\phi_i$ of a generic virtual element basis function.
We apply the continuous local operator $a^E(\cdot,\cdot)$ to 
$$
\PiNabla{k}\phi_i + \left( \phi_i - \PiNabla{k}\phi_i\right)\,, 
$$
then, using the bi-linearity of such operator and the orthogonality property,
we get the form
$$
a\left(\phi_i,\,\phi_j\right) = a\left(\PiNabla{k}\phi_i,\,\PiNabla{k}\phi_j\right) + a\left(\phi_i -\PiNabla{k}\phi_i,\,\phi_j-\PiNabla{k}\phi_j\right)\,. 
$$
Since we are able to compute the projection $\PiNabla{k}$, we can handle the first term.
The second one still contains virtual functions so it can not be computed.
However, such term can be replaced with \emph{any} positive definite bilinear form $s^E(\cdot,\,\cdot)$ that
satisfies for each $v_h\in V_h\cap\ker(\PiNabla{k})$
$$
\alpha_*\,a^E(v_h,\,v_h)\leq s^E(v_h,\,v_h) \leq \alpha^*\,a^E(v_h,\,v_h)\,,
$$
where $\alpha_*$ and $\alpha^*$ are positive constants.
In the virtual element setting, there exist various choices for such an operator, 
and all of them can be computed even if we do not know explicitly the virtual function $v_h$. 
In Section~\ref{sec:localMat} we will give one of the most common choice,
the so-called \say{\texttt{dofi-dofi}} stabilisation~\cite{autostoppisti}.
Based on this understanding, we are now ready to define the local discrete operator
\begin{equation}
a_h^E(\phi_i,\,\phi_j):= 
a\left(\PiNabla{k}\phi_i,\,\PiNabla{k}\phi_j\right) + s^E\left(\phi_i -\PiNabla{k}\phi_i,\,\phi_j-\PiNabla{k}\phi_j\right)\,. 
\label{eqn:localForm}
\end{equation}
The first part is commonly referred to as the \emph{consistency} term
that preserves the accuracy of the method.
The second part is known as the \emph{stability} term and
it guarantees the proper convergence rate of the solution.

Then, to get a suitable approximation of the right hand side, we use the $\Pi_{k-2}^0$ projection operator.
Since the definition of such operator depends on the degree $k$,
we distinguish two cases:
\begin{itemize}
\item if $k=1$ we approximate the integral using the $L^2$ projection defined in Equation~\eqref{eqn:piZeroDefK1} on both $f$ and $\phi_i$:
$$
\int_E \Pi_{0}^0\,f \,\Pi_{0}^0\phi_i\dE= \Pi_{0}^0\,f \,\Pi_{0}^0\phi_i |E|\,, 
$$
\item if $k\geq 2$ we apply the $L^2$ projection defined in Equation~\eqref{eqn:piZeroDefKgen} 
on the load term $f$, $f_h:=\Pi_{k-2}^0\,f$
and we proceed as follows exploiting the properties of this operator:
$$
\int_E \Pi_{k-2}^0\,f \,\phi_i\dE = \int_E \Pi_{k-2}^0\,f\, \Pi_{k-2}^0 \phi_i\dE = \int_E f\, \Pi_{k-2}^0 \phi_i\dE\,.
$$
\end{itemize}

Now we have all the ingredients to assemble the linear system associated with the virtual element approximation of the problem defined in Equation~\eqref{eqn:poiss}. 
Such assembly process is akin to the one used in the FEM.
In both cases we have local basis functions, $\phi_i$,
that defines \emph{local} linear and bilinear forms.
These local forms are then assembled into a global matrix and 
the solution is obtained via the resolution of a linear system.
While Section~\ref{sec:EqualVEM} will delve deeper into the specifics of this procedure, 
it is evident that both VEM and FEM follow a similar workflow 
to obtain a discrete solution.


\section{Implementation details}\label{sec:VEMFEM}

In this section, we will analyse the VEM from a more practical point of view. 
In Section~\ref{sub:DiffVem} we will discuss the main challenges in implementing the VEM compared to the FEM.
Subsequently, in Section~\ref{sec:EqualVEM}, we will highlight the similarities between these two methods.
Both sections will also include details about the structure of the \vpp{} library.

\subsection{Differences with the FEM}\label{sub:DiffVem}

In this subsection we will describe the main difficulties in coding the VEM.
More specifically, in Section~\ref{sec:eleShape} we focus on the difficulties related to the presence of arbitrarily shaped polygons. 
Then, in Section~\ref{sec:dofImpl}, we discuss how in \vpp{} virtual basis functions are handled. Finally, in Section~\ref{sec:integration}, we consider the problem of making integrals over polygons/polyhedrons and
how to deal with the integration of virtual functions.

\subsubsection{Element shape}\label{sec:eleShape}

The VEM can handle a mesh composed by arbitrarily shaped polygons and polyhedrons. 
As a consequence, since we do not have any reference shape,
we need to take into account all possible cases. 
For instance, given a generic polygon, 
we can not a priori say that it has a fixed number of edges or vertices.
Moreover, due to the flexibility of VEM, we can not either say that it has no hole!
As a consequence in \vpp{} polygons and polyhedrons are managed as 2d and 3d \emph{Piece-wise Linear Complex} (PLC), respectively~\cite{miller1996control}. 

\begin{figure}[!htb]
    \centering
    \begin{subfigure}[t]{0.3\textwidth}
    \centering
    \includegraphics[width=0.6\textwidth]{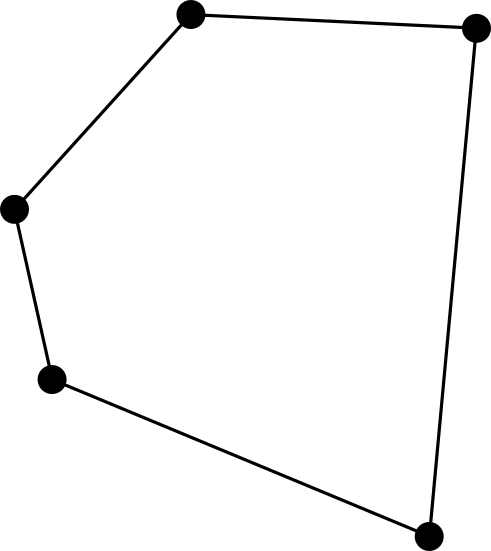}
    \caption{A \texttt{polygon}}
    \label{fig:plc2dPoly}
    \end{subfigure}
    \begin{subfigure}[t]{0.3\textwidth}
    \centering
    \includegraphics[width=0.6\textwidth]{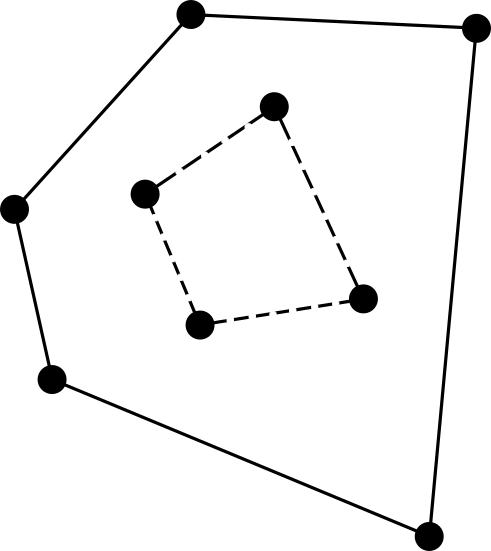}
    \caption{A set of \texttt{polygon} objects.}
    \label{fig:plc2dSetWithOr}
    \end{subfigure}
    \begin{subfigure}[t]{0.3\textwidth}
    \centering
    \includegraphics[width=0.6\textwidth]{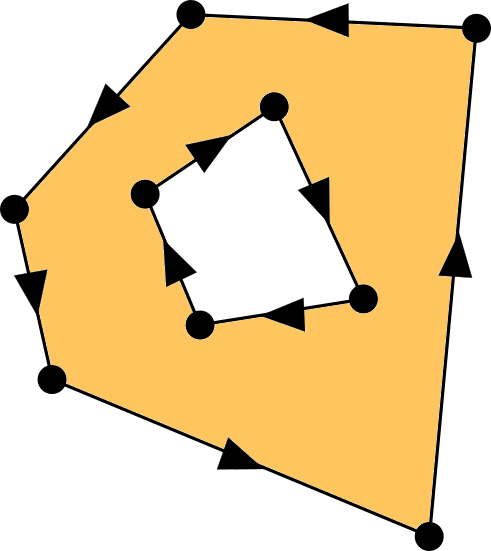}
    \caption{A \texttt{facet}.}
    \label{fig:thefinalFacet}
    \end{subfigure}
    \caption{The data structure for a 2d PLC used in \vpp{}.}
    \label{fig:2dPLC}
\end{figure}

\paragraph{\texttt{Facet}.}
The facet is the topological object used in \vpp{} to handle a generic 2d PLC.
It is composed by a collection of \texttt{polygon}, i.e., 
a sequence of points that identify a close piece-wise linear complex, see Figure~\ref{fig:2dPLC}~(a).
Among this collection we identify the main polygon that encloses all the other ones, see Figure~\ref{fig:2dPLC}~(b).
Then, the generic \texttt{facet} is defined running counter-clock-wise the main polygon and 
clock-wise the polygons which defines the holes, see Figure~\ref{fig:2dPLC}~(c). 

As for the standard FEM, in some cases a virtual element approximation requires the definition of the out-ward pointing normal.
However, thanks to the proposed orientation criterion, 
it is always possible to define out-pointing normal of each \texttt{facet}'s edge although it has holes.
Moreover, one can exploit this orientation to compute the area in a straightforward way using the so-called shoelace formula~\cite{shoelace}.

\begin{code}
The area of a generic \texttt{facet} is given by the main polygon area 
minus the area of all holes.
\end{code}

\begin{figure}[htb]
\centering
\begin{overpic}[unit=1mm,width=0.35\textwidth]{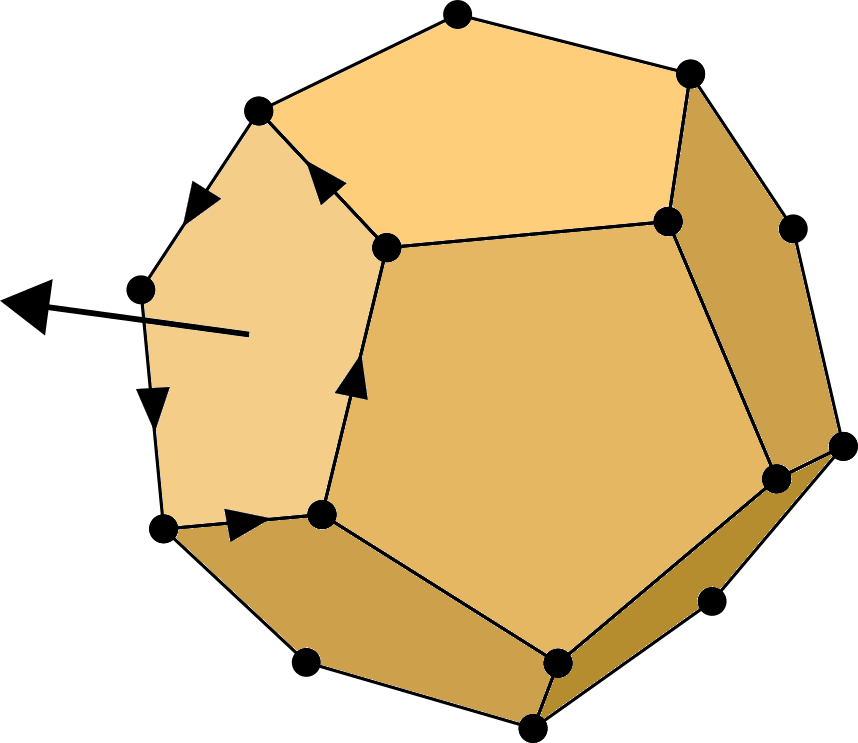}
\put(2,58){{\large $\mathbf{n}_F$}}
\put(28,33){{\Large $F$}}
\end{overpic}
\caption{An example of polyhedron in \vpp{} where we highlight the orientation of one \texttt{facet} so 
that the normal points outside.}
\label{fig:3dPLC}
\end{figure}

\paragraph{\texttt{Polyhedron}.}
To handle a polyhedron in \vpp{},
we use a generic 3d PLC data structure. 
Also in this case we can not make any a priori assumption about the shape of the element:
it can be composed by any number of vertexes, edges and faces.
Furthermore, we may have holes inside the polyhedron. 
As a consequence a polyhedron is handled by a \cpp{} class called \texttt{polyhedron}
that is a collection of \texttt{facet} and, 
eventually, a sequence of points that identify the interior of each hole 
(in both 2d and 3d mesh generation~\cite{triangle,tetgen} this way of detecting an hole via a point is quite common).
Moreover, the \texttt{facet} are oriented in such a way 
that normal constructed via the right-hand rule points outside the polyhedron, see Figure~\ref{fig:3dPLC}.

This kind of orientation plays a key role in the computation of the volume enclosed by the polyhedron.
Indeed, one can use the divergence theorem to exactly compute the volume from face integrals.
Consider a polyhedron $P$ and a vector field $\mathbf{F}$ such that its divergence is 1.
Then, the following relation holds
\begin{equation}
|P| = \int_P \dP = \int_P \text{div}(\mathbf{F})\dP = \sum_{F\in\partial P}\int_{F} \mathbf{n}_F\cdot \mathbf{F}\dF\,.
\label{eqn:vol}
\end{equation}

\begin{code}
As vector field $\mathbf{F}$ in Equation~\eqref{eqn:vol}, 
one can take $\mathbf{F}(x,\,y,\,z)=(x,\,0,\,0)^\top$.
Such choice is particularly appealing from the computational point of view.
Indeed, the vector field we are integrating over faces is a polynomial of degree 1
and consequently we may use a low order quadrature rule to ``exaclty'' compute the face integral with the smallest computational effort.
\end{code}

Before ending this section, 
we underline that the set of polyhedron's face does not have a specific order in a 3d PLC.
However, it is possible to give a specific 3d orientation of polyhedron's faces \emph{only} for some specific polyhedrons, e.g., tetrahedrons, hexagons, prism, etc.~\cite{tetgen,hexaMesh}.
Since the VEM can handle generic polyhedrons, 
finding a specific rule to orient arbitrary shaped polygons is not straightforward 
and the gain in terms code speed could be not relevant.
As a consequence in \vpp{} the \texttt{facet} which define the polyhedron 
are randomly stored inside the object \texttt{polyhedron}. 

\subsubsection{Definition of the degrees of freedom}\label{sec:dofImpl}

In a virtual element framework basis functions are not polynomials as for the FEM, 
see Section~\ref{sec:dof}.
This is the key feature of the VEM:
the basis functions are virtual,
we do not know them or,
it will be better saying,
we \emph{do not want} to have their explicit expression.
Indeed, to solve a partial differential equation via VEM,
the \emph{only} information we need are the values of virtual basis functions on degrees of freedom and the projection operators.

The structure of \vpp{} follows the same philosophy:
there is no an explicit expression of the virtual basis functions,
but they are a collection of information stored in the class \texttt{vemDofDescription}.

To better understand this fact, 
we consider the degrees of freedom described in Figure~\ref{fig:polyWitdofs}~(c), i.e.,
a VEM approximation degree $k=3$,
and we show how the class \texttt{vemDofDescription} describes them in \vpp{}.

First of all each dof has a tag so that we can recognise among each type. 
For instance, in the present case 
we can distinguish among red bullets, blue crosses and green squares. 
Then, according to the tag, the class \texttt{vemDofDescription} stores the main information to deal with the dof:
\begin{itemize}
\item[-] for red bullets it contains the id of the vertex the dof is associated with and the coordinates of such point;
\item[-] for each blue crosses, it contains the id of the edge the dof belongs to and the coordinate of the point;
\item[-] for each green square the class \texttt{vemDofDescription} contains the id of the face the dofs belongs to and the exponents of the monomial it is associated with. Notice that for this simple case we have three internal moments.
Then, the set of exponents are the couples (0,0), (1,0) and (0,1) that represent the following monomial moments:
$$
\int_E \phi_i\dE,\qquad
\int_E \phi_i\left(\frac{x-x_E}{h_E}\right)\dE\quad\text{and}\quad
\int_E \phi_i\left(\frac{y-y_E}{h_E}\right) \dE\,,
$$
respectively. 
\end{itemize}
We remark that such kind of information is related to virtual element space
we are considering. 
If we are considering a different virtual element space the information stored inside \texttt{vemDofDescription} change. 
For instance, if you are considering the VEM spaces used to solve a Stokes problem we have moments associated with the divergence of the virtual basis function~\cite{brick,beirao2020stokes}. 

\subsubsection{Integration}\label{sec:integration}

As for the FEM, also in the virtual element method one needs to compute integrals. 
However, in this framework there are two main issue.
On the one hand one has to integrate over polygons, 
on the other hand one has to deal with integrations of virtual basis function
that are unknown.

\subsubsection*{Integration over polygons.}
There are different ways to achieve this goal.
The more straightforward one is to sub-triangulate polygons and 
use triangle quadrature points and weights.
A more sophisticated strategy is to use quadrature rules built ad-hoc for general polygons/polyhedrons.
Among them we mention Vianello formulas for two dimensional polygons~\cite{vianello},
that was extended for polyhedrons in~\cite{vianello3d}, 
and Chin-Sukumar quadrature rule~\cite{ChinSuku}.
In the majority of the cases such rules result in a collection of more than needed quadrature points and,
as a consequence, an high computational effort. 

To avoid this issue, one possible strategy is to use compression rules 
that are able to reduce the number of points according to the quadrature rule precision required~\cite{compression}.
Such compression procedure are general and it is based on the resolution of a non-negative least squares problem 
that also ensures non-negativity of weights.
To show the effectiveness of the compression procedure, 
we consider a cube and the quadrature formula that exactly integrates polynomials of degree 2 described in~\cite{vianello3d} Section~5.2.
A standard rule has 192 quadrature points while the compressed one has only 10 points.
The key point is that this compression strategy always results in the minimum number of points 
that interpolate a specific polynomial. 
Indeed, in this particular case we are considering a quadrature rule of degree 2 
so the compression procedure selects and modifies the weights of 10 quadrature points 
that are exactly the number of points you need to construct a degree 2 polynomial in three variables.

\subsubsection*{Integration of virtual basis functions.}
In a virtual element setting we do not have an explicit expression of basis functions so we can not evaluate them at any point of their support. 
As a consequence one can think that it is not possible to make integrals of such functions.
However, this fact is not completely true.
Indeed, in a virtual element setting it is possible to exactly compute the essential integrals to set the method.
Moreover, the computation of such an integral is exact.  

In the following paragraphs we will give a better explanation of this fact considering three issues:
integral on edges, on faces and computing the integrals required to compute the $~\PiNabla{k}$-projection.

\paragraph{Edge integrals and virtual functions.}
Suppose that you want to compute the following integral
\begin{equation}
\int_e v_h\,g(t) \de\,,
\label{eqn:edgeIntegral}
\end{equation}
where $v_h$ is a virtual function and $g$ is a generic known function.
In the virtual element setting described in Section~\ref{sec:dics},
$v_h$ is a polynomial of degree $k$ on the boundary.
Moreover, one can reconstruct such a polynomial starting from its degrees of freedom over the edge $e$.
As a consequence it is possible to compute the integral in Equation~\eqref{eqn:edgeIntegral}
and its value is accurate according to the quadrature rule taken into account.

It often happens that the function $g\in\mathbb{P}_s(e)$, i.e., the function $g$ is a polynomial of degree $s$.
As a consequence if one choose a quadrature rule that is exact for polynomial of degree $k+s$,
the integral in Equation~\eqref{eqn:edgeIntegral} is ``exact''.

There are several virtual element spaces where the virtual function is a polynomial over edges, 
see for instance~\cite{antonietti2016c,da2018virtual}.
This is also common in the three dimensional space.
where the virtual functions are defined as polynomials on the skeleton of the polygon,
see, e.g.,~\cite{da2017high,beirao2020stokes}.

\begin{rem}
If $s\leq k-1$ the integration of Equation~\eqref{eqn:edgeIntegral} is straightforward.
Indeed, since the the degrees of freedom of $v_h$ are defined at the Gau\ss{}-Lobatto nodes over edges,
one can exploit such integration rule and get
\begin{equation}
\int_e v_h\,g(t) \de = \sum_{j=0}^{k} \omega_j\,v_h(t_j)\,g(t_j) = \sum_{j=0}^n \omega_j\,\texttt{dof}_{e,j}(v_h)\,g(t_j)\,,    
\label{eqn:integralEdgeTrik}
\end{equation}
where $\texttt{dof}_{e,j}(v_h)$ is the dof value at the $j-$th Gau\ss{}-Lobatto quadrature point on the edge $e$, 
while $t_j$ and $w_j$ are the quadrature points and weights, respectively.
\end{rem}

\begin{rem}
If $\phi_i$ is a virtual basis function associated with the $i-$th Gau\ss{}-Lobatto node of the edge $e$,
Equation~\eqref{eqn:integralEdgeTrik} further simplifies in
\begin{equation}
\int_e \phi_i\,g(t) \de = \sum_{j=0}^{k} \omega_j\,\phi_i(t_j)\,g(t_j) = \sum_{j=0}^n \omega_j\,\texttt{dof}_{e,j}(\phi_i)\,g(t_j) = \omega_j\,g(t_j)\,.    
\label{eqn:edgeIntegralPhi}
\end{equation}
Moreover, if the virtual basis function is not associated with any Gau\ss{}-Lobatto nodes of the edge $e$, i.e., it is associated with a moment \texttt{D3}, 
such integral is zero and no computation are needed.
\end{rem}

\begin{code}
In \vpp{} the degrees of freedom are \cpp{} classes 
where the most important information about the dofs are stored.
More specifically, for a degrees of freedom \texttt{D1} and \texttt{D2}
the point parameter $t_j$, the weight $\omega_j$ and the identifier of the edge are stored 
so that the integral of Equation~\ref{eqn:edgeIntegralPhi} can be computed at once. 
\end{code}

\paragraph{Face integrals and virtual functions.}
When you consider a virtual function $v_h$, 
it is not possible to compute a generic face integral of the form 
\begin{equation}
\int_E v_h\,g(x,y) \dE\,.
\label{eqn:faceIntegral}
\end{equation}
Indeed you do not know the function $v_h$ inside the polygon and
you can not evaluate it at any quadrature point inside $E$. 
However, for particular choices of $g(x,\,y)$ it is possible 
to get the \emph{exact} value of the integral in Equation~\eqref{eqn:faceIntegral} at once \emph{without} resorting to any quadrature rule since 
it is a degrees of freedom of $v_h$.

Consider for instance the virtual element space defined in Section~\ref{sec:dics} and 
the degrees of freedom defined in Section~\ref{sec:dof}.
If the virtual element approximation is $k=3$, 
we know the exact result of the following integrals 
\begin{equation}
\int_E v_h\dE,\qquad
\int_E v_h\left(\frac{x-x_E}{h_E}\right)\dE\quad\text{and}\quad
\int_E v_h\left(\frac{y-y_E}{h_E}\right) \dE\,.
\label{eqn:faceMonoIntegral}    
\end{equation}
Indeed, such integrals are the degrees of freedom of the function $v_h$ and 
since we know only the degrees of freedom of the virtual function $v_h$, 
we know in fact the values of the integrals in Equation~\eqref{eqn:faceMonoIntegral}:
we do not need any quadrature rule to compute them,
we know its exact values for free.

Moreover, if we are considering the basis function associated with such a space,
we deduce the values of the integrals in Equation~\eqref{eqn:faceMonoIntegral} from the type of the dof.
More specifically, if $\phi_i$ is a basis function associated with a dof type $\texttt{D1}$ and $\texttt{D2}$
\begin{equation*}
\int_E \phi_i\dE = 
\int_E \phi_i \left(\frac{x-x_E}{h_E}\right)\dE = 
\int_E \phi_i \left(\frac{y-y_E}{h_E}\right) \dE = 0\,.
\end{equation*}
While if $\phi_i$ is associated with a dof type $\texttt{D3}$ one of these integral will be 1,
see Remark~\ref{rem:dofs}.

\begin{code}
In \vpp{} also the $\texttt{D3}$ dofs are a \cpp{} class 
where it is stored the identifier of the face and the exponent of the monomial it is associated with
so that we are able to identify and consequently compute the integrals of Equation~\eqref{eqn:faceMonoIntegral}
at once.
\end{code}

As you can deduce from the above example,
the face integrals we are able to compute are related to both
the virtual element approximation degree and the virtual element space we are considering.
In particular, if we are considering the spaces Section~\ref{sec:dics} with a VEM approximation degree $k=2$,
we are to compute only 
$$
\int_E v_h\dE\,,
$$
since we do not have any degrees of freedom associated with the moments
$$
\left(\frac{x-x_E}{h_E}\right)\qquad\text{and}\qquad\left(\frac{y-y_E}{h_E}\right)\,.
$$

\begin{rem}
Similar consideration can be done for integral over polyhedrons 
when we are considering a three dimensional setting of the VEM.
In this specific case we will have some basis functions associated with volume moments so 
we can have the exact value of them for free~\cite{da2017high}.
\end{rem}

\paragraph{Integrals related to the computation of $~\PiNabla{k}$-projection.}

In this section we will give an explicit description on 
how it is possible to compute the $\PiNabla{k}$ projection 
although we are dealing with virtual unknown functions.
We will focus in the computation of such projection for a generic basis function $\phi_i$,
the computation of this operator for a generic virtual function follows similar arguments.
We refer the reader to~\cite{autostoppisti} for a more detailed description of this computation.
We consider the scaled monomial basis $\{m_j\}_{j=1}^{\pi_k}$ as basis of the polynomial $\PiNabla{k}\phi_i$, i.e.,
\begin{equation}
\PiNabla{k}\phi_i = \sum_{j=0}^{\pi_k} c_j\,m_j\,,
\label{eqn:projOnMono}
\end{equation}
where $c_j\in\R$ are the projection coefficients and 
$\pi_k$ is the dimension of the polynomials of degree lower or equal to $k$.
To facilitate exposure,
we recall the definition of this projection operator.
\begin{equation}
\left\{
\begin{array}{rcl}
     \int_E \nabla \PiNabla{k}\phi_i\cdot\nabla m_j\dE &=& \int_E \nabla \phi_i\cdot\nabla m_j\dE  \\[1em]
     \int_{\partial E} \PiNabla{k}\phi_i \de&=&  \int_{\partial E}\phi_i \de
\end{array}
\right.\,.
\label{eqn:piNablaDef2}
\end{equation}
Notice that here we do not consider a generic polynomial $p_k$.
Indeed, we consider the $\pi_k$ monomials $m_j$ 
that are a basis of $\P_k(E)$ so that we are able to compute the coefficients $c_j$ in Equation~\eqref{eqn:projOnMono}.

The main issue related to the computation of \emph{all} the quantities in Equation~\eqref{eqn:piNablaDef2} is related to the presence of the virtual function $\phi_i$. Indeed, in both equation the left-hand side contains polynomials so it is computable, c.f. ``Integration over polygons''.
The right-hand sides contain the virtual basis function $\phi_i$ so a priori,
we are not able to compute such an integrals.
We will see that this fact is not true and,
to compute the right hand side the knowledge of the virtual function's degrees of freedom is enough.

Consider the first condition in Equation~\eqref{eqn:piNablaDef2} and 
proceed with an integration by parts
$$
\int_E \nabla \phi_i\cdot\nabla m_j\dE = -\int_E \phi_i\Delta m_j\dE + \int_{\partial E} \phi_i\,(\nn\cdot \nabla m_j)\de\,.
$$
As a consequence the bulk integral is split into two parts:
a boundary part and face part.
The boundary integral is computable since the virtual function $\phi_i$ is a polynomial on $\partial E$, 
c.f. ``Edge integrals and virtual functions''.
However, the bulk integral is computable too. 
Indeed, we have the following identity
$$
\Delta m_j = \frac{1}{h_E^2}\left(\tilde{c}\,\tilde{m}_j+\overline{c}\,\overline{m}_j\right)\,,
$$
where $\tilde{m}_j$ and $\overline{m}_j$ are scaled monomials of degree $k-2$ and 
$\tilde{c}$ and $\overline{c}$ are proper constants depending on $m_j$. 
Then, we observe that this integral becomes the sum of two integrals 
where we recognise the definition of moments, Equation~\eqref{eqn:mom}. 
$$
-\int_E \phi_i\Delta m_j\dE = 
-\frac{\tilde{c}}{h_E^2}\int_{E}\phi_i\,\tilde{m}_j\dE 
-\frac{\overline{c}}{h_E^2}\int_{E}\phi_i\,\overline{m}_j\dE\,.
$$
As a consequence we know the result of such integrals 
and we are able to compute them, c.f. ``Face integrals and virtual functions''. 
To better understand the previous identity,
we show the following example
$$
-\int_E \phi_i\Delta\left[\left(\frac{x-x_E}{h_E}\right)^3
\left(\frac{y-y_E}{h_E}\right)^2\right]\text{d}E = 
-\frac{6}{h_E^2}\int_{E}\phi_i\,\left(\frac{x-x_E}{h_E}\right)\left(\frac{y-y_E}{h_E}\right)^2\text{d}E 
-\frac{2}{h_E^2}\int_{E}\phi_i\,\left(\frac{x-x_E}{h_E}\right)^3\text{d}E\,.
$$
Regarding the second condition in Equation~\eqref{eqn:piNablaDef2},
we have an integral over $\partial E$ of the virtual function $\phi_i$.
Once again we recall that $\phi_i$ is a known polynomial over $\partial E$ so
such we are able to compute such an integral, c.f. ``Edge integrals and virtual functions''.

\begin{code}
From the previous computations it comes to mind
that the VEM is based on monomials.
More specifically, it is clear that it requires 
the computation of their derivatives and recognising monomials. 
As a consequence in \vpp{}, 
we develop two classes \texttt{monomial2d} and \texttt{monomial3d}
that other than defining two and three dimensional monomials,
they furnish a set of basic symbolic calculus. 
For instance, given two instances of the class \texttt{monomial2d}, 
it is possible to get a new \texttt{monomial2d} that is the product between them, see line 6 of Listing~\ref{alg:mono}.
Moreover, given a \texttt{monomial2d} it is possible to get its derivatives, see lines 8 and 9 of Listing~\ref{alg:mono}.
{\footnotesize
\begin{lstlisting}[language=C++, caption=Monomial symbolic calculus, label=alg:mono]
// definition of the monomial 3x^1y^2
monomial2d mono1(1,2,3.);
// definition of the monomial 2x^2y^1
monomial2d mono1(2,1,2.);
//product between them
monomial2d mono3 = mono1 * mono2;
//product derivative of this monomial w.r.t. x and y 
monomial2d mono3Dx = mono3.dx();
monomial2d mono3Dy = mono3.dy();
\end{lstlisting}}
\end{code}

\subsection{Analogies with FEM}\label{sec:EqualVEM}

\begin{figure}[!htb]
\centering
\includegraphics[width=0.9\textwidth]{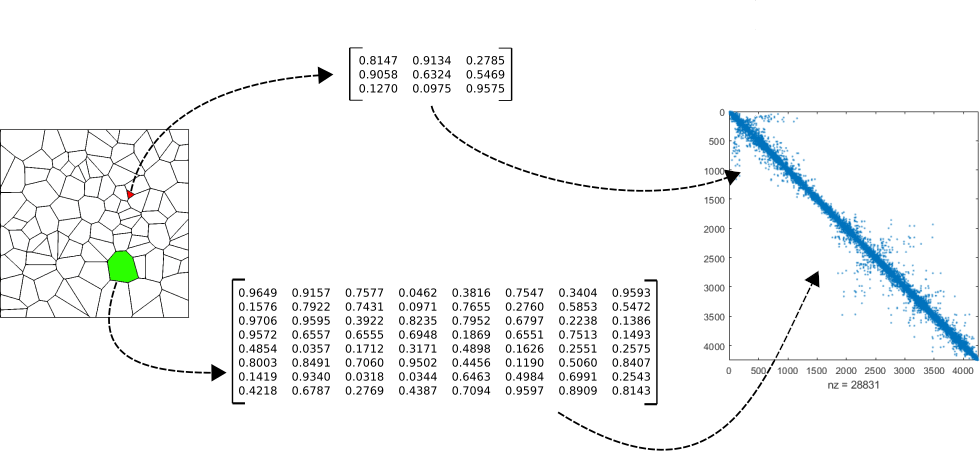}
\caption{General scheme to find the numerical solution of a PDE.}
\label{fig:genAss}
\end{figure}

The VEM is an extension of the FEM to general polytopal meshes~\cite{volley}.
As a consequence it follows the same workflow of the FEM 
to assemble and solve the linear system arising from the discretization of a~PDE, see Figure~\ref{fig:genAss}:
\begin{enumerate}[1.]
\item\label{step:1} create each local matrices that represent each bilinear or linear operators;
\item assemble such local matrix in the global one;
\item solve the linear system via a proper solver.
\end{enumerate}
The only difference stays in the local matrices.
Indeed, in a FEM framework they have always the same size and 
they are computed starting from the reference element.
On the contrary in the VEM their size depends on the shape of the element and 
to compute them one needs suitable projection operators.
From a more practical point of view,
once the code is able to handle such local matrices 
the implementation of a VEM code and a FEM code coincides.

\subsubsection{Local Matrix Definition}\label{sec:localMat}

Since the main difference between a FEM and a VEM assembling procedure 
stays in the definition of the local matrices,
in this subsection we focus on how local matrices are handled in~\vpp{}. 

Before going into the detail we underline that 
given a specific linear operator,
it can be seen as a sequence of matrix products and sums. 
Consider for instance the operator which gives the VEM discretization of the grad-grad operator:
\begin{equation}
a_h^P(\phi_i,\,\phi_j) = \int_E \nabla \Pi_k^\nabla \phi_i\cdot \nabla\Pi_k^\nabla \phi_i\dE + 
\sum_{r=1}^{\#\text{dofs}}\texttt{dof}_r\left(\left(I-\Pi_k^\nabla\right)\phi_i\right)\texttt{dof}_r\left(\left(I-\Pi_k^\nabla\right)\phi_j\right)\,.
\label{eqn:aOperator}
\end{equation}
If we define the following matrix
$$
\mathbf{G}_{ij} = \int_E \nabla m_i\cdot \nabla m_j\dP\,,
$$
and the matrices that contains the projection coefficient with respect to the monomial basis, $\PiNStarMat$, and the virtual basis $\PiNMat$, then the operator defined in Equation~\eqref{eqn:aOperator} becomes a simple sequence of matrix products
\begin{equation}
\left(\PiNStarMat\right)^t\,\mathbf{G}\,\PiNStarMat + \left(\mathbf{I}-\PiNMat\right)^t\left(\mathbf{I}-\PiNMat\right)\,,    
\label{eqn:aOperatorMatForm}
\end{equation}
where $\mathbf{I}$ is the identity matrix. 
In this paper we focus on the implementation of the VEM inside \vpp{} so 
we do not give the explicit details on how matrices $\PiNStarMat$ and $\PiNMat$ are constructed.
If the reader is interested on that,
we refer to~\cite{autostoppisti}.
In this reference, there is a detailed description of the computation of such matrices and 
an some examples where one can check the coefficients via by-hand computation.

\begin{code}
With this respect we further underline 
that Remark 3.3 in~\cite{autostoppisti} is really important from a coding point of view.
Indeed, one can check the correctness of the projection operators via a proper identity,
the so-called ``GBD-identity''.
In \vpp{} this strategy was deeply used and 
once it was verified in the 80\% of the cases the \emph{whole} implementation was correct.
\end{code}

From Equation~\eqref{eqn:aOperatorMatForm} the matrices $\mathbf{G},\,\PiNStarMat$ and $\PiNMat$ are essential to compute the local matrix that represents the grad-grad operator. 
However, if we consider a problem with several bilinear operators,
such matrices can be also used in the computation other problem operators.
Moreover, they may also play a key role in some post processing procedure, like the computation of the error.
As a consequence one needs to be able to compute the local matrix, store it and, eventually, reload it.

In \vpp{} we define a proper data structure that accomplishes all these goals.
Such data structure is called \texttt{vemMatrixHandler}.
Before describing the implementation detail,
we have to spend few words on the ``reloading process''.
Indeed, to reload a matrix it is necessary to recognise it. 
In \vpp{} each matrix is associated with a tag, i.e.,
a \cpp{} enum variable called \texttt{matrixType}.
As a consequence each matrix can be uniquely detected by this tag.

The data structure \texttt{vemMatrixHandler} has a ``compartment'' for each mesh face.
Then, during the loop over faces, 
once a new matrix is computed,
it is stored inside the actual face compartment with is tag.
As a consequence, during the computation of one specific local matrix,
it is possible to access to the compartment of the actual face and
reload the matrixes needed.
To have a clearer idea of the whole procedure
we summarise it in Listing~\ref{alg:loopOnEle}.

{\footnotesize
\begin{lstlisting}[language=C++, caption=Example of loop to assemble a linear system, label=alg:loopOnEle]
std::vector<matrixType> listOfMat(4);
listOfMatrix[0] = GNODALMAT;
listOfMatrix[1] = PINABLASTARMAT;
listOfMatrix[2] = PINABLAMAT;
listOfMatrix[3] = STIFFOPERATOR;

vemMatrixHandler allFaceMats(numFaces);

for(UInt faceId=0; faceId<numFaces; ++faceId)
{
    vemFace2d * actualFace = mesh.getVemFace(faceId);
    for(UInt mId=0; mId<listOfMatrix.size(); ++mId)
    {
        Eigen::MatrixXd tmpMat;
        computeTheMatrix(*actualFace, listOfMat[mId], allFaceMats, tmpMat);
        
        allFaceMats.insertMatrix(*actualFace, listOfMat[mId], tmpMat);
    }
}
\end{lstlisting}}

Notice that the function \texttt{computeTheMatrix}, line 15 of Listing~\ref{alg:loopOnEle},  
takes as input the container itself: 
the current matrix computation may requires other matrices
that are already inside \texttt{vemMatrixHandler}.
From this specific example, this point becomes more evident.
Indeed, the computation of the fourth matrix, \texttt{STIFFOPERATOR}, 
requires \texttt{GNODALMAT}, \texttt{PINABLASTARMAT} and \texttt{PINABLAMAT}, 
see Equation~\eqref{eqn:aOperatorMatForm}.
As a consequence we are forced to compute local matrices in a precise sequence 
otherwise we are not able to compute all of them.

On line 17 of Listing~\ref{alg:loopOnEle} we also observe 
that the local matrix \texttt{tmpMat} is inserted in the container \texttt{allFaceMats} which takes in 
input: the actual face, \texttt{actualFace}, and its tag, \texttt{listOfMatrix[matId]}.
The former information is necessary to put the matrix in the compartment associated with the current face
while the latter one is used to uniquely identify it.


\section{How play with \vpp{}}\label{sec:play}

As already mentioned in the introduction,
the idea behind~\vpp{} is not to give a code that is used as a black box,
but users can play with it developing new features and testing their own research.
The aim of this section is to further underline this aspect and 
give the reader more hint about that.
In Section~\ref{sec:mat} we show how to define a new local matrix in the code.
Another interesting features that a user can study with~\vpp{} is
the impact of different choices of quadrature formulas and solvers.
In Section~\ref{sec:qr} and~\ref{sec:solvers} we will show how to plug in a different 
quadrature rule or solver in~\vpp{}, respectively.
Moreover, inside the library there are ten tutorials whose 
aim is to make the reader understand the code and 
let him plug in new features.

\subsection{New local matrices}\label{sec:mat}
In Section~\ref{sec:EqualVEM} we already saw that 
bilinear and projection operators are build via elementary operation among matrices.
As a consequence, if one has to develop a new operator or projection operator,
it has to define a new local matrix.
Creating a new matrix in~\vpp{} is a straightforward procedure and 
it is deeply described in \texttt{tutorial 5}.

\paragraph{Define the name of the new local matrix.}
Each matrix has a unique tag, i.e., 
a value of the enum \texttt{vemMatrixType}. 
If one would like to create a new matrix, 
it has to define a new tag, for instance \texttt{MYNEWMAT}, and
add it to the \texttt{vemMatrixType} enum, see Listings~\ref{alg:enumMat}.
{\footnotesize
\begin{lstlisting}[language=C++, caption=New matrix definition, label=alg:enumMat]

enum vemMatrixType{
  NULLMATRIXTYPE,
  HNODALMAT,                    
  HVECTNODALMAT,
  ...
  MYNEWMAT,
  ...
}

\end{lstlisting}}

\paragraph{Define how to compute a matrix.}
In~\vpp{} there is a one-to-one correspondence 
between \cpp{} enum \texttt{vemMatrixType} and~\cpp{} classes.
To develop a new local matrix the user has 
to define a public child of the~\vpp{} main class \texttt{vemLocalMatrixDefinition}.
For instance in Listing~\ref{alg:newMatrix},
we are defining a new class \texttt{vem2dMyNewMat} that defines a new matrix whose tag is \texttt{MYNEWMAT}. 
Moreover, it has to define the method \texttt{computeMatrix} of the class 
\texttt{vemLocalMatrixDefinition}.
Such method is virtual so it can be redefined on each child class of \texttt{vemLocalMatrixDefinition}.
{\footnotesize
\begin{lstlisting}[language=C++, caption=New matrix definition, label=alg:newMatrix]

class myNewMatrix : public vemLocalMatrixDefinition
{
	public:

	vem2dMyNewMat(UInt _degreeK)
	  :  vemLocalMatrixDefinition(2, _degreeK, vem::MYNEWMAT) 
	{
	...
	...
	
	void computeMatrix(const vemFace2d  & element, 
	                   vemMatrixHandler & matrixHandler, 
	                   Eigen::MatrixXd & localMatrix) const;

	}
};
\end{lstlisting}}

When you are looping over the mesh elements to compute the local matrices the method 
\texttt{findOrComputeMatrix} of the class \texttt{vemLocalMatrixDefinition} is called, 
see Linsitngs~\ref{alg:findOrCompute}.
Such method first look for the matrix inside the \texttt{vemMatrixHandler}.
Then, if it does not find the matrix inside such container,
it calls the virtual function \texttt{computeMatrix} 
whose implementation is in the child class and,
as a consequence, the desired matrix is build
and then put in the container.

{\footnotesize
\begin{lstlisting}[language=C++, caption=\texttt{findOrComputeMatrix} method developed in \texttt{vemLocalMatrixDefinition}, label=alg:findOrCompute]
 class vemLocalMatrixDefinition
{
    ...
    ...
    template<class ELEMENT>
    void findOrComputeLocalMatrix(const ELEMENT & element, 
                                 vemMatrixHandler & matHandler, 
                                 Eigen::MatrixXd & localMatrix) const
    {
        if(!matHandler.getLocalMatrix(element, matrixName, localMatrix))
        {
            computeMatrix(element, matrixHandler, localMatrix);
            matHandler.insertMatrix(element, matrixName, localMatrix);
        }
    }
    ...
    private:
    
    vemMatrixType matrixName;
};
\end{lstlisting}}

\paragraph{Add such matrix in the \texttt{localAssembler}.}
The class \texttt{localAssembler} manages the construction of matrices in \vpp{}. 
As a consequence, once a new matrix is added, 
the method that computes the matrices in \texttt{localAssembler} has to be properly updated, see Listings~\ref{alg:newMatrixLocalAss}.
{\footnotesize
\begin{lstlisting}[language=C++, caption=Method of the \texttt{localAssembler} updated, label=alg:newMatrixLocalAss]

vemLocalMatrixDefinition * localMatrixAssembler::makeVemLocalMatrix(
                UInt localMatrixId, 
                const geometricalEntityType & entityType) const
{
    ...
    
    vemLocalMatrixDefinition * result = NULL;
    switch(matrixName)
    {
        case(HNODALMAT):
            if(entityType==GEOENTITY2D) result = new vem2dLocalMatrixH(degreeK);
            else                        result = new vem3dLocalMatrixH(degreeK);
            break;
        ...
        case(MYNEWMAT):
            if(entityType==GEOENTITY2D) result = new vem2dMyNewMat(degreeK);
            else                        result = NULL;
            break;
        ...
        
}
\end{lstlisting}}

These three steps are required 
to build a new matrix inside \vpp{}. 
In \texttt{tutorial 5} and \texttt{tutorial 7} there are a more precise description of the procedure with a specific reference to the lines of the \texttt{localAssembler} to be modified.

\subsection{Quadrature rules}\label{sec:qr}

As already mentioned before, 
one of the challenges in a virtual element code is 
dealing with the integration over polygons and inside polyhedrons.
Independently on the shape of the domain,
an integration formula is defined via a set of integration points, $\{\xx_i\}_{i=1}^n$, 
and weights, $\{\omega_i\}_{i=1}^n$, i.e.,
\begin{equation}
\int_P f(\xx)\dx \approx \sum_{i=1}^n \omega_i f(\xx_i)\,.
\label{eqn:intDef}
\end{equation}
In \vpp{} once $\{\xx_i,\,\omega_i\}_{i=1}^n$ are defined, 
they are stored inside polyhedrons, polygons or edges.
Then, in \vpp{} there is a specific \cpp{} class 
that will automatically use such list of quadrature points and weights 
to compute an integral.

For edges in~\vpp{} there are already implemented both Gau\ss{}-Lobatto and Gau\ss{} quadrature formulas.
For polygons a quadrature rule based on sub-triangulation and Vianello's quadrature rule~\cite{vianello}.
Finally for polyhedrons it is implemented one based on sub-tetralisation and
an extension of Vianello's formula to polyhedrons~\cite{vianello3d}.

Moreover, it is also implemented a general rule to make a compression of them~\cite{compression}.
Starting from the set of quadrature points such compression algorithm solve a non negative least squares problem that also ensures that the weights are positive.

Since the integration over polytopes is still an open research field,
\vpp{} was designed in such a way that adding a new quadrature rule is straightforward.
Indeed, if one would like to introduce a new quadrature rule,
it has to follow these steps 
which are similar to the ones described in the case of adding a new matrix.

\paragraph{Define the name of the new quadrature formula.}
In \vpp{} there is an enum used to identify the quadrature points called \texttt{quadraturePointsType}. 
As a consequence, if you would like to add a new quadrature formula,
you have to give it a name,
for instance \texttt{MYQUADFORMULA4PPOLYGONS},
and add in such enum, see Linstings~\ref{alg:quadEnum}.
{\footnotesize
\begin{lstlisting}[language=C++, caption=Enum \texttt{quadraturePointsType} to be updated, label=alg:quadEnum]
enum quadraturePointsType{NONERULE, 
	                        GAUSSEDGEQUADRATURE, 
	                        GAUSSLOBATTOEDGEQUADRATURE, 
	                        STANDARDTRIANGLEQUADRATURE, 
	                        STANDARDTETRAHEDRONQUADRATURE, 
	                        MYQUADFORMULA4PPOLYGONS,
	                        VIANELLO4POLYGONS,
	                        VIANELLO4POLYHEDRONS};
\end{lstlisting}}

\paragraph{Define the rule} 
In \vpp{} there are three main classes to build quadrature rules according to the dimension of the integration domain.
Indeed we have
\begin{itemize}
 \item \texttt{edgeQuadraturePointsMaker}: to compute quadrature rule for edges;
 \item \texttt{faceQuadraturePointsMaker}: to compute quadrature rule for polygons;
 \item \texttt{volumeQuadraturePointsMaker}: to compute quadrature rule for polyhedrons.
\end{itemize}
Each of this class has a method to select which quadrature rule put inside the element.
For instance in Listing~\ref{alg:quadRule}, we show the method \texttt{} 
to compute the quadrature rule, i.e., the set $\{\xx_i,\,\omega_i\}_{i=1}^n$ in the case of polygons.
{\footnotesize
\begin{lstlisting}[language=C++, caption=Method to define polygon quadrature rule, label=alg:quadRule]
...
switch(typeQuad)
{
    case(STANDARDTRIANGLEQUADRATURE):
         makeQuadraturePointsAndWeightsViaTriangulation(ele, deg, pts, weights);
         break;
         
    case(VIANELLO4POLYGONS):
         makeQuadraturePointsAndWeightsVianello(ele, deg, pts, weights);
         break;
}
...
\end{lstlisting}}
As a consequence, if one would like to add a new quadrature,
it has to enrich such method with the case associated with its quadrature rule, see Listing~\ref{alg:new}.
Such an enrichment can be done adding a proper method in the class \texttt{faceQuadraturePointsMaker}, 
creating an external class within \vpp{} library or 
via a wrapper of an external library.
{\footnotesize
\begin{lstlisting}[language=C++, caption=Method updated with a new rule, label=alg:new]
...
switch(typeQuad)
{
    case(STANDARDTRIANGLEQUADRATURE):
         makeQuadraturePointsAndWeightsViaTriangulation(ele, deg, pts, weights);
         break;
         
    case(VIANELLO4POLYGONS):
         makeQuadraturePointsAndWeightsVianello(ele, deg, pts, weights);
         break;
    
    case(NEWQUADRULEFORPOLYGONS):
        makeMyNewQuadraturePointsAndWeights(ele, deg, pts, weights);
        break;
}
...
\end{lstlisting}}

Moreover, inside such class there it is also implemented the compression routine proposed in~\cite{compression},
see the function \texttt{compressPoints} in the class \texttt{faceQuadraturePointsMaker}. 
As a consequence, once a set $\{\xx_i,\,\omega_i\}_{i=1}^n$ is build it can be always compressed, i.e., 
the number of quadrature points can be reduced by properly changing the weights.

\subsection{Solvers}\label{sec:solvers}

Another active research field in the framework of partial differential equations is 
finding new ways to solve a linear system arising from the PDE discretisation.
\vpp{} is designed to facilitate as more as possible researchers interested in such issue.

First of all it has already an interface to the most common linear algebra \cpp{} libraries. 
More specifically, there are two wrapper classes \texttt{eigenSolver} and \texttt{petscSolver} 
for the \cpp{} libraries Eigen~\cite{eigen} and Petsc~\cite{petsc}, respectively.
As a consequence, if one developed a particular solver inside such a libraries,
it can look inside such wrappers and call all the routines needed. 

Then, since a method to solve a linear system may require additional information about the mesh or the degrees of freedom,
in~\vpp{ } there is a class \texttt{globalSparseSolverData} that stores some useful information
and they are \emph{directly} passed to the solver, see Listing~\ref{alg:dataSolver}.

Other than the basic information such as the sparse matrix, \texttt{matrixSparse},
the right hand side, \texttt{rhs}, and the solution vector,   \texttt{solution},
This class contains some additional information that can be filled to set particular solvers.

For instance, in \vpp{} preconditioners based on field split approach~\cite{fieldsplit} were developed
for a Maxwell equations~\cite{BEIRAODAVEIGA202282},
the variable \texttt{listOfFieldData} contains all the details about electric and magnetic fields 
to set the preconditioner. 
Moreover, when a BDDC preconditioner is set, 
the method requires the coordinates of pressure dofs~\cite{BDDCStokes} and 
such list of points is stored inside the variable \texttt{listoOfDofData}.
{\footnotesize
\begin{lstlisting}[language=C++, caption=Class that contains all the information on solvers, label=alg:dataSolver]
class globalSparseSolverData
{
    public:

    solverOption solverOpt;
    Eigen::SparseMatrix<REAL> * matrixSparse;
    Eigen::SparseMatrix<REAL> * massMatrixSparse;
    Eigen::VectorXd  * rhs;
    Eigen::VectorXd  * solution;
    bool symmetricPattern;
    UInt numIterationDone;
    REAL lastResidualError;
    REAL timeToSolve;

    std::vector<fieldData> listOfFieldData;
    dofData listOfDofData;

    meshData listOfMeshData;
    
    globalSparseSolverData();
...
\end{lstlisting}}

\subsubsection{Adding a new solver library}

The process of adding a new library that solve a linear system is technical but not difficult.
The factory class \texttt{solverHandler} is able to distinguish among solvers from different libraries, 
see Linstings~\ref{alg:solver}.
Then, the methods \texttt{solveWithEigen} and \texttt{solveWithPetsc} call the solve method from the 
wrapper classes of the \cpp{} libraries Eigen and Petsc, respectively.

{\footnotesize
\begin{lstlisting}[language=C++, caption=Solve method of the class \texttt{solverHandler}, label=alg:solver]
void solverHandler::solve(globalSparseSolverData & problemDescription) const 
{
    setSymmetry(problemDescription);
    printInitialSystemData(problemDescription);
    
    switch(problemDescription.getSolverLibrary())
    {
        case(EIGENLIBRARY):
                           solveWithEigen(problemDescription);
                           break;
        case(PETSCLIBRARY):
                           solveWithPetsc(problemDescription);
                           break;
        default:
                           std::cout << "!! !! !! Unknown solver" << std::endl;
                           exit(1);
    }
    printFinalSystemData(problemDescription);
}
\end{lstlisting}
}

If one would like to add another solver library, 
it can get ideas on the implementation of the wrappers for \texttt{Eigen} or \texttt{Petsc},
but in general it has to make the following steps.

\paragraph{Add linking}
First of all a proper link to the desired solver library is needed.
To achieve this goal, one has to properly modify the \texttt{CMakeList.txt}.

\paragraph{Define a wrapper class}
At this step one has to create a proper wrapper class that starting from the information provided by \texttt{globalSparseSolverData}
solve the linear system and store the solution inside the \texttt{Eigen::VectorXd solution}, see Listings~\ref{alg:dataSolver}.
To build a suitable wrapper,
one can follow the implementation of \texttt{eigenSolver} and \texttt{petscSolver} classes.

\paragraph{Update \texttt{solverEnum}}
Inside the header file \texttt{solverEnum.h} there are two \cpp{} enums to be updated with the tag associated with the new library, 
\texttt{solverType} and \texttt{solverLibrary}. 
Then, there are the methods  \texttt{translateToSolverType}, \texttt{getLibraryFromSolver} and
\texttt{getKindOfMethodFromSolver}. that has to be updated with the new enum added in \texttt{solverType}.

\paragraph{Update \texttt{solverHandler}}
The last step is to add the new class in the \texttt{solverHandler}.
More specifically, one has to add the method \texttt{solve} a call to the wrapper class,
compare Listing~\ref{alg:solver} and Listing~\ref{alg:newSolver}.

{\footnotesize
\begin{lstlisting}[language=C++, caption=Solve method of the class \texttt{solverHandler}, label=alg:newSolver]
void solverHandler::solve(globalSparseSolverData & problemDescription) const 
{
    setSymmetry(problemDescription);
    printInitialSystemData(problemDescription);
    
    switch(problemDescription.getSolverLibrary())
    {
        case(EIGENLIBRARY):
                           solveWithEigen(problemDescription);
                           break;
        case(PETSCLIBRARY):
                           solveWithPetsc(problemDescription);
                           break;
        case(NEWSOLVERLIBRARY):
                           solveWithNewSolverLibrary(problemDescription);
                           break;
        default:
                           std::cout << "!! !! !! Unknown solver" << std::endl;
                           exit(1);
    }
    printFinalSystemData(problemDescription);
}
\end{lstlisting}
}


\section{Partial Differential Equations in \vpp{}}\label{sec:exe}

In this section we will show some types of partial differential equation handled by \vpp{}.
More specifically, in Section~\ref{sec:acc} we provide some standard academical problems,
while in Section~\ref{sec:prac} we see some more practical applications 
where \vpp{} was successfully used.
For each example we furnish the reference where it was taken 
so that if the reader is interested in such topic,
it will have more detail.

\subsection{Academic Applications}\label{sec:acc}

\subsubsection*{Laplacian problems}
The first step in the resolution of partial differential equation is the resolution of a Laplacian problem, i.e.,
\begin{equation}
\left\{
\begin{array}{rll}
-\Delta u=&f     &\text{in}\:\Omega\\
u=&g             &\text{on}\:\partial\Omega 
\end{array}
\right.\,.
\label{eqn:poiss2}
\end{equation}
In~\vpp{} there it is possible to solve such a problem both in two and three dimensional case and for a generic approximation degree $k$. 
In~\cite{da2017high} it is possible to see all the theory and the results related to such an equation.
Moreover, in the library there is the folder \texttt{bin/ApolloTest} 
that contains all the executable associated with each example of~\cite{da2017high}.

\subsubsection*{Stokes and Navier-Stokes problems}
The~\vpp{} library contains also a module related to the resolution of Stokes 
\begin{equation}
\left\{\begin{array}{rll}
       {\bm \Delta}\uu - \nabla p &= \textbf{f} &\text{in}\:\Omega \\
       \div(\uu)&=0 &\text{in}\:\Omega\\
       \uu &= {\bm g}&\text{on}\:\partial\Omega
       \end{array}
\right.
\label{eqn:stokes}
\end{equation}
and the Navier-Stokes problem
\begin{equation}
\left\{\begin{array}{rll}
       {\bm \Delta}\uu + \uu \nabla \uu - \nabla p &= \textbf{f} &\text{in}\:\Omega \\
       \div(\uu)&=0 &\text{in}\:\Omega\\
       \uu &= {\bm g}&\text{on}\:\partial\Omega
       \end{array}
\right.\,.
\label{eqn:ns}
\end{equation}
Using the virtual element method to solve such equation is particularly appealing.
Indeed, in the VEM it is possible to construct a discrete vector field $u_h$ whose divergence is pointwise zero.
This a property is not common to the majority of the finite element approaches and 
it has a lot of advantages such as the decupling of the error on velocities and pressure~\cite{vaccaDivFree}.
We refer the reader to~\cite{brick, beirao2020stokes} to have a deeper analysis on the results obtained by~\vpp{} for both the two and three dimensional case.
Moreover, if on is interested in the implementation of the code we refer to the folders \texttt{bin/Stokes2dPaperTest} and \texttt{bin/Stokes3dPaperTest}
where there are the executables associated with the examples associated with the above mentioned papers.

\subsubsection*{Static Maxwell Equation in Kichuci formulation}
The virtual element method was successfully applied to solve a magnetostatic problem
\begin{equation}
\left\{
\begin{array}{rll}
\bcurl(\HH)&=\jj &\text{in}\:\Omega \\
\div(\mu\HH)&=0   &\text{in}\:\Omega \\
\HH\wedge\nn&={\bm g} &\text{on}\:\partial\Omega . 
\end{array}
\right.
\end{equation}
in the so-called Kikuchi formulation, see, e.g.,~\cite{KikuchiIEEE}.
More specifically in~\cite{max2d} VEM was successfully applied to general order two dimensional magnetostatic problems.
In~\cite{lowestMax3d} a lowest order degree for the three dimensional case was defined,
while in~\cite{max3d} the general order case for three dimension was introduced.
In \vpp{} there are all the executables for such papers.
More specifically the ones of~\cite{max2d} are in \texttt{bin/maxwell2d},
while the examples of~\cite{lowestMax3d} and~\cite{max3d} are all in \texttt{bin/maxwell3dPaperTest}.

\subsubsection*{Elasticity Equation with the Hellinger-Ressiner principle}

Another intersting problem solved via the virtual element method and implemented inside \vpp{} is 
the linear elasticity problems based on the Hellinger–Reissner variational principle.
This mixed formulation describes a linear elasticity via displacement and the stress fields.
Indeed, we solve the following partial differential equation
\begin{equation}
\left\{\begin{array}{rll}
       -\div({\bm{\sigma}}) &= \textbf{f} &\text{in}\:\Omega \\
       {\bm{\sigma}} &= \mathbb{C}\epsilon(\uu) &\text{in}\:\Omega\\
       \uu &= {\bm g}&\text{on}\:\Gamma_1\\
       \bm{\sigma}\nn &= {\bm \psi}&\text{on}\:\Gamma_2
       \end{array}
\right.\,.
\label{eqn:hr}
\end{equation}
where $\uu$ represents the displacement field and ${\bm{\sigma}}$ is the stress tensor,
$\Gamma_1$ and $\Gamma_2$ are two partitions of the boundary of the domain $\Omega$,
where the displacement and the normal component of the stress field are fixed.
In a finite element framework one of the main issue related to such a problem is obtaining a stable and accurate 
scheme that does preserve both the symmetry of the stress tensor and the continuity of the tractions at the inter-elements.
However, we are able to avoid these two drawbacks by exploiting the flexibility of the VEM.
For more detail about the mixed linear elasticity based on the Hellinger–Reissner variational principle
we refer the reader to~\cite{DASSI2020112910}.
Moreover, the executable of the tests shown in such paper are available in the folder \texttt{bin/HellingerReissnerPaper}.

\subsection{Practical Applications}\label{sec:prac}

\subsubsection*{Elasticity with curved edges}
In~\cite{curviMec} \vpp{} was successfully used in two dimensional solid mechanics applications
that are characterised by curve boundaries.
More specifically they test the virtual element method on 
two classical benchmark examples involving inelastic materials:
the thick-walled viscoelastic cylinder subjected to internal pressure~\cite{Zienkiewicz_Taylor_Fox} and 
the perforated plastic plate~\cite{zienkiewicz2005finite}. 

\begin{figure}[!htb]
\centering
\begin{tabular}{cc}
\includegraphics[width=0.40\textwidth]{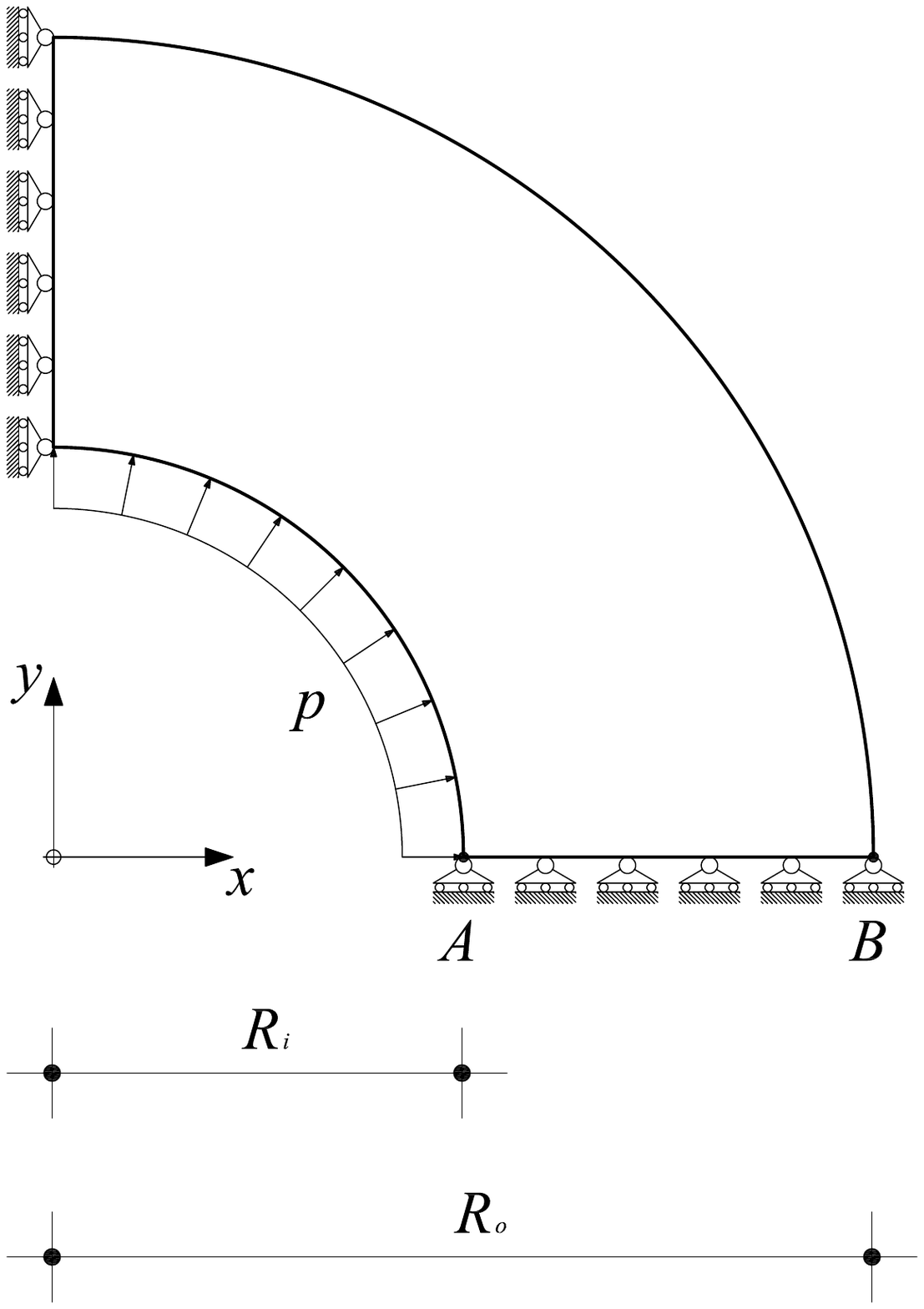}&
\includegraphics[width=0.40\textwidth]{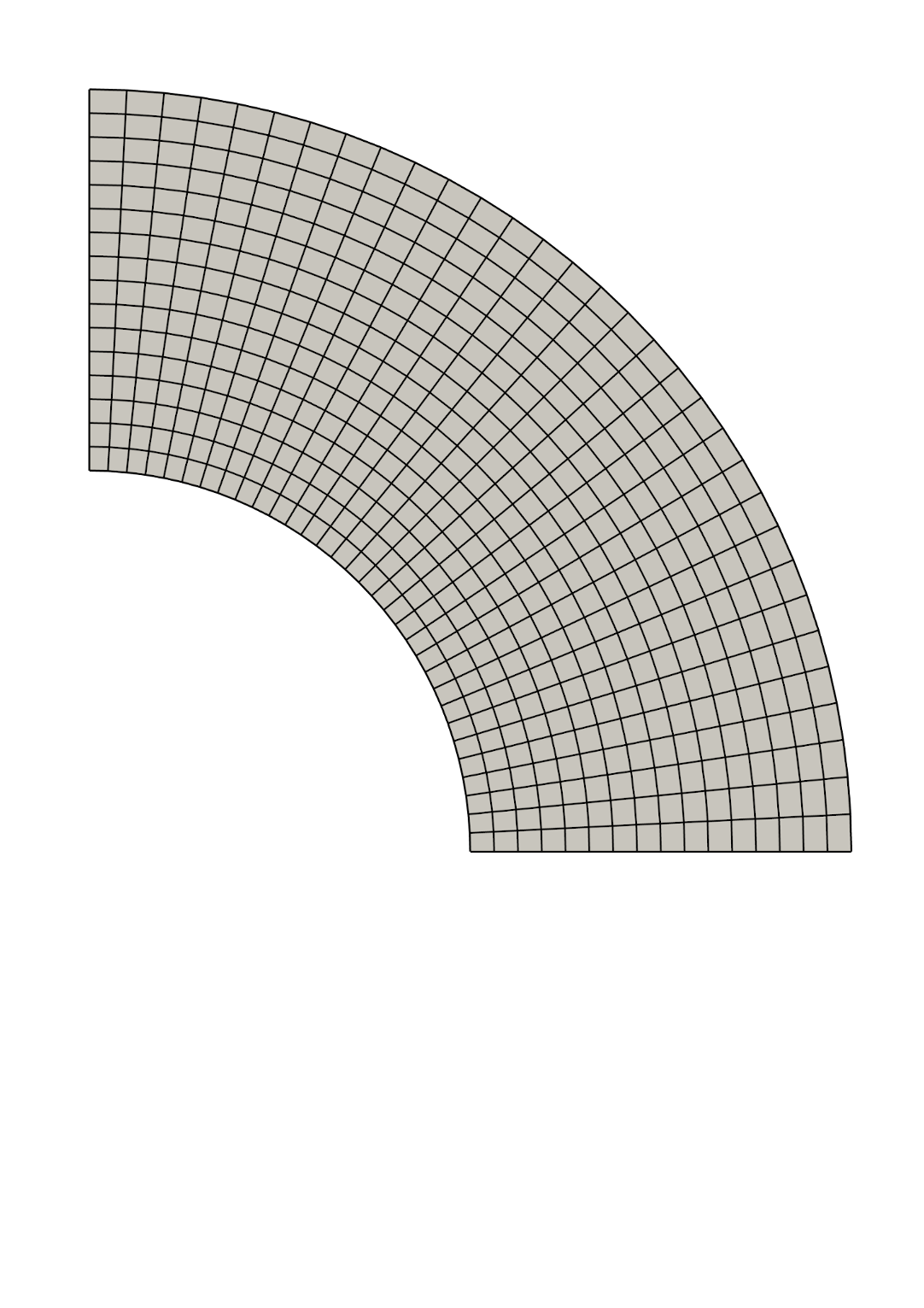}
\end{tabular}
\caption{Thick-walled viscoelastic cylinder subjected to internal pressure. Geometry, boundary conditions, applied load (left), mesh used (right).}
\label{fig:pressure_cylinder_geom}
\end{figure}

\paragraph{The thick-walled viscoelastic cylinder subjected to internal pressure}
The cylinder has inner [resp. outer] radius $R_i = 2$ [$R_o=4$] and 
it is subjected to uniform pressure $p$ on the inner surface, 
see Figure~\ref{fig:pressure_cylinder_geom}. 

The material is isotropic and it follows a viscoelastic constitutive law as outlined 
Section~3.1 of reference~\cite{PartII}. 
We consider two sets of viscoelastic parameters 
$\left( \mu_0 , \mu_1 \right)_{\veOne} = \left( 0.01, 0.99 \right )$ and
$\left (\mu_0 , \mu_1 \right)_{\veTwo} = \left( 0.3, 0.7 \right )$, respectively. 
The former case is devised in a way that the bulk-to-shear moduli ratio for instantaneous loading is given by $K/G(0) = 2.167$ and for long time loading, say at $t=8$, by $K/G(8) = 216.7$, which indicates a nearly incompressible behavior for sustained loading (for instance, at $t = \infty$ the Poisson ratio results $0.498$). The second material set indicates an intermediate response at $t = \infty$ after loading is  applied. Plane strain assumption is applied in this numerical simulation.

In this test case we make a displacement accuracy analysis between the curved VEM formulation with an approximation degree $k=1,2,3$ and a reference solution.
Such reference solution is obtained with quadratic quadrilateral displacement based finite elements with nine nodes $Q9$ 
running in the FEAP platform using an extremely fine mesh~\cite{Zienkiewicz_Taylor_Fox}.

\begin{figure}
\begin{center}
\begin{tabular}{cc}
\includegraphics[width=0.48\textwidth]{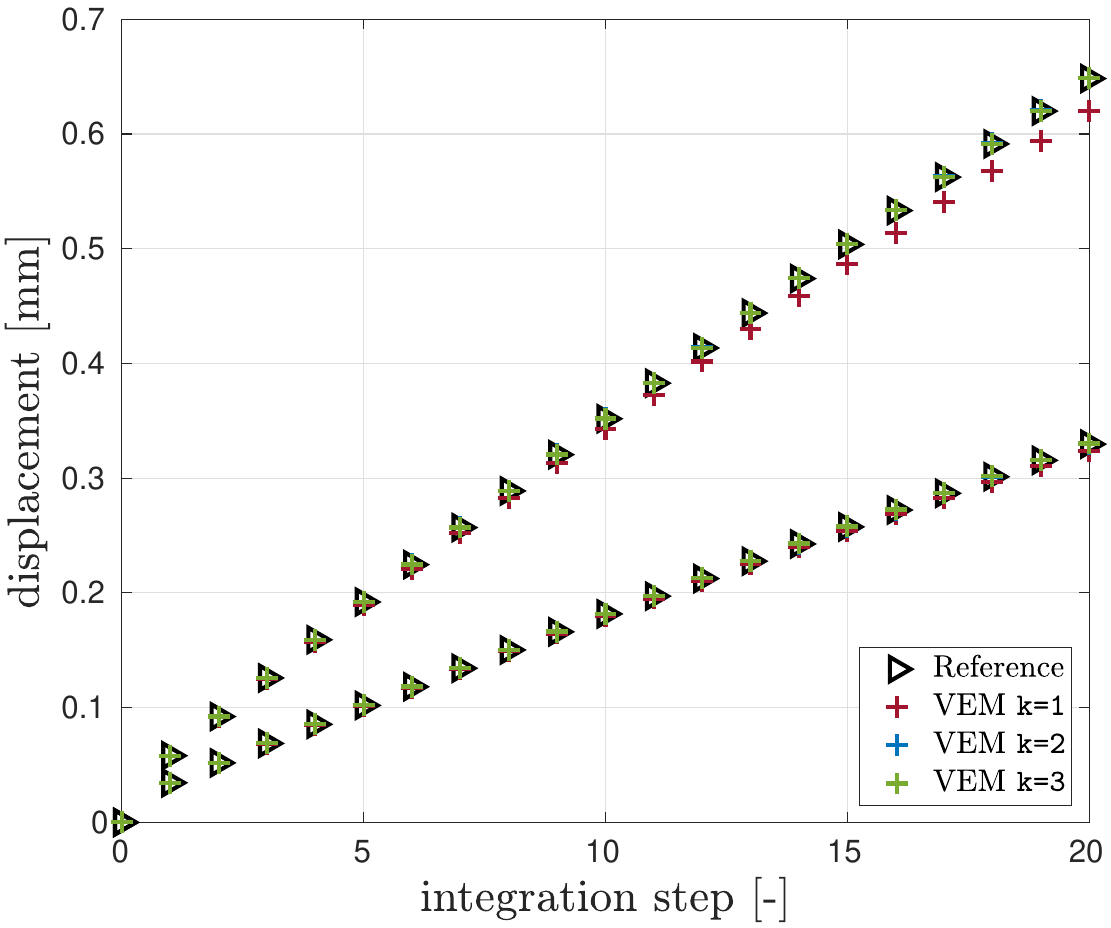} &
\includegraphics[width=0.48\textwidth]{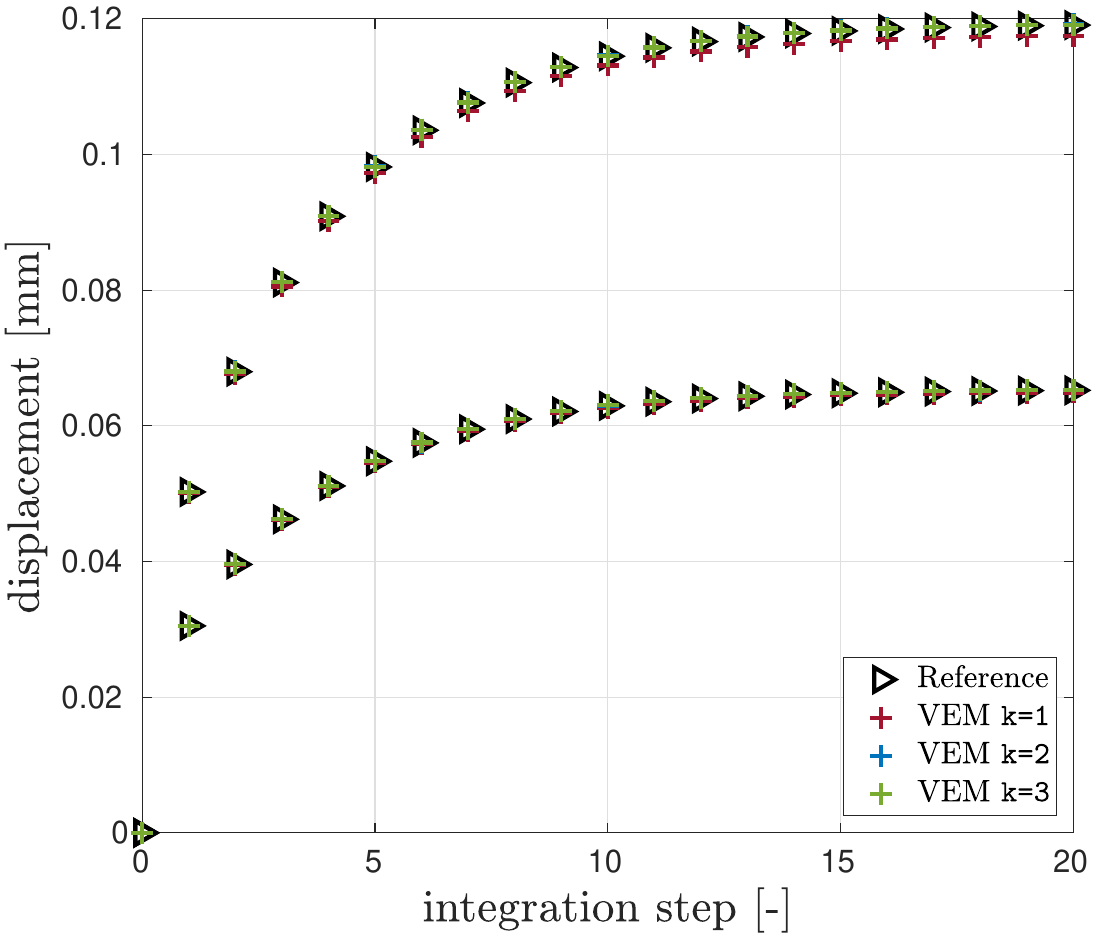}\\
(a) & (b)
\end{tabular}
\end{center}
\caption{Thick-walled viscoelastic cylinder with internal pressure. Integration step {\it vs.} radial displacement curves for control points $A$ (higher curve), and $B$ (lower curve). (a) case $\left( \mu_0 , \mu_1 \right)_{\veOne} = \left( 0.01, 0.99 \right )$; (b) case $\left (\mu_0 , \mu_1 \right)_{\veTwo} = \left( 0.3, 0.7 \right )$.}
\label{fig:visco_disp}
\end{figure}

The integration-step versus displacement curves for control points $A$ and $B$, Figure~\ref{fig:pressure_cylinder_geom}, is shown in Figure~\ref{fig:visco_disp}. 
For both configurations we observe that the solutions provided by \vpp{} match the reference one although 
it decidedly has fewer degrees of freedom. 

\begin{figure}[!htb]
\centering
\begin{tabular}{cc}
\includegraphics[bb=50 0 600 825, clip, angle=0, scale=0.3]{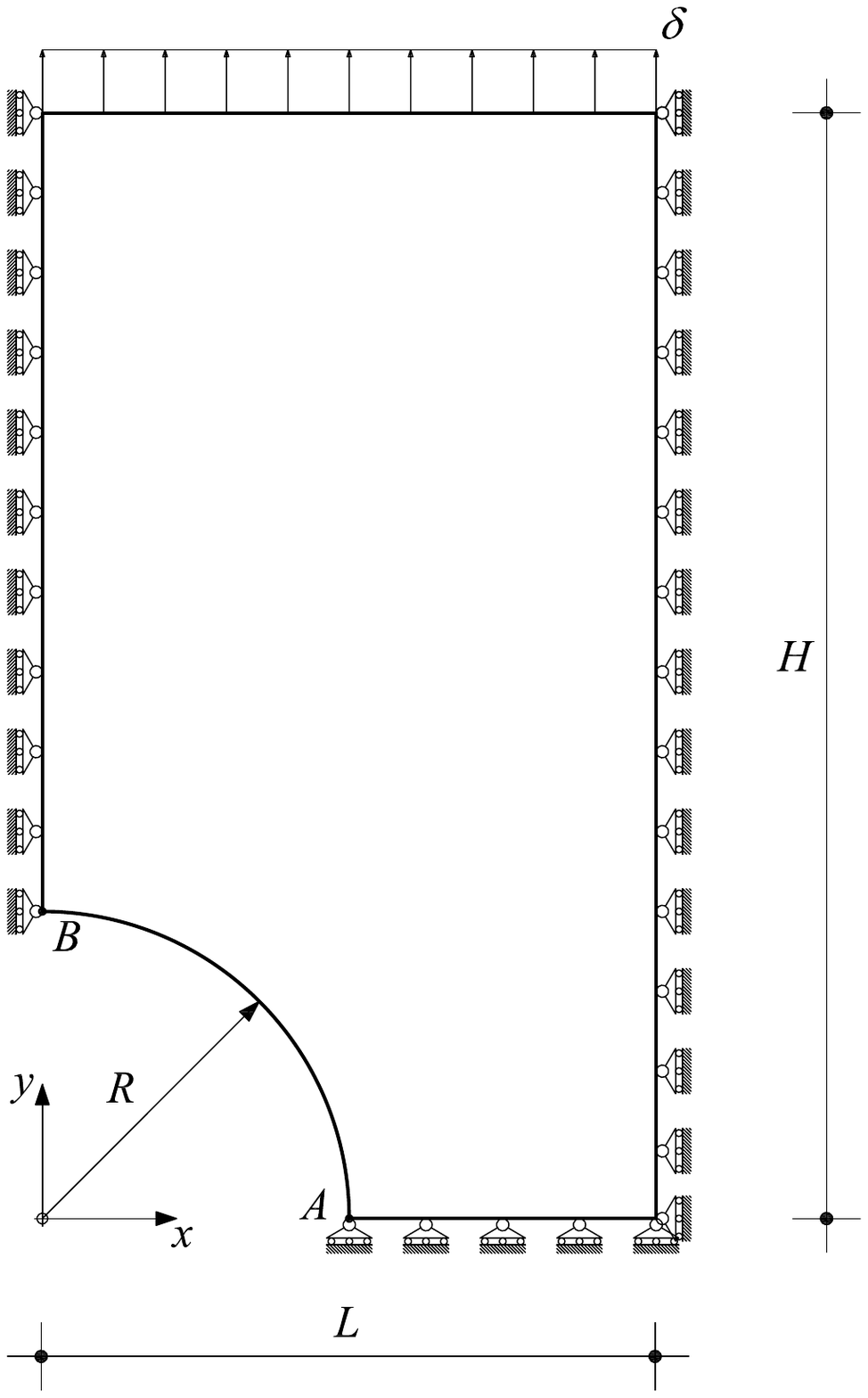} & 
\includegraphics[scale=0.3]{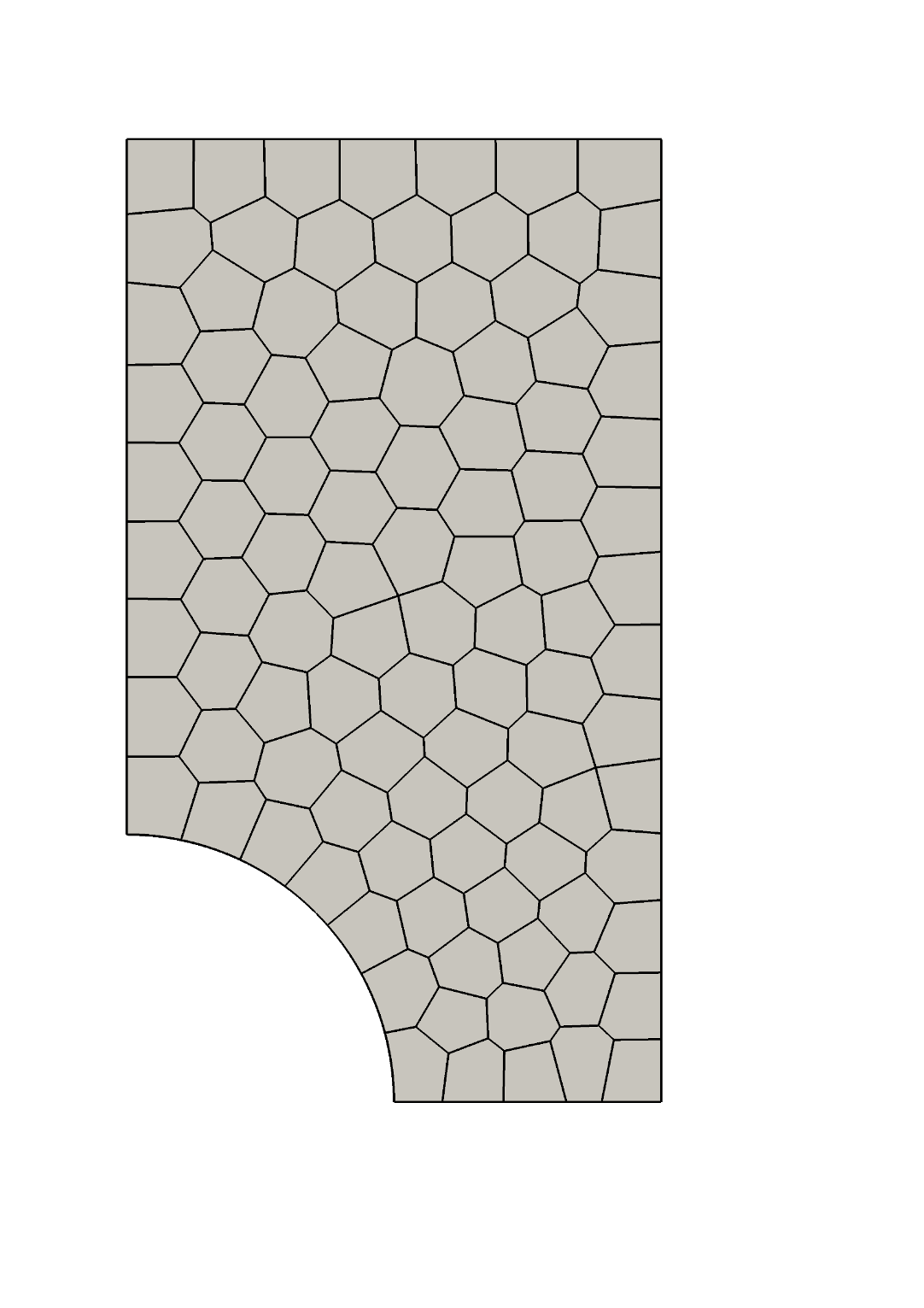}
\end{tabular}
\caption{Perforated plastic plate. Geometry, boundary conditions, loading (left), mesh used (right).}
\label{fig:strip_wt_hole}
\end{figure}

\paragraph{The perforated plastic plate}
In this numerical experiment we consider a rectangular strip with width of $2L = 200$ $\textrm{mm}$ width and length of $2H = 360$ $\textrm{mm}$ with a circular hole of $2R = 100$ $\textrm{mm}$ diameter in its center, see Figure~\ref{fig:strip_wt_hole}). 
Material response here follows classical von Mises plastic constitutive model, 
with material parameters: $E = 7000$ $\textrm{kg} / \textrm{mm}^2$, $\nu = 0.3$,  and yield stress $\sigma_{\textrm{y},0} = 24.3$ $\textrm{kg} / \textrm{mm}^2$~\cite{Zienkiewicz_Taylor_Zhu13}. 
We assume plane strain and we use a standard backward Euler scheme with 
return map projection stress and material moduli computation at the quadrature point level~\cite{simo_computational_1998}. 
Displacement boundary restraints are prescribed for normal components on symmetry boundaries and on top and lateral boundaries. 
Loading is applied by a uniform normal displacement $\delta = 2$ $\textrm{mm}$ with $400$ equal increments on the upper edge, 
see Figure~\ref{fig:strip_wt_hole}. 

\begin{figure}
\begin{center}
\includegraphics[width=0.43\textwidth]{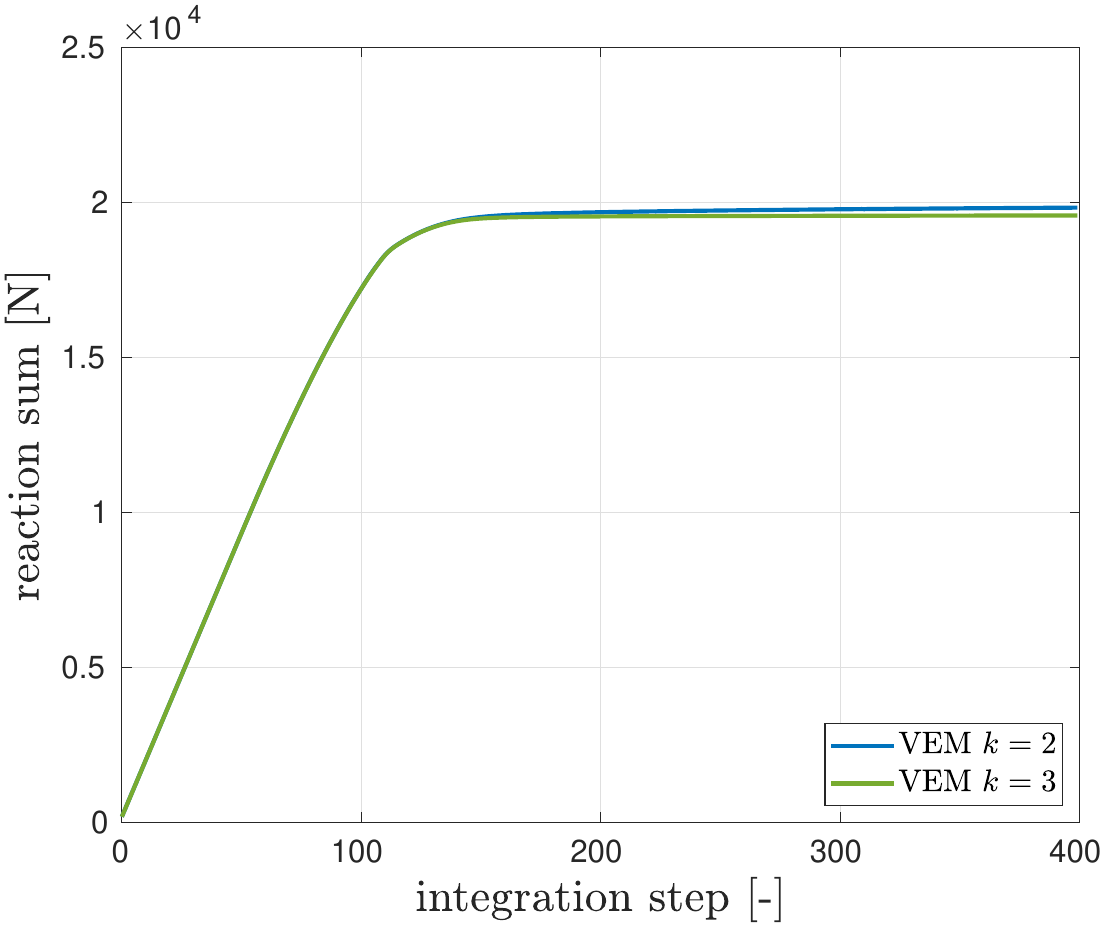}\\
\end{center}
\caption{Perforated plastic plate. Structural response computed by \vpp{}.}
\label{fig:strip_resp}
\end{figure}

We proceed with a virtual element simulation via \vpp{} and we consider order $k=2,3$ in the curved variant.
Accuracy and robustness is assessed by plotting the structural response of the structure, 
in terms of force reaction sum at the imposed displacement top edge {\it vs.} integration step, 
which seems correct for all compared methods and mesh types, 
showing no significant spurious locking phenomena, see Figure~\ref{fig:strip_resp}. 

\subsubsection*{Magnetostatic}

In this subsection \vpp{} is tested in applicative magnetostatic problems. 
We consider a series of numerical experiments and 
we compare the results obtained via the VEM and a standard FEM.
To achieve this goal we take the results provided by MagNet code as a reference solutions~\cite{magnetSoft}.

In \vpp{} both Kikuchi and the potential formulations are implemented.
Then, in the following examples we show the result for each of these approaches and  
we refer to them as VEM Kik and VEM pot.

\paragraph{C-core actuator}
In this section we present a numerical simulation of a C-core electromagnet made in \vpp{}.
In this example the electromagnet is composed of a fixed C-shaped core and a movable plunger. 
A DC current of 1 A supplies the winding around the core limb which in turn excites the magnetic field lines. 
The overall size of the electromagnet is 80 mm$\times\,$60 mm, while the cross-sectional area of the winding,
which incorporates 1000 turns, is equal to 400 mm$^2$.

We focus on the computation of the following physical quantities
\begin{enumerate}
    \item $B_x$: the $x$-component of the magnetic field $\B$ in the mid-point of the air gap;
    \item $F_x$: the $x$-component of the force $\mathbf{F}$ acting on the plunger, via Maxwell’s stress tensor method,
\end{enumerate}
for different position of the plunger with respect to the C-core.
Since we are considering small air-gap widths,
a non-linear approximation of the magnetic permeability is required.
In the following experiments we define the material of both C-core and plunger as 
a standard laminated iron with 5 mm thickness~\cite{bianchi2017electrical}.

\begin{figure}[!htb]
    \centering
    \includegraphics[width=0.5\textwidth]{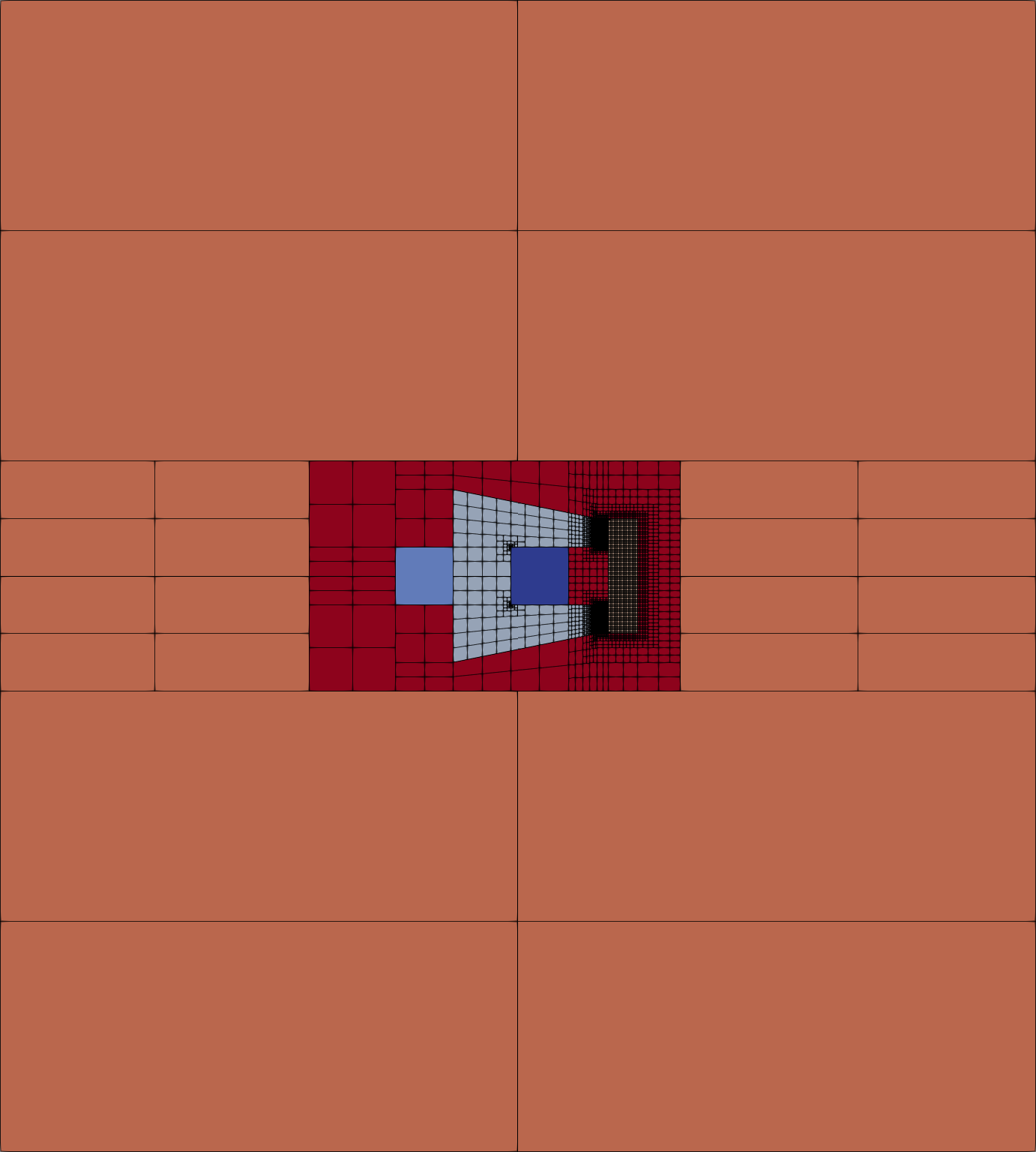}
    \caption{C-core actuator: the whole mesh}
    \label{fig:meshCCore}
\end{figure}

Since the physics behind such problem is well known, 
we can \emph{a priori} identify regions 
which are particularly interesting from the physical point of view.
As a consequence we can generate a mesh refined 
in the ``hot'' regions and coarse in the other ones.
Such procedure is particularly straightforward with VEM 
exploiting as much as possible the possibility to add hanging-nodes. 
Starting from a really coarse mesh,
it was possible to generate the mesh represented in Figure~\ref{fig:meshCCore}.

\begin{figure}[!htb]
    \centering
    \begin{tabular}{cc}
    \includegraphics[width=0.45\textwidth]{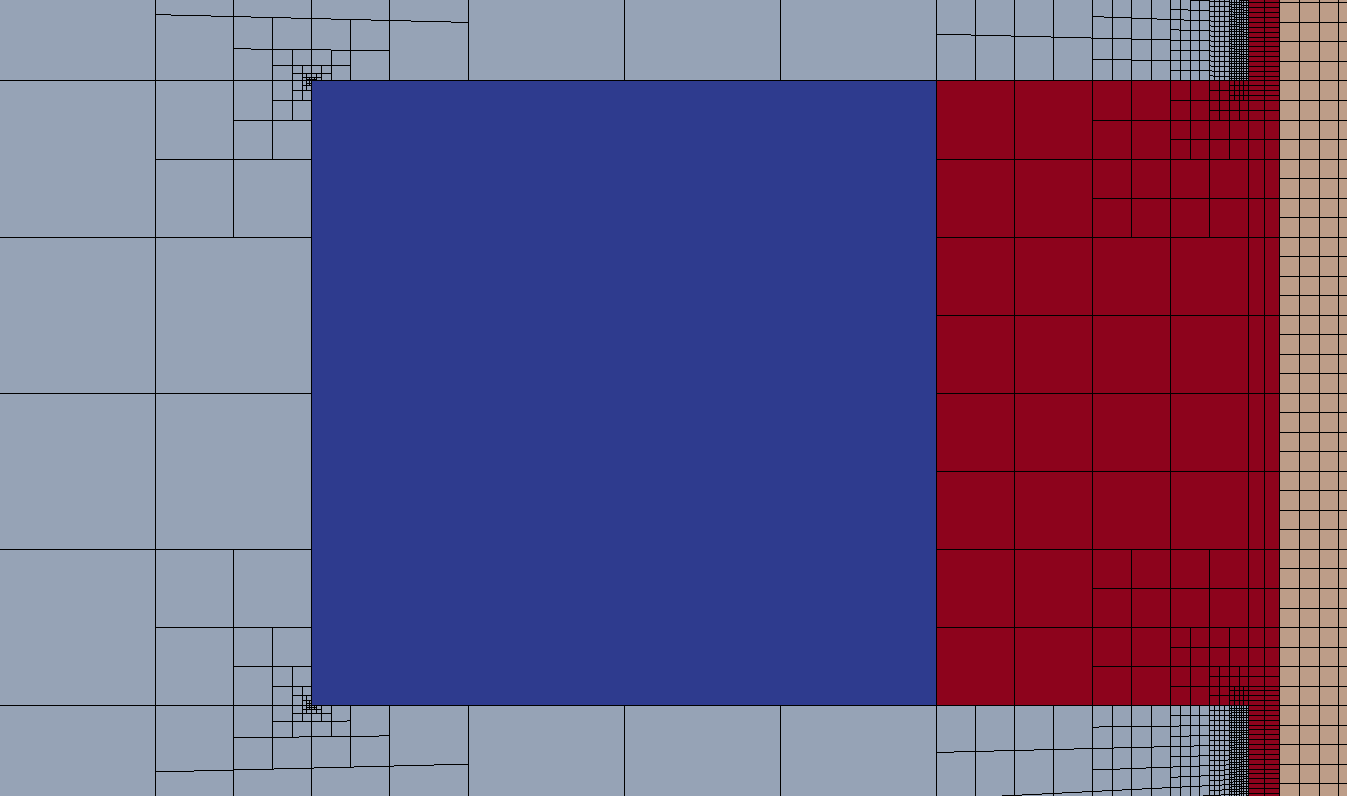} &
    \includegraphics[width=0.45\textwidth]{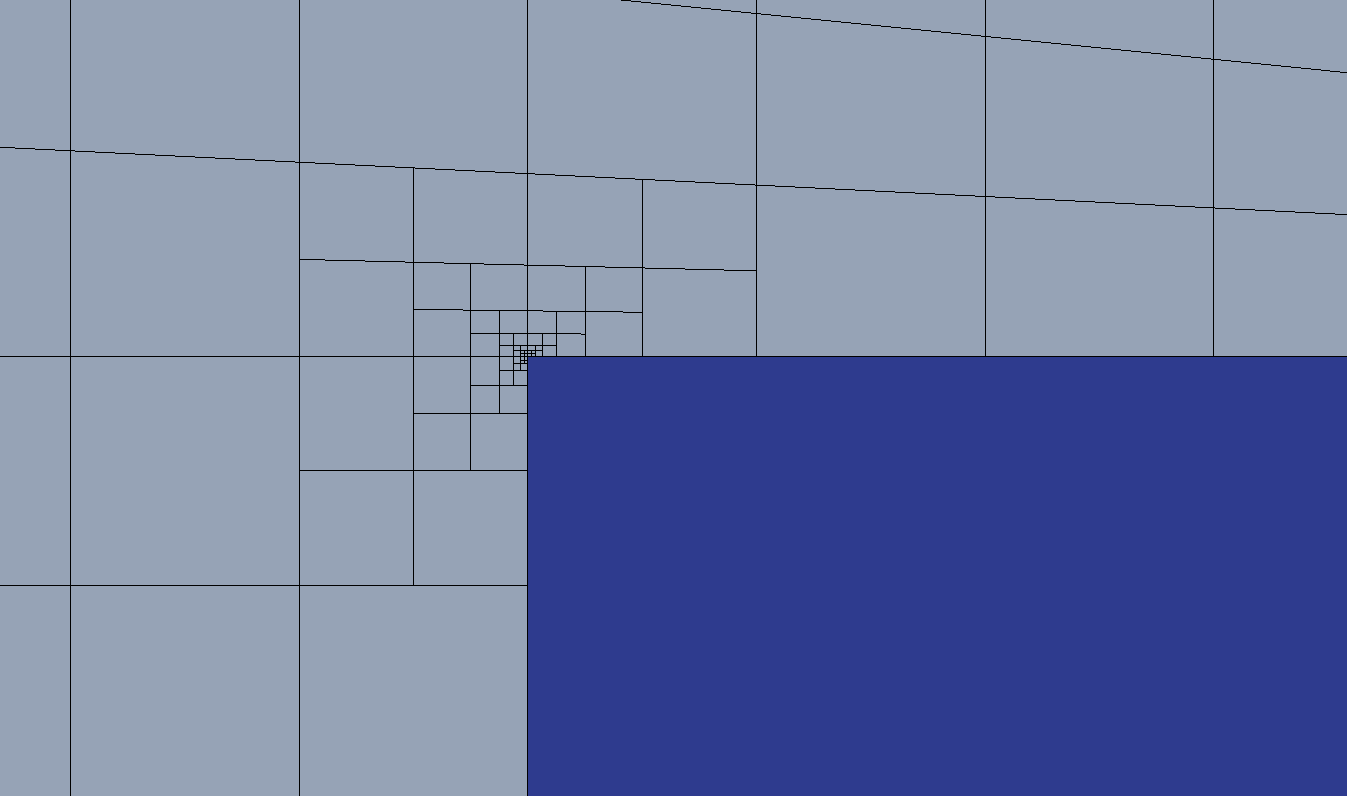} \\
    \includegraphics[width=0.45\textwidth]{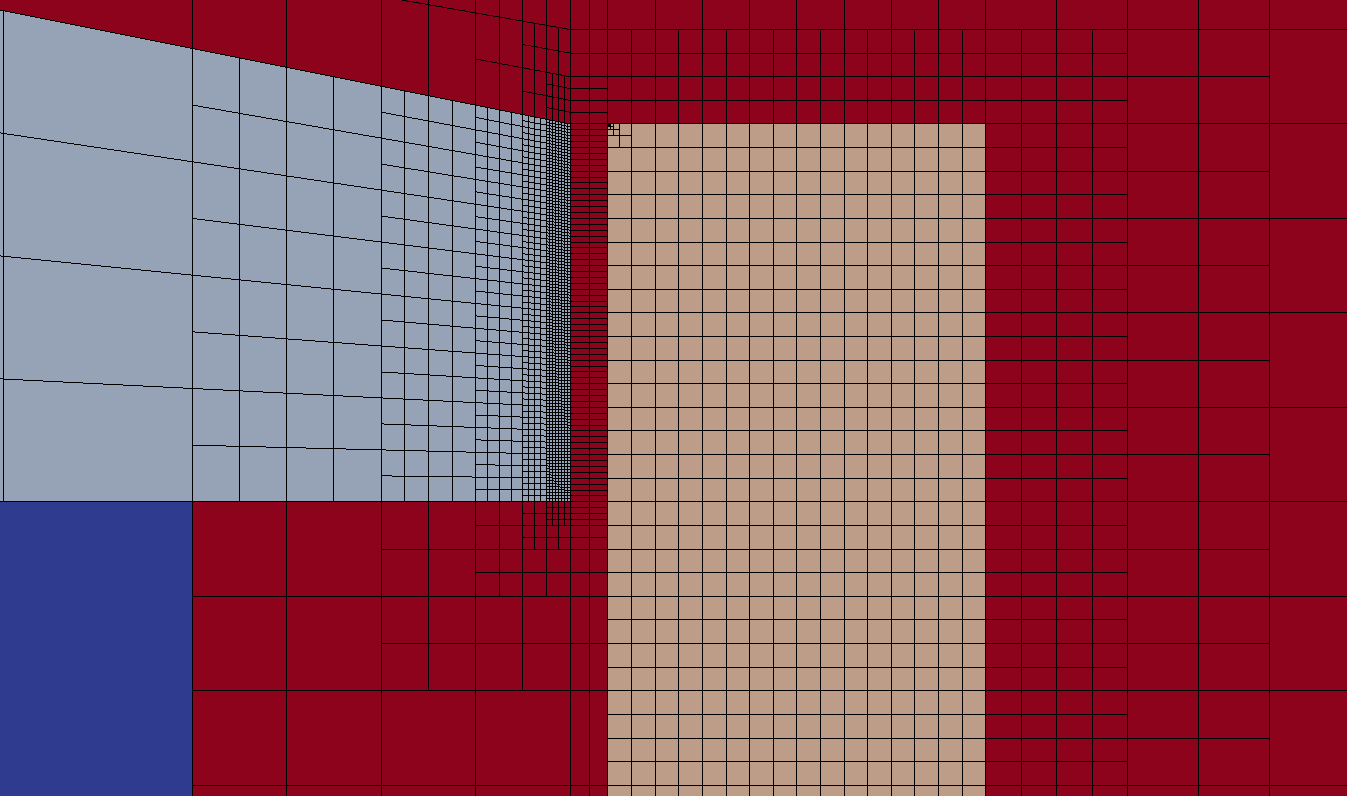} &
    \includegraphics[width=0.45\textwidth]{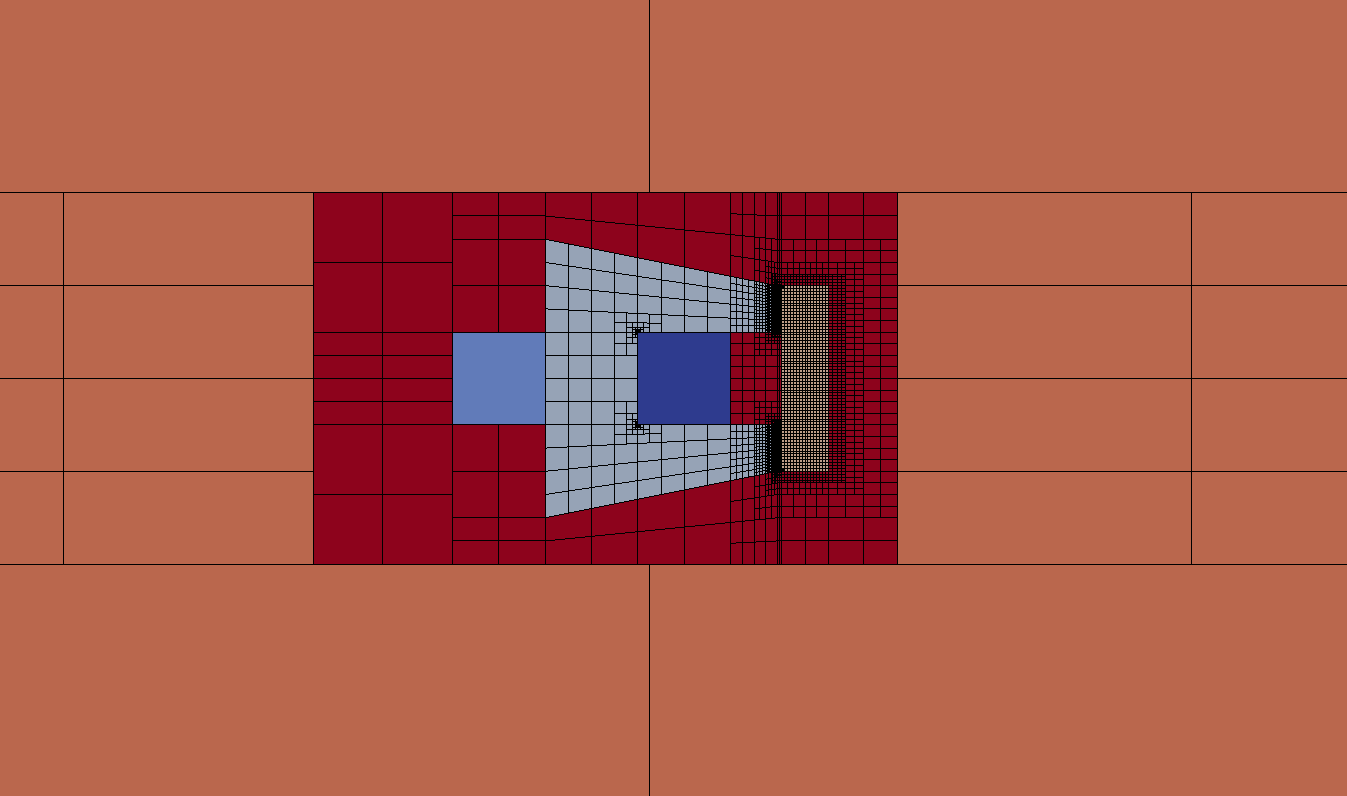} \\
    \end{tabular}
    \caption{C-core actuator: the whole mesh}
    \label{fig:meshCCoreDet}
\end{figure}

Notice that such mesh is refined in the most important regions:
at the corners of the C-shaped core 
where there are singularities of the magnetic permeability;
on the boundary of the plunger 
where there are singularities of the magnetic field and 
inside the air-gap between the C-core and the plunger, see the details in Figure~\ref{fig:meshCCoreDet}.
Then, other regions where there is no need to have an accurate representation of the solution 
are not refined and they are charaterised by the presence of a lot of hanging nodes.
For instance the mesh inside the two current-carrying regions forming the winding is a really big square
element with a lot of hanging nodes, see the details in Figure~\ref{fig:meshCCoreDet}.
Moreover the rectangles far away from the plunger are also big elements 
that represents the truncated air domain, 
see the whole domain in Figure~\ref{fig:meshCCore} and the details in Figure~\ref{fig:meshCCoreDet}.

We consider a fixed air-gap of 5 mm,
we use the same VEM mesh for both Kikuchi and potential formulation, 
while we consider a uniform mesh composed of triangles to get a FEM solution via MagNet.

In Table~\ref{tab:comVEMFEMBx} we show the values of $B_x$ and $F_x$ for each method and approximation degree.
On the one hand, if we fix the method and we vary the approximation degree, such values becomes stable 
(the first two digits of both $B_x$ and $F_x$ stay the same).
On the other hand if we fix the approximation degree and we vary the method,
the computed values are close to each other.
This fact is a numerical evidence that the solutions provided by \vpp{} 
are compatible with the one obtained by FEM.

\begin{table}[!htb]
\centering 
\begin{tabular}{crrr}
       &VEM Kik      &VEM Pot     &FEM     \\
degree &$B_x\,$[T]  &$B_x\,$[T]  &$B_x\,$[T]  \\  
\hline
1 &0.1225     &0.1228    &0.1228    \Tstrut
2 &0.1250     &0.1250    &0.1227    \Bstrut
3 &0.1250     &0.1250    &0.1227    \Bstrut
\hline     
&&&\\
       &VEM Kik      &VEM Pot     &FEM     \\
degree &$F_x\,$[N] &$F_x\,$[N] &$F_x\,$[N] \\  
\hline
1  &-5.7325  &-5.9803  &-5.9898  \Tstrut
2  &-6.0334  &-6.2475  &-6.1091  \Bstrut
3  &-6.0976  &-6.2717  &-6.1254  \Bstrut
\hline     
\end{tabular}
\caption{C-core actuator: $B_x\,$[T] and $F_x\,$[N] computed by FEM and~VEM.}
\label{tab:comVEMFEMBx}
\end{table}

Form another point of view we can also state that 
such results are comparable to a commercial software like MagNet 
that is one of the best software used in the magnetostatic field.

However, we admit that such comparison is not completely clean.
Indeed, the FEM solution is computed on a different mesh, 
we can not obtain exactly the same values.
Moreover, MagNet also makes some post processing on the computed solution 
to get smooth values of both $B_x$ and $F_x$ and, 
consequently, its solution converges to different values with respect to ones of \vpp{}.

Now we proceed with the evaluation of $B_x$ e $F_x$ varying the air-gap width.
We compute such solution \emph{only} via \vpp{} 
since we already validate the method with the previous comparison.
Such data are collected in Table~\ref{tab:convVEM}.
\begin{table}[!htb]
\centering 
\begin{tabular}{cccccc}
\multicolumn{2}{c}{$d$ [mm]}           &0.25   &0.50   &1.00   &2.00       \Bstrut
\hline
\multirow{2}{*}{$B_x$ [T]}       &VEM Kik   &2.0428 &1.1176 &0.5870 &0.3015     \Tstrut
                                 &VEM Pot   &2.0449 &1.1185 &0.5875 &0.3017     \Bstrut
\hline
\multirow{2}{*}{$F_x$ [N]}       &VEM Kik   &-1360.5 &-413.00 &-117.09 &-32.498 \Tstrut
                                 &VEM Pot   &-1365.5 &-416.32 &-119.26 &-33.433 \Bstrut
\hline     
\end{tabular}
\caption{C-core actuator: values of $B_x$ and $F_x$ by varying air-gap widths.}
\label{tab:convVEM}
\end{table}

\begin{figure}[!htb]
\begin{center}
\begin{tabular}{cc}
\multicolumn{2}{c}{$|\B|$}\\
\multicolumn{2}{c}{\includegraphics[width=0.42\textwidth]{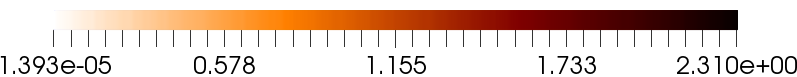}}\\[1em]
\includegraphics[width=0.4\textwidth]{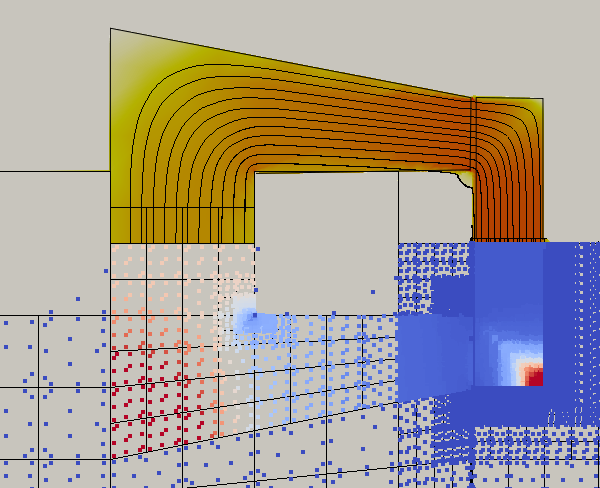} &
\includegraphics[width=0.4\textwidth]{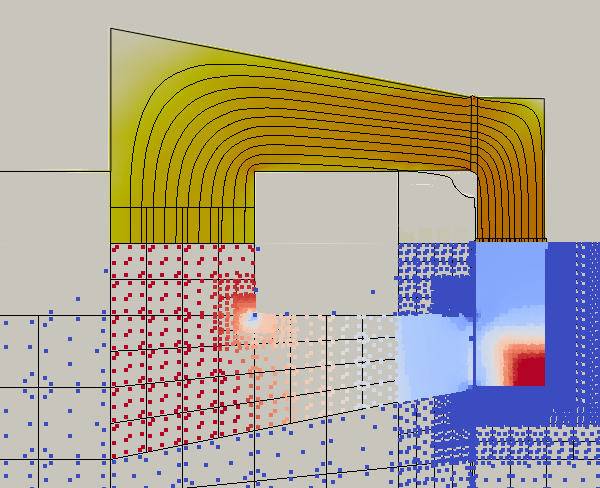} \\
$d=0.25$ mm & $d=0.50$ mm \\[1em]
\includegraphics[width=0.4\textwidth]{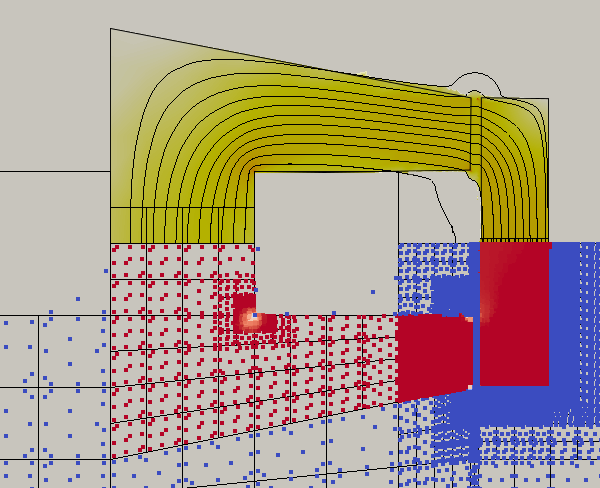} &
\includegraphics[width=0.4\textwidth]{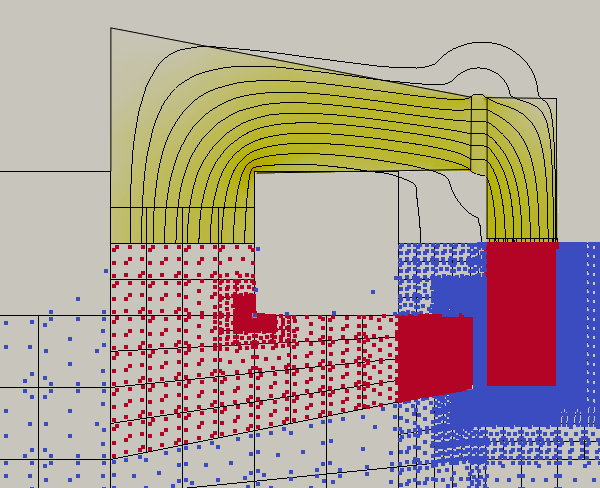} \\
$d=1.00$ mm &$d=2.00$ mm \\
\multicolumn{2}{c}{$\mu_r$}\\
\multicolumn{2}{c}{\includegraphics[width=0.42\textwidth]{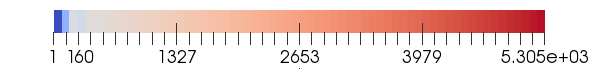}}\\
\end{tabular}
\end{center}
\caption{C-core actuator: the vector field $\B$, on top of each figure,
and the magnetic permeability $\mu_r$ at the quadrature points, on bottom of each figure.}
\label{fig:muAndB}
\end{figure}

Both formulations give similar results by varying the air-gap width, $d$.
More specifically, we observe that the trend of $B_x$ and $F_x$ is approximately $d^{-1}$ and $d^{-2}$, respectively, that is perfectly aligned with the theory and the real-life experiments.
Finally, in Figure~\ref{fig:muAndB} we show both the vector field $\B$ 
at mesh vertices and the value of $\mu_r$ at the quadrature points.

Such data are computed via the VEM potential formulation (Kik VEM is similar so we do not show it).
The behaviour of magnetic field and the magnetic permeability are the expected ones.
Indeed, the vectors are properly aligned and the strength of $\B$ increases for small air-gap widths.
Moreover, both the C-core and the plunger have a uniform $\mu_r$ for large air-gap widths, $d=1.00$ mm and $d=2.00$ mm,
while the material starts to saturate when $d$ is small.

\paragraph{Interior permanent magnet motor}

\begin{figure}[!htb]
\centering
\begin{tabular}{rlcrl}
\multicolumn{5}{c}{\includegraphics[width=0.32\textwidth]{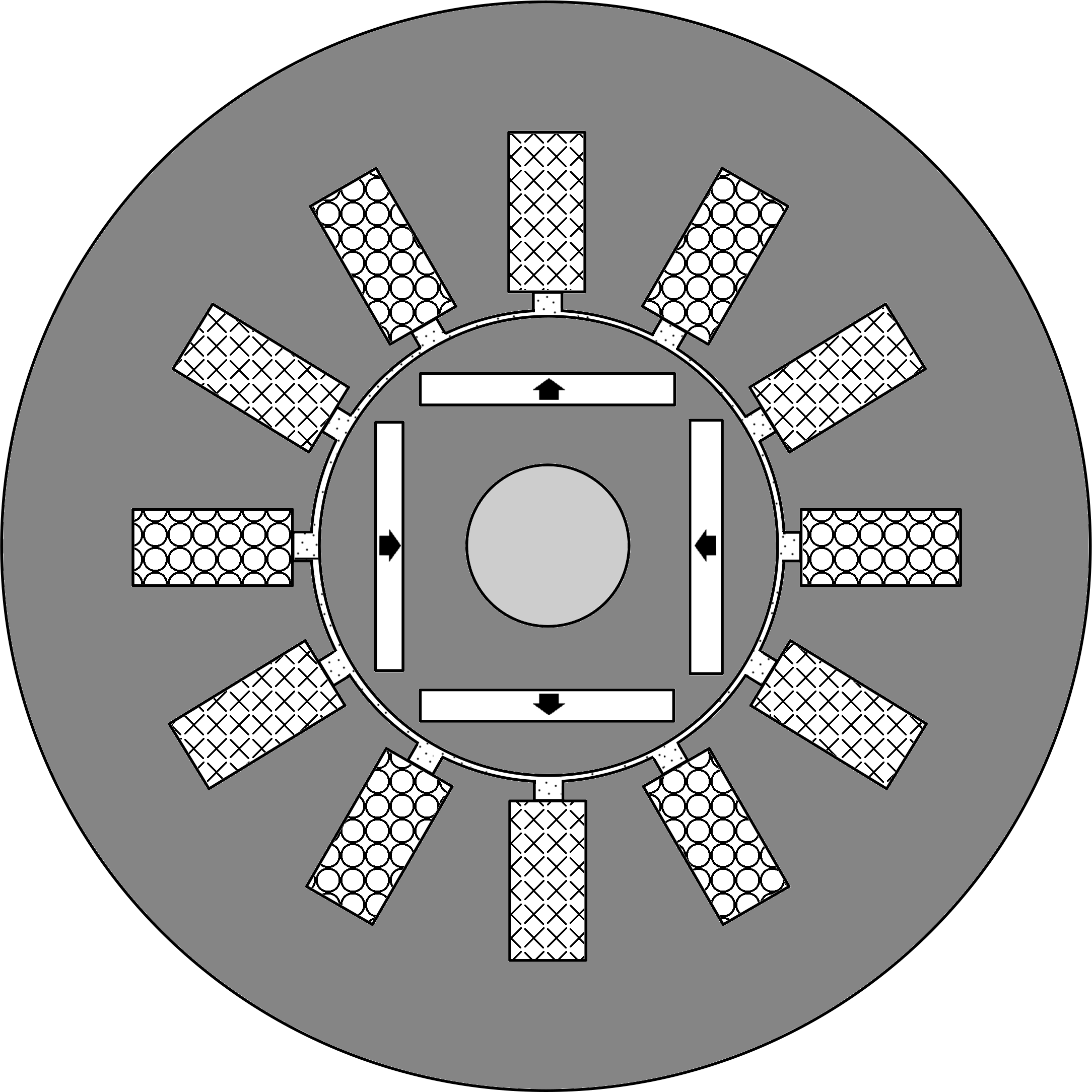}}\\[0.5em]
\includegraphics[width=0.02\textwidth]{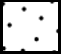}\phantom{,}:&\hspace{-1em}air &\phantom{mm}&
\includegraphics[width=0.02\textwidth]{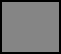}\phantom{,}:&\hspace{-1em}iron \\
\includegraphics[width=0.02\textwidth]{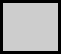}\phantom{,}:&\hspace{-1em}non-magnet steel &\phantom{mm}& 
\includegraphics[width=0.02\textwidth]{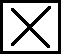}\phantom{,}
\includegraphics[width=0.02\textwidth]{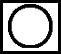}\phantom{,}:&\hspace{-1em}3-phase winding 
\end{tabular}
\caption{Interior permanent magnet motor: the geometry of the 4-pole motor taken into account.}
\label{fig:motorGeo}
\end{figure}

In this section we are using \vpp{} to make a numerical simulation of an Interior-Permanent-Magnet (IPM) motor characterized by four poles and 12 stator slots~see Figure~\ref{fig:motorGeo}.
The external and the rotor diameters are $68.~\text{mm}$ and $30.~\text{mm}$,
respectively, 
while the air-gap width is $0.5~\text{mm}$.
The permanent magnet exhibits a radial magnetization with remanent field equal to $\B_0=1.~\text{T}$ and 
coercive field equal to $\HH_c=7.957\,10^5~\text{Am}^{-1}$.

In this example we are interested in the computation of the cogging torque, i.e.,
the torque acting on the rotor 
when the three-phase current in the rotor slots is zero (no-load operation).
Such quantity is crucial to design a permanent magnet motor since 
it takes into account the tendency of the permanent magnet axis to align with the direction 
that corresponds to the minimum energy stored in the motor. 

\begin{figure}[!htb]
\centering
\begin{tabular}{cc}
\includegraphics[width=0.45\textwidth]{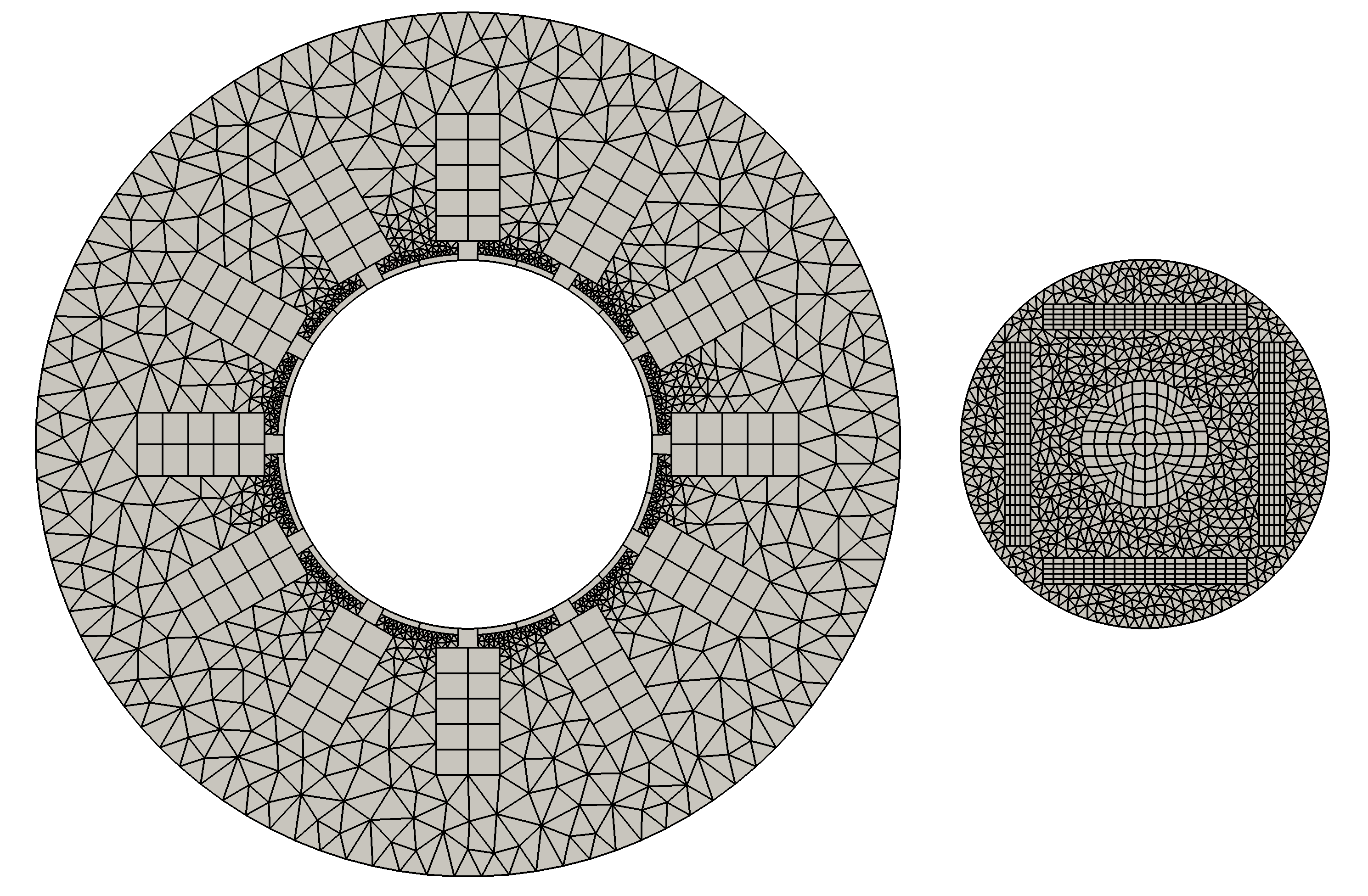} & 
\includegraphics[width=0.45\textwidth]{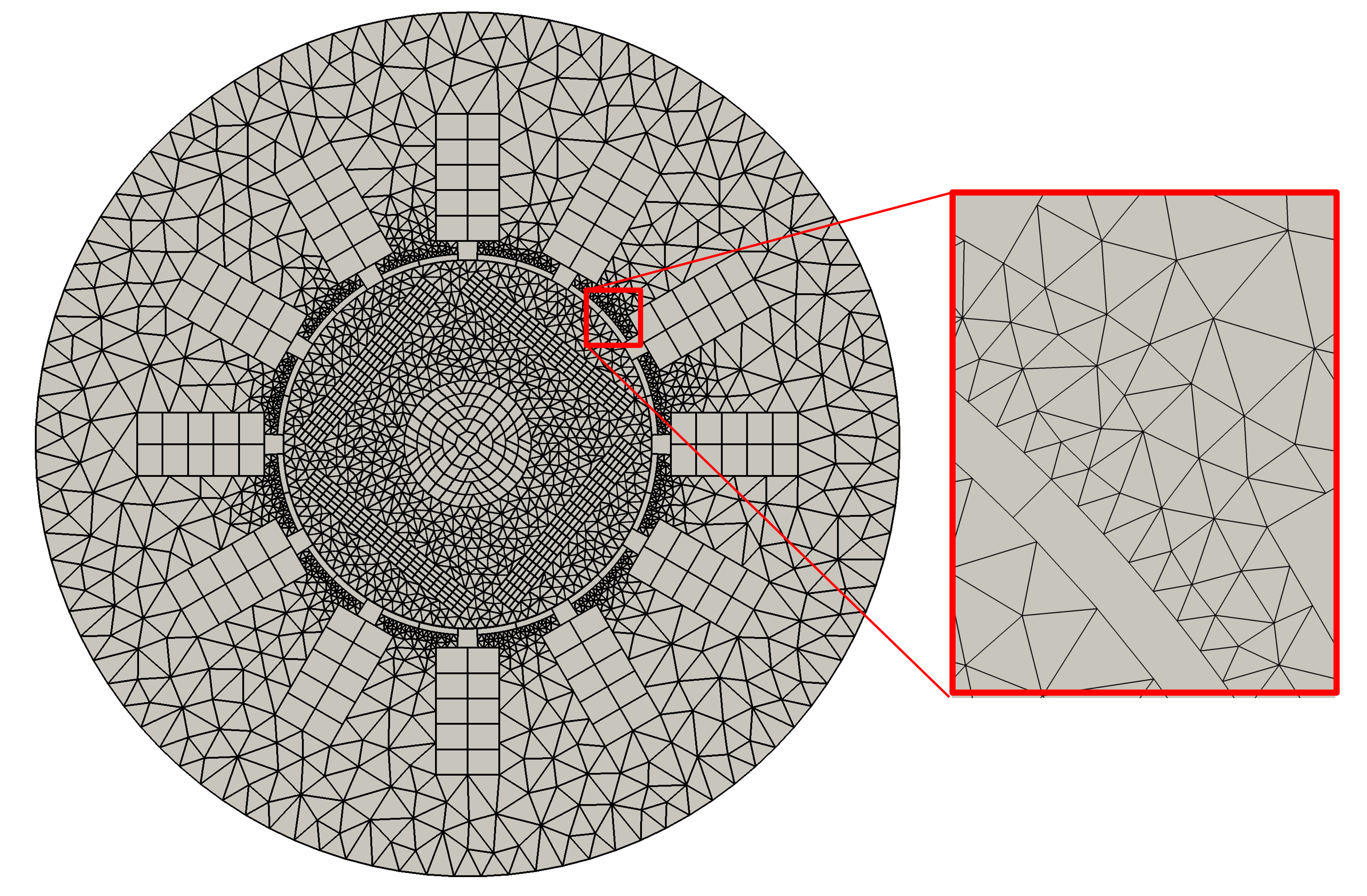} \\
(a) &(b)
\end{tabular}
\caption{Interior permanent magnet motor: (a) stator and rotor mesh, (b) mesh glued with a detail of the hanging-nodes generated by the gluing procedure.}
\label{fig:meshStep}
\end{figure}

The cogging torque has to be computed considering different angular positions of the rotor.
To generate \emph{all} these meshes, 
we can exploit the flexibility of VEM in gluing meshes.
More specifically, we generate stator and rotor mesh separately only one time.
Then, to compute the torque cogging torque with different angles,
we glue them together with the desired angle to run the computation, 
see Figure~\ref{fig:meshStep}.
To achieve this goal in \vpp{} it is implemented the class \texttt{mesh2dMerger}.
Given two meshes $\Omega_{h,1}$ and $\Omega_{h,2}$ in $\R^2$,
such class is able to understand if they have a common boundary
and merge them adding hanging nodes, 
see the detail highlighted in Figure~\ref{fig:meshStep} (b).
Moreover, as you can see from this example,
\texttt{mesh2dMerger} is able to merge two meshes along a curved common boundary
and this process is geometry free. 
Indeed, there is no need of a geometrical description of the curve,
\texttt{mesh2dMerger} recognise piece-wise segments on curved boundaries and 
refine them so stat the two meshes are joined together,
see the detail highlighted in Figure~\ref{fig:meshStep} (b).

\begin{figure}[!htb]
\centering
\begin{tabular}{cc}
$|\B|~[\text{T}]$ &\\
\includegraphics[width=0.32\textwidth]{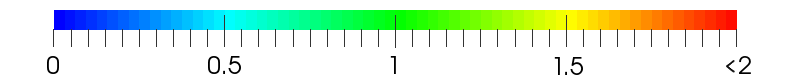} &\multirow{2}{*}{\includegraphics[width=0.55\textwidth]{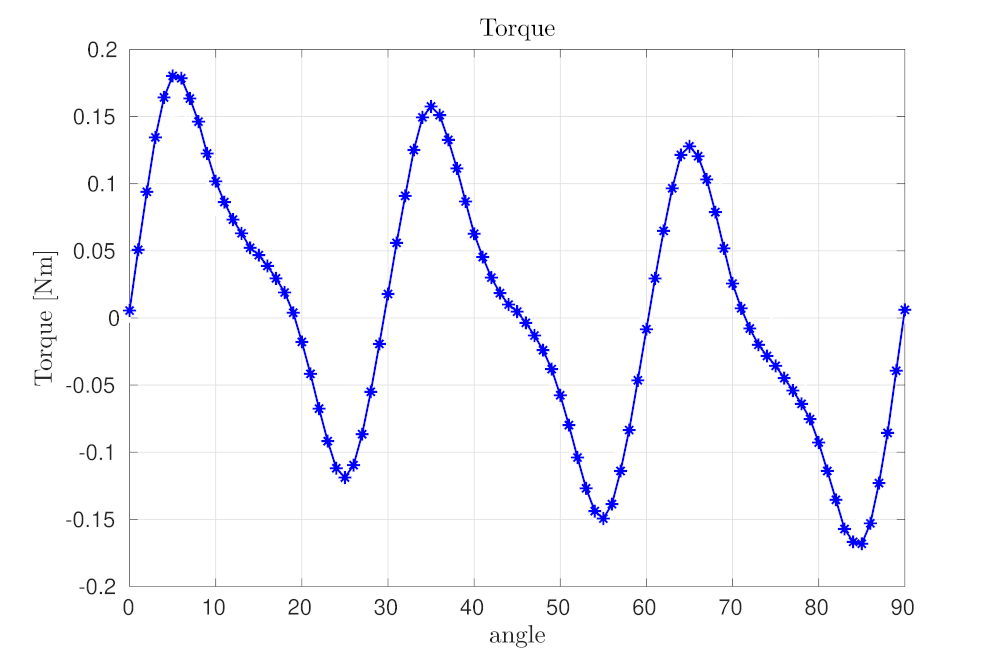}} \\[0.2em]
\includegraphics[width=0.32\textwidth]{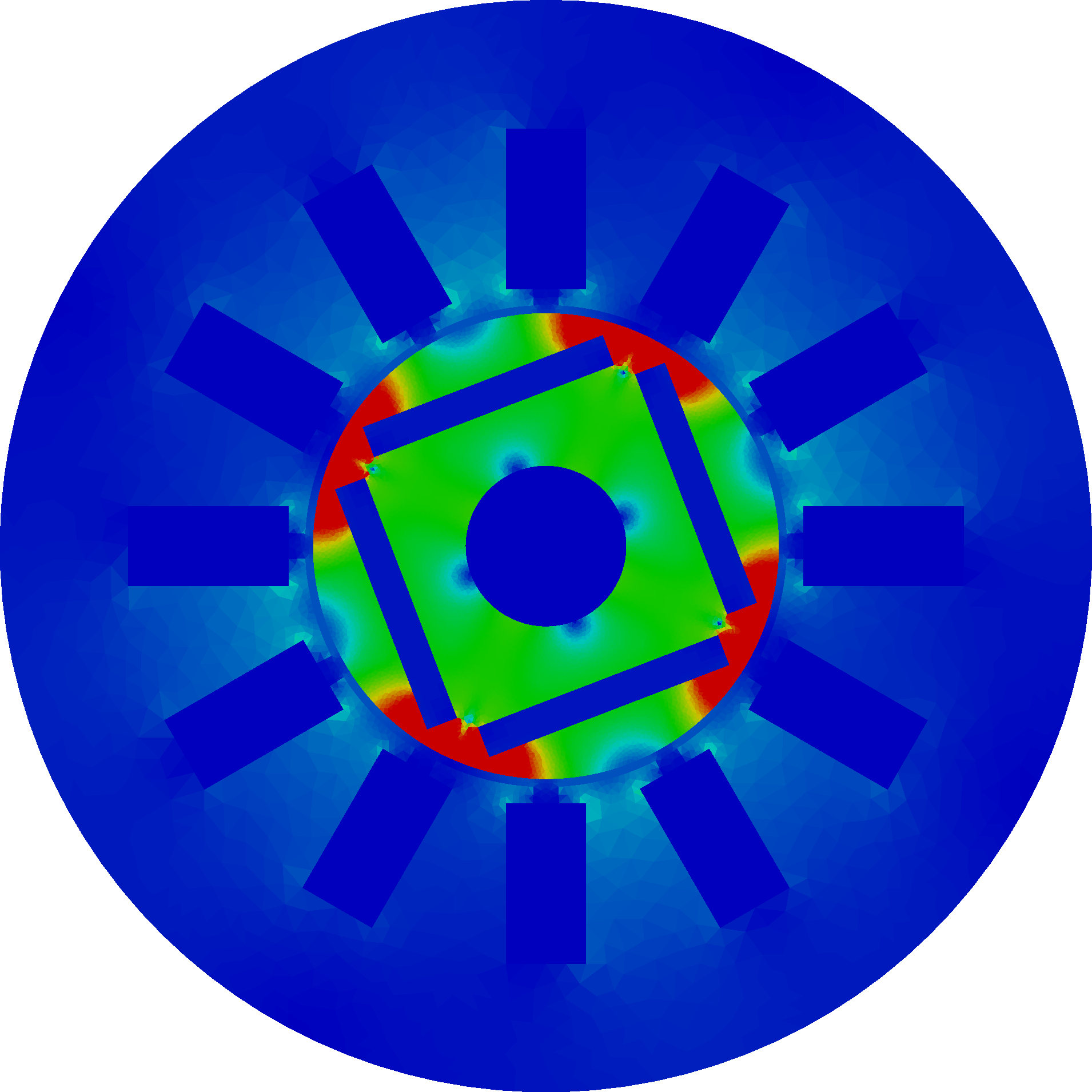}&\\
(a) & (b)\\
\end{tabular}
\caption{Interior permanent magnet motor: (a) magnetic induction map with an angle of $21^\circ$ (b) values of torque for different position of the rotor.}
\label{fig:ese3Torque}
\end{figure}

To compute the cogging torque, 
we use the Maxwell stress tensor approach considering a cylindrical surface co-axially 
located with respect to the rotation axis as the integration surface accordingly.
In Figure~\ref{fig:ese3Torque} (a) we depict the magnitude of the magnetic field $\B$ when 
the torque-angle is $21^\circ$ degrees while
in Figure~\ref{fig:ese3Torque} (b) we show the torque-angle curve for step equal to $1^\circ$.
As expected, the torque period is equal to $30^\circ$, indeed
$$
360^\circ/\text{LCM}(4,12) = 30^\circ\,,
$$  
where $\text{LCM}$ is the least common multiple operator and 
$4$ and $12$ are the number of permanent magnets and slots, respectively.
Moreover, it exhibits zero mean value over the period. 
Once again there is a good agreement between Kikuchi and potential virtual element formulations.

\paragraph{Optimal shape design of an electromagnet}

In this section \vpp{} was used for handling the shape variation of a magnet.
More specifically, we consider the optimal design of a magnet for applications in clinical hypethermia.

A typical device is characterised by a magnetic core made of ferrite
and it exhibits three limbs. 
Two series-connected current-carrying coils are wound on the central limb 
which has a wide air-gap, 
where the patient is accommodated during the treatment.

\begin{figure}[htb]
\centering
     \begin{tabular}{rlcrl}
     \multicolumn{5}{c}{
     \begin{overpic}[abs,unit=1mm,scale=.25]{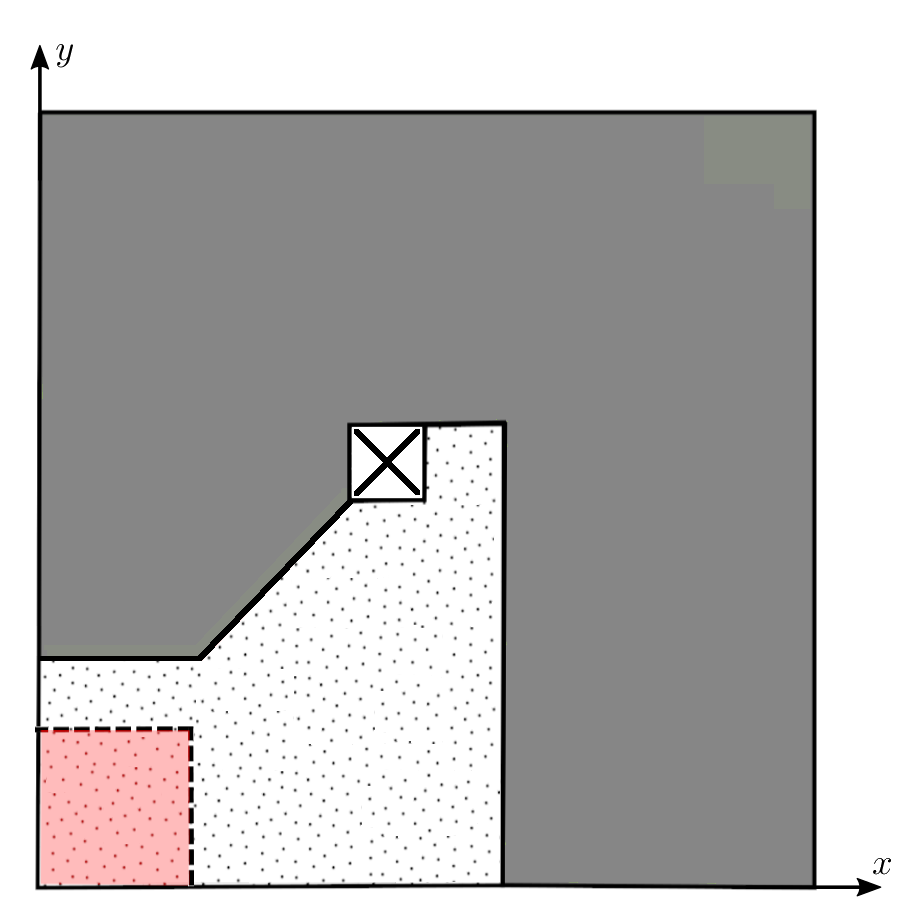}
     \put(-1,-1){{\Large $o$}}
     \put(5,5){{\Large $\Omega_{\mathcal{I}}$}}
     \put(47,47){{\Large $\Omega_{\mathcal{F}}$}}
     \put(27,24){{\Large $\Omega_{\mathcal{J}}$}}
     \put(19,19){{\Large $\Gamma_{\mathcal{P}}$}}
     \put(27,6){{\Large $\Omega_{\mathcal{A}}$}}
     \put(-2, 58){{\Large $l_2$}}
     \put(61, -4){{\Large $l_1$}}
     \end{overpic}}\\[1.0em]
     \includegraphics[width=0.02\textwidth]{ironIcon.png}\phantom{,}:&\hspace{-1em}iron core, $\Omega_\mathcal{F}$
     &&\includegraphics[width=0.02\textwidth]{current1Icon.png}\phantom{,}:&\hspace{-1.0em}winding, $\Omega_\mathcal{J}$\\
     \includegraphics[width=0.02\textwidth]{airIcon.png}\phantom{,}:&\hspace{-1em}air, $\Omega_\mathcal{A}$ &\phantom{m}&\includegraphics[width=0.02\textwidth]{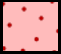}\phantom{,}:&\hspace{-1em}domain of interest, $\Omega_\mathcal{I}\subset\Omega_\mathcal{A}$\\
     \end{tabular}
\caption{Optimal shape design of an electromagnet: domain to consider.}
\label{fig:optDom}
\end{figure}

A quarter of the model geometry here considered is shown in Figure~\ref{fig:optDom}.
A ferrite core fills in region $\Omega_\mathcal{F}$, 
while an air-gap 30 cm high and 20 cm long incorporates the region of interest $\Omega_{\mathcal{I}}$, i.e.,
a path along which the degree of uniformity of flux density is controlled.
The non-linear $\B$-$\HH$ curve of the ferrite considered in the model exhibits an initial value of relative permeability equal to 1800,
and a saturation flux density of 490 mT. 
The complementary domain includes the winding cross-section $\Omega_\mathcal{J}$
which is composed of 16 turns and carries a sinusoidal current of 150 Arms at 100 kHz 
and an air region, $\Omega_\mathcal{A}$.

The design challenge is to shape the magnetic pole, $\Gamma_{\mathcal{P}}$, in such a way 
that the prescribed field takes place in the region of interest, $\Omega_{\mathcal{I}}$.
Such problem is usually reformulated as an inverse problem solved by means of a numerical optimisation technique. 
Indeed, a suitable functional is defined to measure the discrepancy between actual and
the prescribed flux density in the region of interest,~$\Omega_{\mathcal{I}}$.

\begin{figure}[!htb]
\centering
\begin{subfigure}[t]{0.30\textwidth}
\centering
\includegraphics[width=0.99\textwidth]{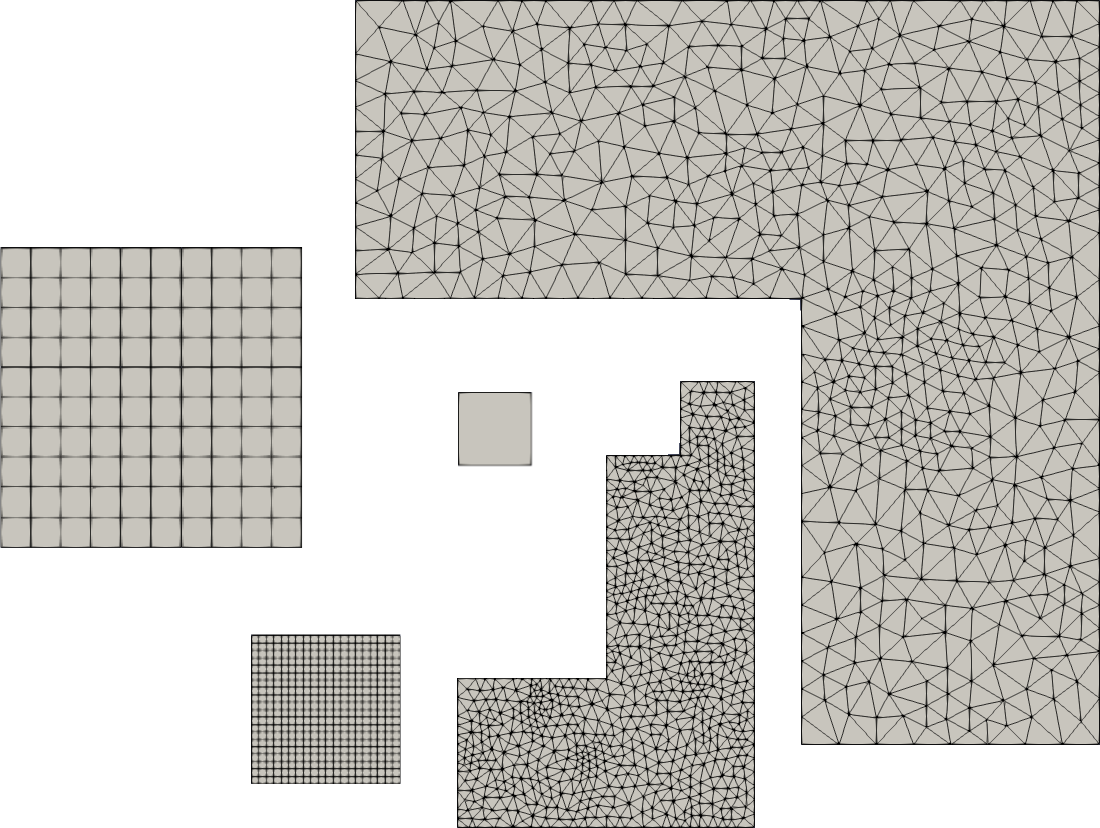}
\caption{Initial meshes.}
\label{fig:wholeProc1}
\end{subfigure}\hfill
\begin{subfigure}[t]{0.30\textwidth}
\centering
\includegraphics[width=0.99\textwidth]{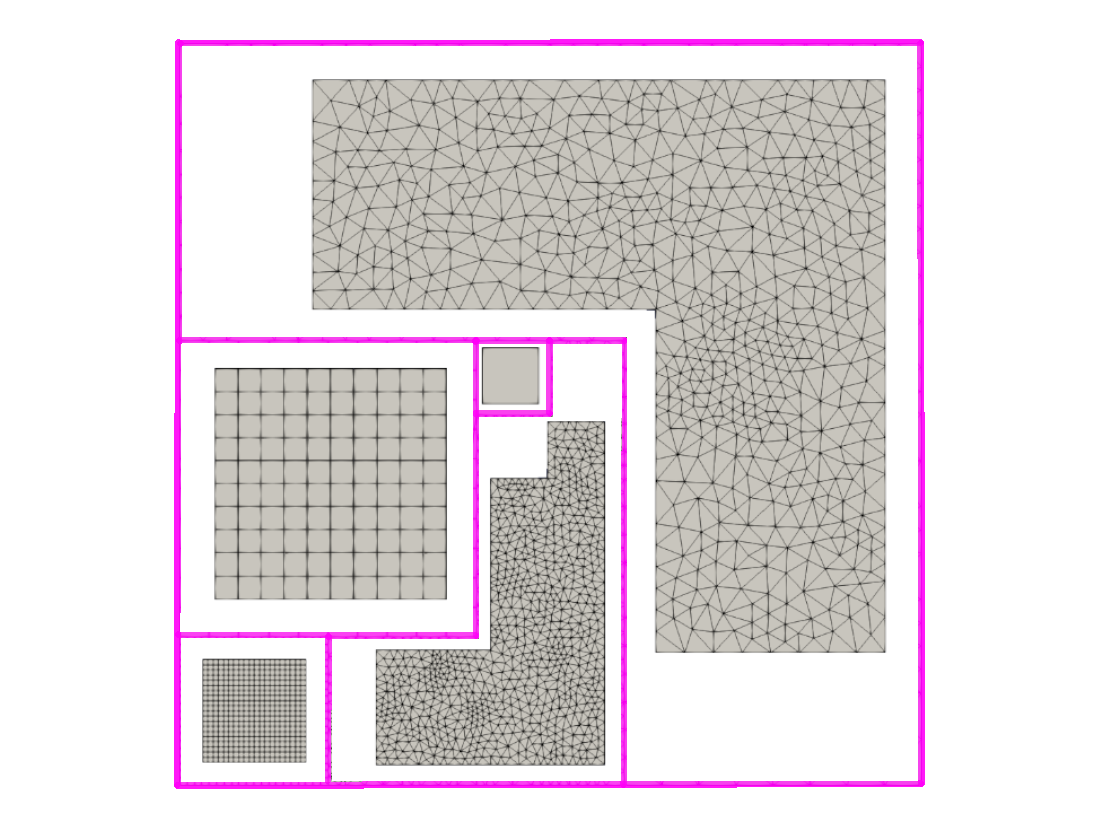}
\caption{Gluing guidelines.}
\label{fig:wholeProc1.5}
\end{subfigure}\hfill
\begin{subfigure}[t]{0.30\textwidth}
\centering
\includegraphics[width=0.99\textwidth]{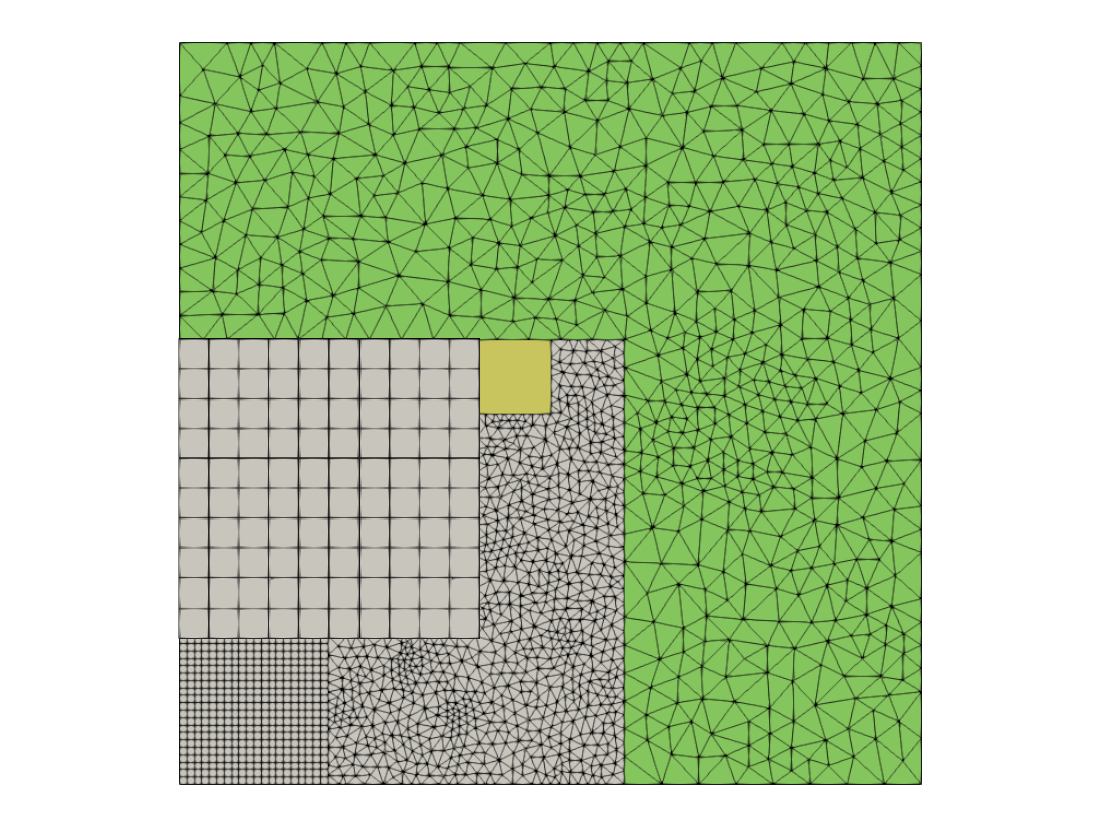}
\caption{The resulting mesh glued.}
\label{fig:wholeProc2}
\end{subfigure}\\[1em]
\begin{subfigure}[t]{0.30\textwidth}
\centering
\includegraphics[width=0.99\textwidth]{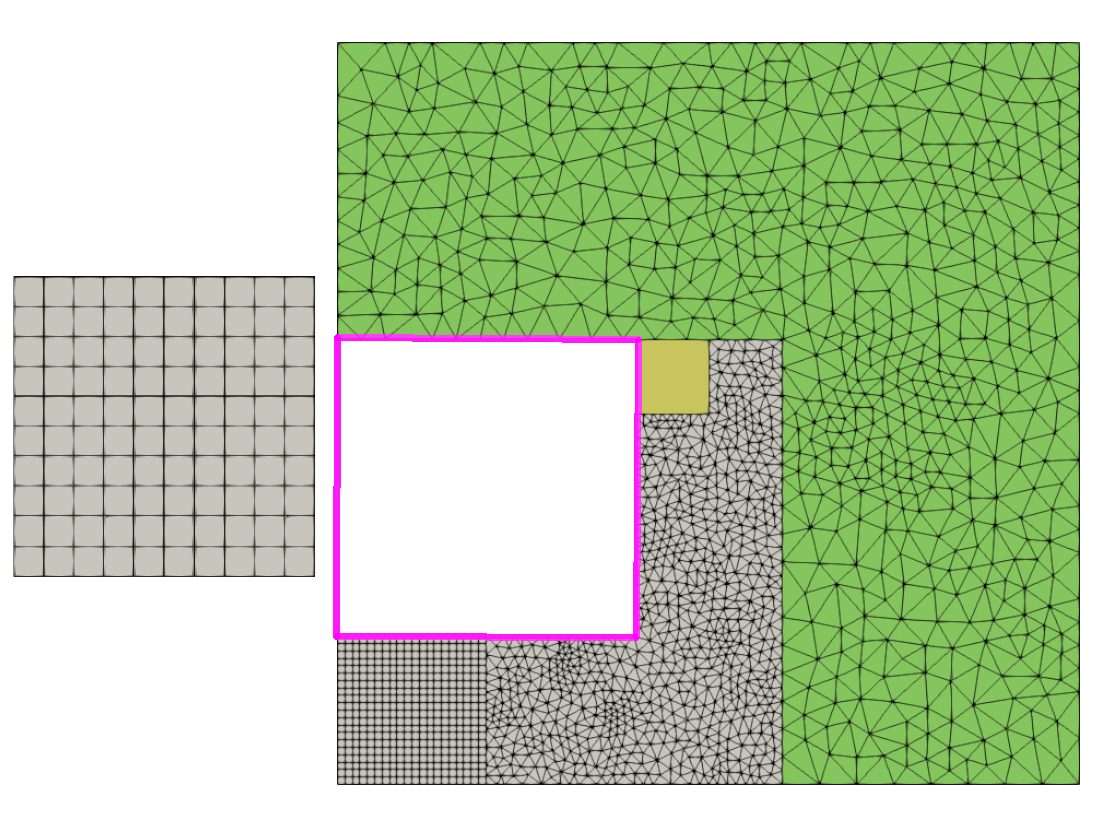}
\caption{Extraction of the domain where to draw the magnet profile.}
\label{fig:wholeProc3}
\end{subfigure}\hfill\hfill
\begin{subfigure}[t]{0.30\textwidth}
\centering
\includegraphics[width=0.99\textwidth]{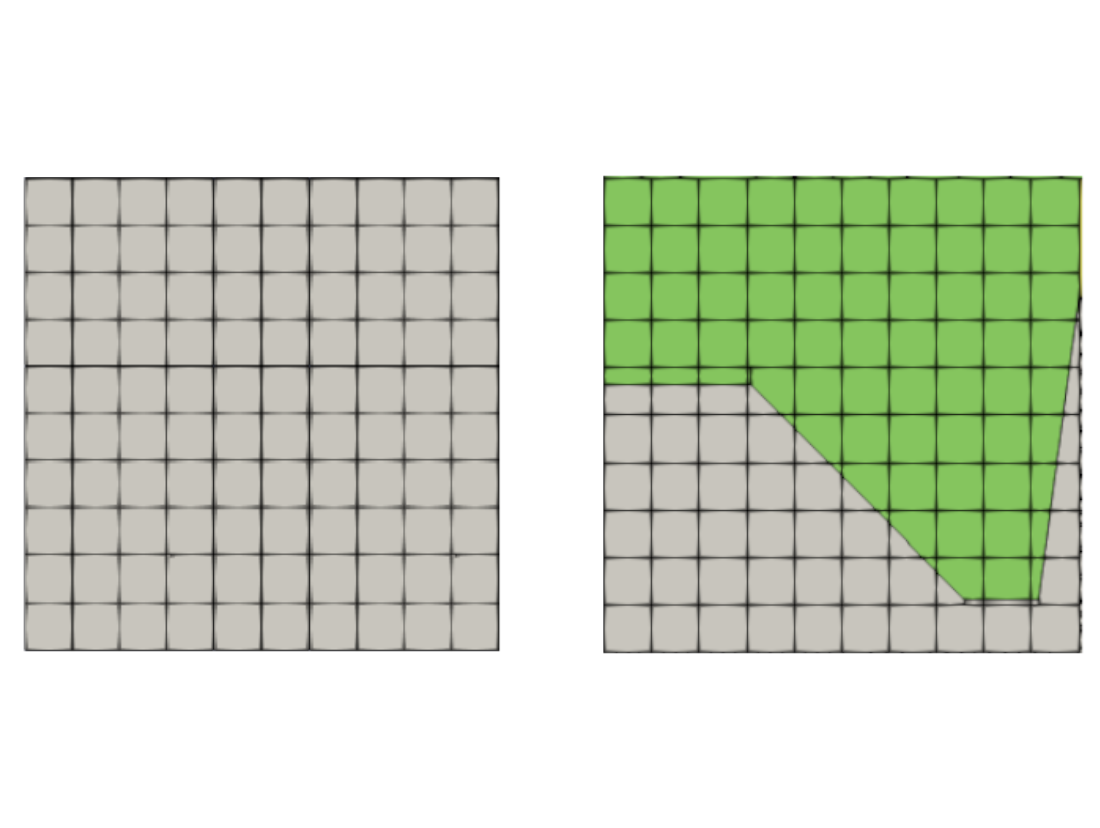}
\caption{Draw the computed magnet profile on this small mesh and cut elements.}
\label{fig:wholeProc4}
\end{subfigure}\hfill
\begin{subfigure}[t]{0.30\textwidth}
\centering
\includegraphics[width=0.99\textwidth]{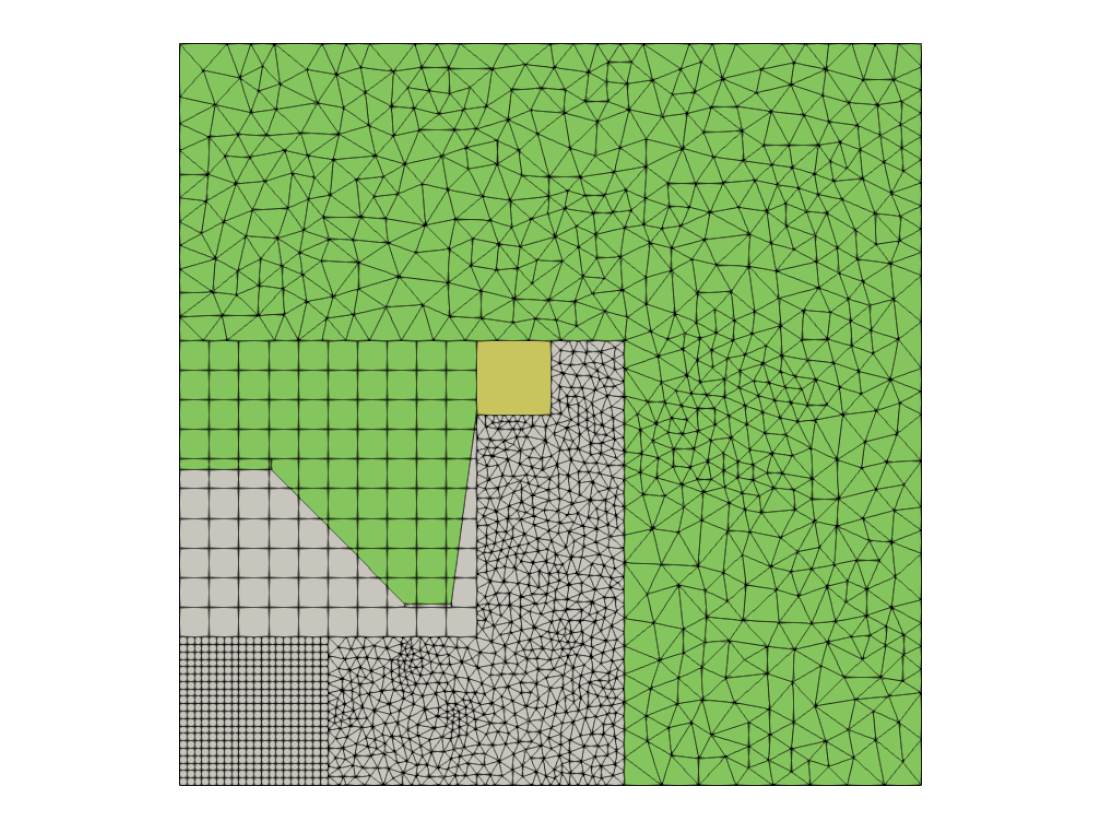}
\caption{The mesh is glued back.}
\label{fig:wholeProc5}
\end{subfigure}
\caption{Optimal shape design of an electromagnet: the whole meshing procedure during optimisation.}
\label{fig:wholeProc}
\end{figure}

The proper optimisation procedure is deeply described in~\cite{freeCutting}.
However, in this section we show how to exploit the flexibility in mesh generation of VEM to tackle such procedure. 
In Figure~\ref{fig:wholeProc} we summarise the whole procedure
we use to generate the mesh.
Each piece of the domain is discretised at the beginning and then glued together using \texttt{mesh2dMerger},
see Figure~\ref{fig:wholeProc} (a), (b) and~(c). 
Then, the mesh associated with region where we have to draw the magnet profile 
is extracted and the structured quadrilateral mesh is stored, see Figure~\ref{fig:wholeProc} (d).
In Figure~\ref{fig:wholeProc} (e) we show the key step where VEM plays an important role.
We start from a uniform mesh composed, e.g. of quadrilateral elements regularly spaced, 
see Figure~\ref{fig:wholeProc} (e) (left),
and then introduce suitable cuts by segments, see Figure~\ref{fig:wholeProc} (e) (right).
Consequently the original quadrilateral element is cut in a pair of polygons that are naturally handled by VEM and 
there is no need of sub-triangulate or re-build the mesh.
To cut a mesh via a affine segments,
in \vpp{} there is a specific class called \texttt{addStraightLineOnMesh2d}.
Finally, this piece of mesh is glued back to the mesh, Figure~\ref{fig:wholeProc} (f).

Thus, the optimisation procedure looks for the best profile by changing the locations of the points that define it.
However, at each optimisation step the mesh is not built from scratch:
the meshing part jumps from Figure~\ref{fig:wholeProc} (e) and~(f).
As a consequence, the meshing part is quite fast since 
it has to make a cut of a small part of the domain without touching the other parts
and the background mesh is structured.

\begin{figure}[!htb]
    \centering
    \begin{tabular}{cc}
    \includegraphics[width=0.35\textwidth]{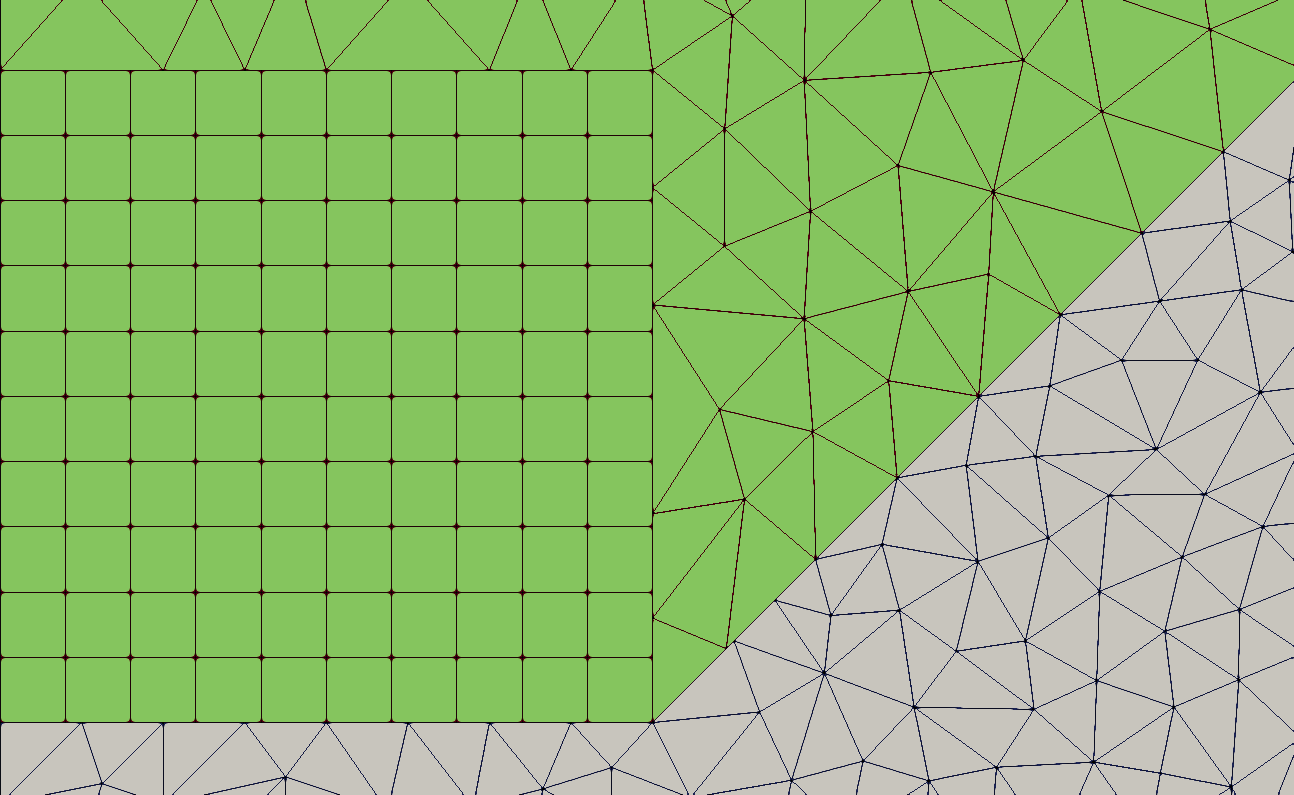} &
    \includegraphics[width=0.35\textwidth]{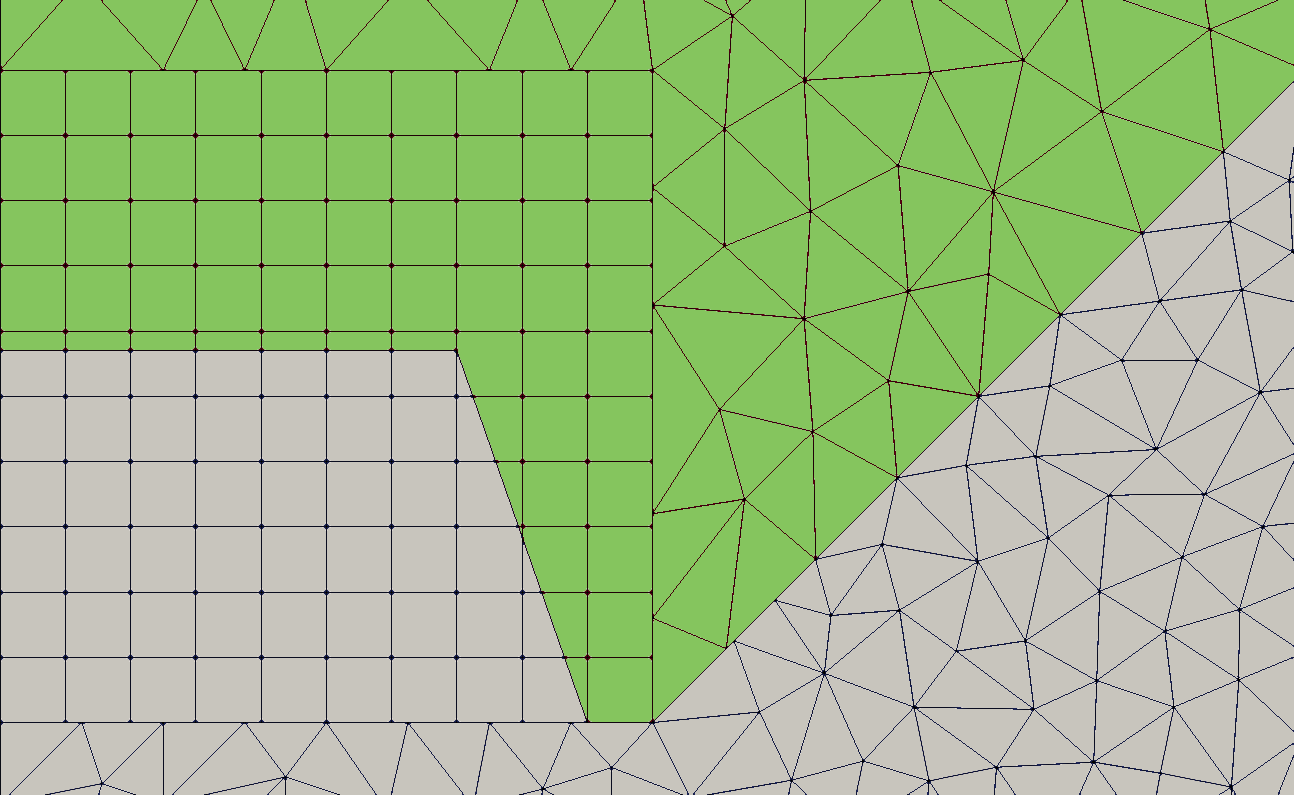}   \\
    (a) & (b)
    \end{tabular}
    \caption{Optimal shape design of an electromagnet: the initial profile (a) and the optimised profile (b).}
    \label{fig:opt1}
\end{figure}

\begin{figure}[!htb]
    \centering
    \begin{tabular}{cc}
    \multicolumn{2}{c}{$|\B|$ [T]}\\
    \multicolumn{2}{c}{\includegraphics[width=0.25\textwidth]{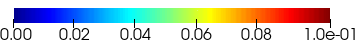}}\\[0.5em]
    \includegraphics[width=0.4\textwidth]{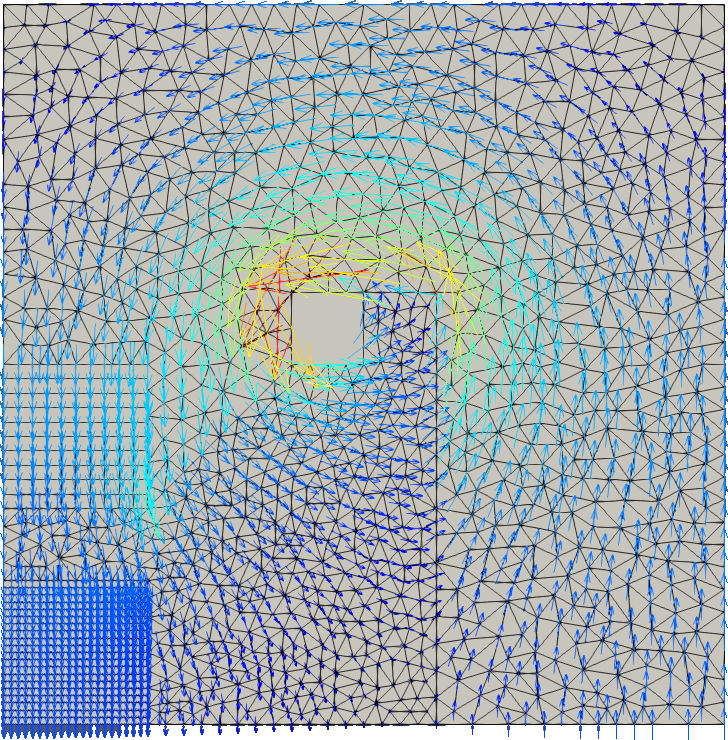} &
    \includegraphics[width=0.4\textwidth]{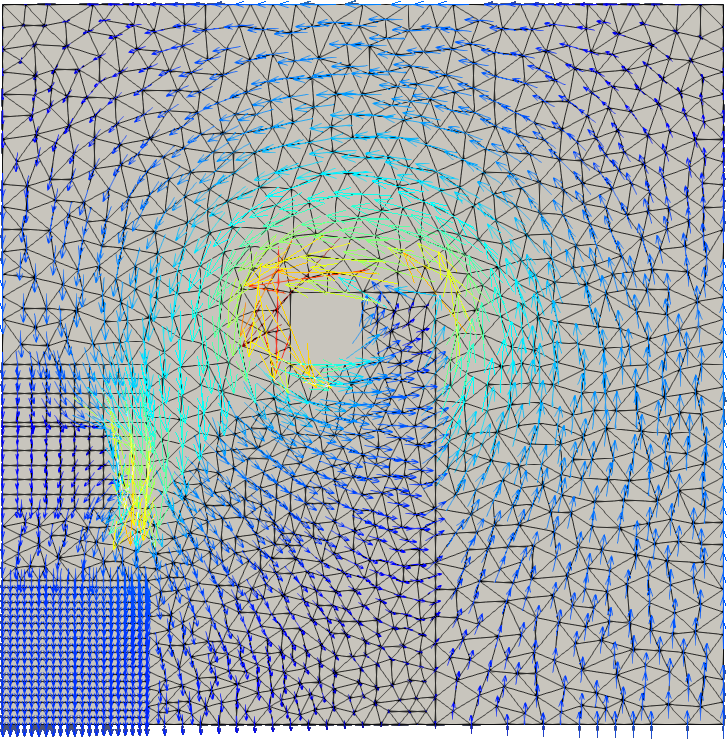}   \\
    (a) &    (b) \\
    \end{tabular}
    \caption{Optimal shape design: field $\B$ arrow maps for the designs (a) and (b).}
    \label{fig:profile}
\end{figure}

In Figure~\ref{fig:opt1} (a) we show the initial $\B$ computed on the profile depicted in Figure~\ref{fig:profile} (a),
while in Figure~\ref{fig:opt1} (b) the resulting $\B$ field on an optimised geometry that is represented in Figure~\ref{fig:profile} (b).

\begin{figure}
    \begin{center}
    \includegraphics[width=0.55\textwidth]{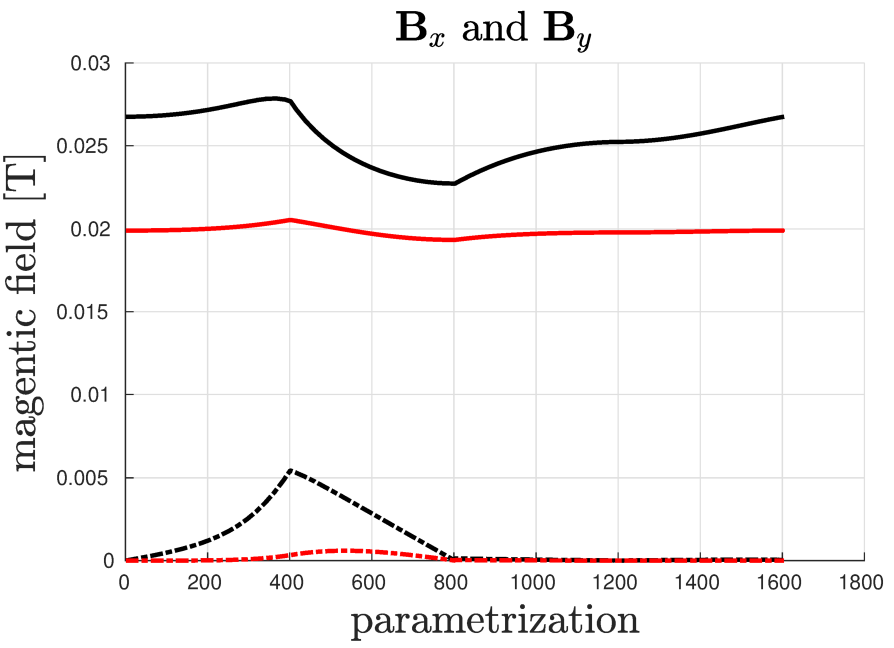}
    \end{center}
    \caption{Optimal shape design: comparison between the initial and final values of $\B$ in black the values associated with the configuration (a), in red the ones associated with (b).}
    \label{fig:linesObBy}
\end{figure}

To have a more qualitative proof of the effectiveness of the optimisation procedure, 
in Figure~\ref{fig:linesObBy} we plot both the $x$ and $y$ component of the magnetic filed $\B$ along
the boundary of the domain of interests, $\Omega_{\mathcal{I}}$.
If we verify that such components are almost constant on this boundary,
the vector field $\B$ is also constant inside $\Omega_{\mathcal{I}}$.
This fact is due to the presence of homogeneous Neumann and Dirichlet boundary conditions on $y=0$ and $x=0$, respectively, and 
the fact that there is no current source inside $\Omega_{\mathcal{I}}$.
The optimised profile does create a uniform magnetic field inside the region of interests.
Indeed the red lines that are associated with the $x$ and $y$ components of the field $\B$ are almost flat.

\subsubsection*{Acknowledgement}
The development of this \cpp{} library is not a work made by only one person. 
The author would like to thank all the other professors and researchers 
that give a contribution on such a code. 
In particular he would like to thank L.~{Beir{\~a}o da Veiga}, L.~Mascotto, A.~Russo, G.~Vacca
and M.~Visinoni.

\bibliographystyle{plain}
\bibliography{bibliography}

\end{document}